\pdfoutput=1
\PassOptionsToPackage{capitalise}{cleveref}
\documentclass{lmcs}
\pdfoutput=1

\usepackage{lastpage}
\lmcsdoi{22}{1}{25}
\lmcsheading{}{\pageref{LastPage}}{}{}%
{Mar.~12,~2024}{Mar.~12,~2026}{}


\let\parencite\cite
\usepackage{amsthm, amsmath}
\usepackage{thmtools} 
\usepackage{hyperref}
\usepackage{cleveref}
\usepackage{my-preamble}
\usepackage[utf8]{inputenc}

\newenvironment{example}{\begin{exa}}{\end{exa}}
\Crefname{exa}{Example}{Examples}
\newenvironment{remark}{\begin{rem}}{\end{rem}}
\Crefname{rem}{Remark}{Remarks}
\newenvironment{proposition}{\begin{prop}}{\end{prop}}
\Crefname{prop}{Proposition}{Propositions}
\newenvironment{corollary}{\begin{cor}}{\end{cor}}
\Crefname{cor}{Corollary}{Corollaries}
\newenvironment{theorem}{\begin{thm}}{\end{thm}}
\Crefname{thm}{Theorem}{Theorems}
\newenvironment{lemma}{\begin{lem}}{\end{lem}}
\Crefname{lem}{Lemma}{Lemmas}
\newenvironment{definition}{\begin{defi}}{\end{defi}}
\Crefname{defi}{Definition}{Definitions}
\newenvironment{notation}{\begin{nota}}{\end{nota}}
\Crefname{nota}{Notation}{Notations}

\theoremstyle{definition}
\newtheorem{construction}[thm]{Construction}
\Crefname{construction}{Construction}{Constructions}
\newtheorem{terminology}[thm]{Terminology}
\Crefname{terminology}{Terminology}{Terminologies}
\newtheorem{axiom}{Axiom}

\Crefname{axiom}{Axiom}{Axioms}

\theoremstyle{defC}
\newtheorem{factC}[thm]{Fact}
\Crefname{factC}{Fact}{Facts}
\renewenvironment{fact}{\begin{factC}}{\end{factC}}

\newtheoremstyle{consCf}%
  {6pt}
  {6pt}
  {\normalfont}
  {}
  {\bfseries}
  {{\bfseries.}}
  {5pt plus 1pt minus 1pt}
  {\thmname{#1} \thmnumber{#2} \thmnote{\normalfont#3}}
\theoremstyle{consCf}

\newtheorem{consC}[thm]{Construction}
\Crefname{consC}{Construction}{Constructions}
\Crefname{exaC}{Example}{Examples}

\begin{document}

\title[Homotopy type theory as a language for diagrams]{Homotopy type theory as a language for diagrams of $\infty$-logoses}

\titlecomment{This article is an extended version of \parencite{uemura2023diagrams}}

\author[T.~Uemura]{Taichi Uemura\lmcsorcid{https://orcid.org/0000-0003-4930-1384}}
\address{Stockholms universitet}

\begin{abstract}
  We show that certain diagrams of $\infty$-logoses are reconstructed in homotopy type theory extended with some lex, accessible modalities, which enables us to use plain homotopy type theory to reason about not only a single $\infty$-logos but also a diagram of $\infty$-logoses.
  This also provides a higher dimensional version of Sterling's synthetic Tait computability---a type theory for higher dimensional logical relations.
\end{abstract}

\maketitle

\section{Introduction}
\label{sec:introduction}

An \emph{\(\infty\)-logos}, also known as an \emph{\(\infty\)-topos}
\parencite{lurie2009higher,anel2021topo-logie}\footnote{The term
  \(\infty\)-logos is Anel and Joyal's terminology
  \parencite{anel2021topo-logie} for \(\infty\)-topos considered as an
  algebraic structure rather than a geometric object. A morphism of
  \(\infty\)-logoses is always considered in the direction of the inverse
  image functor. We use this terminology to clarify the direction of
  morphisms when speaking about (co)limits of \(\infty\)-logoses.}, is an
\((\infty, 1)\)-category that looks like the \((\infty, 1)\)-category of spaces,
among other aspects of it. An \(\infty\)-logos is a place where one can do
homotopy theory just as an ordinary logos is a place where one can do
set-level mathematics.

\emph{Homotopy type theory} \parencite{hottbook} is another place to
do homotopy theory. It is a type theory in the style of
Martin-L\"{o}f \parencite{martin-lof1975intuitionistic} extended by the
\emph{univalence axiom} and \emph{higher inductive types}. The former
forces types to behave like spaces rather than sets, and the latter
allow us to build types representing spaces such as spheres and tori.

\(\infty\)-logoses are conjectured to admit interpretations of homotopy
type theory so that theorems proved in homotopy type theory can be
translated in an arbitrary \(\infty\)-logos. Although the conjecture has
not yet been fully solved (see, for example, \parencite{shulman2019toposes}
for substantial progress), homotopy type theory has brought insight to
\(\infty\)-logos theory. For example, the proof of the Blakers-Massey
connectivity theorem in homotopy type theory
\parencite{hou2016mechanization} has led to a new generalized
Blakers-Massey theorem that holds in an arbitrary \(\infty\)-logos
\parencite{anel2020generalized}.

An \(\infty\)-logos, however, does not live alone. \(\infty\)-logoses are often
connected by functors which are also connected by natural
transformations. Plain homotopy type theory is, at first sight, not
sufficient to reason about a diagram of \(\infty\)-logoses, because the
actions of the functors and natural transformations are not
internalized to type theory. Even worse, it is impossible to naively
internalize some diagrams: some internal adjunction leads a
contradiction \parencite{licata2018internal}; there are only trivial
internal idempotent comonads \parencite{shulman2018brouwer}.
\noindent 
While there is no chance of naive internalization of such interesting
but problematic diagrams to plain homotopy type theory, some other
diagrams can be internalized in a clever way pointed out by
Shulman\footnote{\url{https://golem.ph.utexas.edu/category/2011/11/internalizing_the_external_or.html}}. A
minimal non-trivial example is a diagram consisting of two
\(\infty\)-logoses and a lex, accessible functor between them in one
direction. The two \(\infty\)-logoses are \emph{lex, accessible
  localizations} of another \(\infty\)-logos obtained by the \emph{Artin
  gluing} for the functor, and the functor is reconstructed by
composing the inclusion from one localization and the reflector to the
other. Moreover, this reconstruction is \emph{internal} to the glued
\(\infty\)-logos, because lex, accessible localizations of an
\(\infty\)-logos are expected to correspond to \emph{lex, accessible
  modalities} in its internal language. Hence, plain homotopy type
theory as a language for the glued \(\infty\)-logos is sufficient to reason
about the original diagram.

In this paper, we propose a class of shapes of diagrams of
\(\infty\)-logoses for which the internal reconstruction technique
explained in the previous paragraph works. We call shapes in the
proposed class \emph{mode sketches}. Our main results are summarized
as follows. Let \(\modesketch\) be a mode sketch.
\begin{enumerate}
\item We associate to \(\modesketch\) certain axioms in type theory,
  one of which is to postulate some lex, accessible modalities from
  which one can construct a diagram of \(\infty\)-logoses internally to
  type theory (\cref{sec:modes,sec:mode-sketches-1}).
\item We show that the \((\infty,1)\)-category of models of the axioms associated to
  \(\modesketch\) in \(\infty\)-logoses is equivalent to the \((\infty,1)\)-category of
  diagrams of \(\infty\)-logoses indexed over \(\modesketch\)
  (\cref{thm:main-theorem}), where \(\modesketch\) is regarded as a
  presentation of an \((\infty, 2)\)-category. The right to left
  construction is given by \emph{oplax limits}, a generalization of
  the Artin gluing \parencite{wraith1974artin,shulman2015inverse}.
\end{enumerate}
In \cref{sec:synth-tait-comp-1,sec:modal-local,sec:oplax-limits-1}
below, we explain other new results
(\cref{prop:canonical-join-accessible,thm:equiv-two-mode-sketch-axioms,prop:logos-accessible-oplax-limit,prop:oplax-univalence,prop:mate-correspondence}).

This paper is an extended version of our conference paper
\parencite{uemura2023diagrams}. The conference version presents
\cref{prop:canonical-join-accessible} and proof sketches of
\cref{thm:main-theorem,thm:equiv-two-mode-sketch-axioms}.  The current
version includes every detail of those results. The statement of
\cref{thm:main-theorem} is modified not to rely on an interpretation
of type theory in \(\infty\)-logoses.

\subsection{Modalities in homotopy type theory}
\label{sec:modal-local}

A \emph{modality} in homotopy type theory
\parencite{rijke2020modalities,christensen2020localization,christensen2022characterizations}
is a subuniverse satisfying certain conditions. A modality that is
moreover lex and accessible is expected to correspond to a
sub-\(\infty\)-topos or a localization of a \(\infty\)-logos
\parencite{anel2022left-exact,vergura2019localization-arxiv}. The
\emph{fracture and gluing theorem} of Rijke, Shulman, and Spitters \parencite[Theorem
3.50]{rijke2020modalities} gives a construction of the join
\(\mode \lor \modeI\) of two lex modalities \(\mode\) and \(\modeI\)
under some assumption. The join obtained by this theorem satisfies
that every type \(\ty\) in \(\mode \lor \modeI\) is canonically
fractured into a type \(\ty_{\modeI}\) in \(\modeI\) and a type family
\(\ty_{\mode}\) on \(\ty_{\modeI}\) valued in \(\mode\), and \(\ty\)
is reconstructed as
\(\ty \simeq \ilexists_{\var : \ty_{\modeI}} \ty_{\mode}(\var)\).

In this paper, we improve the fracture and gluing theorem. We show
that the construction of joins of lex modalities preserves
accessibility as well (\cref{prop:canonical-join-accessible}).

\subsection{Synthetic Tait computability}
\label{sec:synth-tait-comp-1}

The fracture and gluing theorem is used in
Sterling's \emph{synthetic Tait computability}
\parencite{sterling2021logical,sterling2021thesis}. It is a technique
of working with \emph{logical relations}, which are used in the study
of type theories and programming languages, in an internal language
for the Artin gluing and has applications to, for example,
normalization theorems for complex type theories
\parencite{sterling2021normalization,gratzer2021normalization}. There
the fracture and gluing theorem is instantiated by the \emph{closed}
and \emph{open} modalities associated to a proposition. Then every
type in the internal language is canonically fractured into an open
type and a closed (unary, proof-relevant) relation on it which are
glued back together. The internal language for the Artin gluing is
thus a type theory with an indeterminate proposition in which
\emph{types are relations} and provides a synthetic method of working
with logical relations.

In this paper, we relate synthetic Tait computability and mode
sketches. The core axiom for synthetic Tait computability is to
postulate some indeterminate propositions. We show that part of the
axioms associated to a mode sketch is equivalent to postulating a
lattice of propositions (\cref{thm:equiv-two-mode-sketch-axioms}).

Mode sketches thus provide a synthetic method of working with logical
relations that is alternative to and generalizes synthetic Tait
computability. This is also natural from
Shulman's point of view
\parencite{shulman2015inverse} that interpretations of type theory in
oplax limits are generalized logical relations. Since we work in
homotopy type theory, what we get is actually \emph{higher-dimensional
  logical relations}, and our primary application of mode sketches in
upcoming paper(s) \parencite{uemura2022normalization-arxiv} will be
normalization for \emph{\(\infty\)-type theories} introduced by
Nguyen and Uemura \parencite{nguyen2022type-arxiv} as a higher-dimensional generalization
of type theories.

\subsection{Oplax limits of \((\infty, 1)\)-categories}
\label{sec:oplax-limits-1}

\emph{Oplax limits} are special \((\infty, 2)\)-categorical limits and
analogous to oplax limits in \(2\)-category theory
\parencite{street1976limits,kelly1989elementary,johnson2021categories}. Oplax
limits indexed over \((\infty, 1)\)-categories are studied by
Gepner, Haugseng, and Nikolaus \parencite{gepner2017lax} and generalized to arbitrary indexing
\((\infty, 2)\)-categories by
Gagna, Harpaz, and Lanari \parencite{gagna2020fibrations-arxiv}.

Oplax limits of \(\infty\)-logoses are of our
interest. Wraith \parencite{wraith1974artin} shows that the oplax limit of a
diagram of (elementary) logoses and lex functors is a
logos. Lurie \parencite[Proposition 6.3.2.3]{lurie2009higher} shows that the
conical limit of a diagram of \(\infty\)-logoses is an \(\infty\)-logos.

The oplax limit of a diagram classifies \emph{oplax natural
  transformations} from a constant diagram to the given diagram just
as a (conical) limit classifies natural transformations from a
constant diagram. It is known \parencite[Theorem
3.8.1]{lurie2009goodwillie} that natural transformations between
diagrams of \((\infty ,1)\)-categories correspond to fibered functors
between the \((\infty, 2)\)-categories of elements of the diagrams. Oplax
natural transformations between diagrams of \((\infty, 1)\)-categories are
to correspond to not necessarily fibered functors between the
\((\infty, 2)\)-categories of elements. Some special cases of this have
already been proved: Haugseng et al.\@ \parencite[Theorem E]{haugseng2021mate-arxiv} show
the case when the diagrams are indexed over an
\((\infty ,1)\)-category; Gagna, Harpaz, and Lanari \parencite[Corollary
4.4.3]{gagna2020fibrations-arxiv} show the case when the domain is
constant on the point.

The \emph{mate correspondence} is a useful source of oplax natural
transformations. In the \(2\)-categorical case
\parencite{street1972two}, it asserts that given two diagrams \(\fun\)
and \(\funI\) of categories, oplax natural transformations
\(\fun \to \funI\) that are point-wise left adjoints correspond to lax
natural transformations \(\funI \to \fun\) that are point-wise right
adjoints. An \((\infty, 2)\)-categorical version is proved by
Haugseng et al.\@ \parencite[Corollary F]{haugseng2021mate-arxiv} in the form of an
equivalence between the \((\infty, 1)\)-categories of oplax/lax natural
transformations in the special case when the diagrams are indexed over
an \((\infty, 1)\)-category.

In this paper, we show some new results on oplax limits of
\((\infty, 1)\)-categories. All of them are consequences of results in the
literature, but it is worth stating them explicitly. We prove an
\(\infty\)-analogue of the result of Wraith \parencite{wraith1974artin}: the oplax
limit of a diagram of \(\infty\)-logoses and lex, accessible functors is an
\(\infty\)-logos (\cref{prop:logos-accessible-oplax-limit}). We show that
oplax natural transformations between diagrams of
\((\infty, 1)\)-categories correspond to arbitrary functors between the
\((\infty, 2)\)-categories of elements (\cref{prop:oplax-univalence}) by
reducing it to the special cases proved by
Haugseng et al.\@ \parencite{haugseng2021mate-arxiv} and Gagna, Harpaz, and Lanari \parencite{gagna2020fibrations-arxiv}. We show
the mate correspondence for diagrams indexed over an arbitrary
\((\infty, 2)\)-category in the form of an equivalence between the spaces
of oplax/lax natural transformations
(\cref{prop:mate-correspondence}).

\subsection{Organization}
\label{sec:organization}

The paper is split into two parts. The first part
(\cref{sec:modal-homot-type,sec:mode-sketches,sec:mode-sketch-synth})
provides the theory of mode sketches internally to type theory. The
second part (\cref{sec:high-categ-theory,sec:semant-mode-sketch}) is
devoted to semantics of mode sketches in \(\infty\)-logoses. The two parts
are independent of each other on a technical level.

In \cref{sec:modal-homot-type}, we review the theory of modalities in
homotopy type theory \parencite{rijke2020modalities}. Our focus is on
the poset of lex, accessible modalities and on the open and closed
modalities associated to propositions.

\Cref{sec:mode-sketches,sec:mode-sketch-synth} are the core of the
paper. We introduce the notion of a \emph{mode sketch}
(\cref{def:mode-sketch}). For every mode sketch, we introduce two
equivalent sets of axioms to encode a certain diagram of
universes. One postulates some lex, accessible modalities while the
other postulates a lattice of propositions. The open and closed
modalities give a construction of the former from the latter which we
show is an equivalence (\cref{thm:equiv-two-mode-sketch-axioms}). The
latter is a higher dimensional analogue of
Sterling's synthetic Tait computability
\parencite{sterling2021thesis}.

\Cref{sec:high-categ-theory} is a preliminary section needed for the
semantics of mode sketches.

Finally in \cref{sec:semant-mode-sketch}, we show our main result
(\cref{thm:main-theorem}): for any mode sketch, the \((\infty,1)\)-category of models of
the axioms associated to the mode sketch in \(\infty\)-logoses is
equivalent to the \((\infty,1)\)-category of diagrams of \(\infty\)-logoses and lex,
accessible functors indexed over the mode sketch.

\subsection{Related work}
\label{sec:related-work}

An earlier version of \emph{cohesive homotopy type theory}
\parencite{schreiber2012quantum} uses modalities in plain homotopy
type theory to internalize a series of adjunctions that arises in
Lawvere's axiomatic cohesion
\parencite{lawvere2007axiomatic}. However, because naive
internalization of adjunctions do not work well
\parencite{licata2018internal,shulman2018brouwer}, the axiomatization
is tricky and not ideal to work with. The newer version of cohesive
homotopy type theory \parencite{shulman2018brouwer} instead extends
homotopy type theory by another layer of context and new modal
operators. The resulting type theory works well for axiomatic cohesion
but is complicated compared to plain homotopy type theory. It is also
too optimized for axiomatic cohesion.

A more general framework for internal diagrams is \emph{multimodal
  dependent type theory} \parencite{gratzer2021multimodal}. It is
roughly a family of type theories related to each other via modal
operators and interpreted in a diagram of presheaf categories. The
shape of diagram is specified directly by an arbitrary \(2\)-category
which is called a \emph{mode theory} in this context. Our terminology
``mode sketch'' is chosen to mean a sketch of a mode
theory. Multimodal dependent type theory is potentially an internal
language for diagrams of \(\infty\)-logoses, but for this one would have to
rectify not only \(\infty\)-logoses but also functors and natural
transformations between them.

Our work brings back the ideas of earlier cohesive homotopy type
theory. Although it might not be the best type theory, it has a lot of
advantages: modalities are internal to plain homotopy type theory, and
thus all results are ready to formalize in existing proof assistants;
keeping type theory simple is also important in informal use of type
theory in which the correctness of application of inference rules is
not checked by computer; the semantics is no more complicated than the
\(\infty\)-logos semantics of homotopy type theory; it also opens the door
to internalization of more general diagrams in a uniform way, which is
the motivation for the current work.

\section{Modalities in homotopy type theory}
\label{sec:modal-homot-type}

We review the theory of \emph{modalities} in homotopy type theory
\parencite{rijke2020modalities}. In this section, we work in homotopy
type theory. By \emph{homotopy type theory} we mean dependent type
theory with (dependent) function types, (dependent) pair types, a unit
type, identity types (without equality reflection), at least two
univalent universes
\begin{math}
  \univ : \enlarge \univ,
\end{math}
an empty type, pushouts, and localizations \parencite[Section
2.2]{rijke2020modalities}. Note that truncations are instances of
localization. We mainly follow the HoTT Book \parencite{hottbook} for
terminologies and notations in homotopy type theory.

A modality is in short a reflective subuniverse closed
under pair types.

\begin{definition}
  A \emph{subuniverse} \(\mode\) is a function
  \(\ilIn_{\mode} : \univ \to \enlarge \univ\) such that
  \(\ilIn_{\mode}(\ty)\) is a proposition for all \(\ty : \univ\). A
  type \(\ty\) satisfying \(\ilIn_{\mode}(\ty)\) is called
  \emph{\(\mode\)-modal}. We define a subtype
  \(\univ_{\mode} \subset \univ\) to be
  \(\{\ty : \univ \mid \ilIn_{\mode}(\ty)\}\).
\end{definition}

\begin{definition}
  A subuniverse \(\mode\) is \emph{reflective} if it is equipped with
  functions \(\opModality_{\mode} : \univ \to \univ_{\mode}\) and
  \(\unitModality_{\mode} : \ilforall_{\ty : \univ} \ty \to
  \opModality_{\mode} \ty\) such that that the precomposition
  \(\ilAbs \map. \map \comp \unitModality_{\mode}(\ty) :
  (\opModality_{\mode} \ty \to \tyI) \to (\ty \to \tyI)\) is an
  equivalence for any \(\tyI : \univ_{\mode}\). Note that such a pair
  \((\opModality_{\mode}, \unitModality_{\mode})\) is unique.
\end{definition}

\begin{definition}
  A reflective subuniverse \(\mode\) is a \emph{modality} if
  \(\ilIn_{\mode}\) is closed under pair types, that is, for
  \(\ty : \univ\) and \(\tyI : \ty \to \univ\), if
  \(\ilIn_{\mode}(\ty)\) and
  \(\ilforall_{\el : \ty} \ilIn_{\mode}(\tyI(\el))\), then
  \(\ilIn_{\mode}(\ilexists_{\el : \ty} \tyI(\el))\).
\end{definition}

An important class of modalities is \emph{accessible} modalities which
are roughly modalities ``presented by small data''.

\begin{definition}
  For types \(\ty, \tyI : \univ\), we say \(\ty\) is \emph{left
    orthogonal} to \(\tyI\) or \(\tyI\) is \emph{right orthogonal} to
  \(\ty\) and write \(\ty \relOrth \tyI\) if the function
  \(\ilAbs (\elI : \tyI). \ilAbs (\blank : \ty). \elI : \tyI \to (\ty \to
  \tyI)\) is an equivalence. For a subuniverse \(\mode\), we define
  subuniverses \(\mode^{\orthMark}\) and \({}^{\orthMark}\mode\) by
  \begin{align*}
    \begin{autobreak}
      \ilIn_{\mode^{\orthMark}}(\tyI) \defeq
      \ilforall_{\ty : \univ_{\mode}}
      \ty \relOrth \tyI
    \end{autobreak}
    \\
    \begin{autobreak}
      \ilIn_{{}^{\orthMark}\mode}(\ty) \defeq
      \ilforall_{\tyI : \univ_{\mode}}
      \ty \relOrth \tyI.
    \end{autobreak}
  \end{align*}
\end{definition}

\begin{definition}
  A \emph{null generator} \(\nullgen\) consists of
  \(\idxNullGen_{\nullgen} : \univ\) and
  \(\tyNullGen_{\nullgen} : \idxNullGen_{\nullgen} \to \univ\). We write
  \(\ilNullGen\) for the type of null generators. Given a null
  generator \(\nullgen\), we define a subuniverse
  \(\modeNull(\nullgen)\) by
  \begin{math}
    \ilIn_{\modeNull(\nullgen)}(\ty) \defeq \ilforall_{\idx :
      \idxNullGen_{\nullgen}} \tyNullGen_{\nullgen}(\idx) \relOrth \ty.
  \end{math}
  It is shown that \(\modeNull(\nullgen)\) is a modality using a
  higher inductive type \parencite[Theorem
  2.19]{rijke2020modalities}. A modality \(\mode\) is
  \emph{accessible} if it is in the image of \(\modeNull\), that is,
  \begin{math}
    \trunc{\ilexists_{\nullgen : \ilNullGen} \mode = \modeNull(\nullgen)}.
  \end{math}
\end{definition}

Another important class of modalities is \emph{lex} modalities.

\begin{definition}
  For a modality \(\mode\), a type \(\ty : \univ\) is
  \emph{\(\mode\)-connected} if \(\opModality_{\mode} \ty\) is
  contractible. This is equivalent to
  \(\ilIn_{{}^{\orthMark}\mode}(\ty)\) by \parencite[Corollary
  1.37]{rijke2020modalities}.
\end{definition}

\begin{definition}
  A modality \(\mode\) is
  \emph{lex} if for any \(\mode\)-connected type \(\ty : \univ\), the
  identity type \(\el_{1} = \el_{2}\) is \(\mode\)-connected for any
  \(\el_{1}, \el_{2} : \ty\).
\end{definition}

Modalities that are both lex and accessible are of particular
importance because they correspond to subtoposes of an
\(\infty\)-topos under the interpretation of types as sheaves on the
\(\infty\)-topos. From now on, we are mostly interested lex, accessible
modalities, so we give them a short name.

\begin{terminology}
  {\acrLAM} is an acronym for lex, accessible modality.
\end{terminology}

Fundamental examples of {\acrLAMs} are \emph{open} and \emph{closed}
modalities which correspond to open and closed, respectively,
subtoposes.

\begin{construction}
  Let \(\propo\) be a proposition. We define the
  \emph{open modality} \(\modeOpen(\propo)\) by
  \begin{math}
    \opModality_{\modeOpen(\propo)} \ty \defeq (\propo \to \ty)
  \end{math}
  and
  \begin{math}
    \unitModality_{\modeOpen(\propo)}(\ty, \el)
    \defeq \ilAbs \blank. \el.
  \end{math}
  It is lex and accessible by \parencite[Example
  2.24 and Example 3.10]{rijke2020modalities}. We also define the \emph{closed
    modality} \(\modeClosed(\propo)\) by
  \begin{math}
    \ilIn_{\modeClosed(\propo)}(\ty) \defeq (\propo \to \ilIsContr(\ty)).
  \end{math}
  It is lex and accessible by \parencite[Example 2.25 and Example
  3.14]{rijke2020modalities}. Note that
  \(\modeClosed(\propo) = {}^{\orthMark} \modeOpen(\propo)\)
  \parencite[Example 1.31]{rijke2020modalities}.
\end{construction}

\subsection{The poset of lex, accessible modalities}
\label{sec:poset-lex-accessible}

Let \(\ilSU\) denote the poset of subuniverses where
\(\mode \le \modeI\) if \(\ilIn_{\mode}(\ty) \to \ilIn_{\modeI}(\ty)\) for
every \(\ty : \univ\). We have the full subposets of \(\ilSU\)
\begin{equation*}
  \ilRSU
  \supset \ilModality
  \supset \ilAccModality
  \supset \ilLexAcc
\end{equation*}
consisting of reflective subuniverses, modalities, accessible
modalities, and lex, accessible modalities, respectively. We also have
the full subposet
\begin{math}
  \ilLex \subset \ilModality
\end{math}
of lex modalities. By definition,
\(\ilLexAcc = \ilLex \cap \ilAccModality\). We study the poset
\(\ilLexAcc\) in more detail.
\begin{defiC}[{\cite[Theorem 3.25]{rijke2020modalities}}]
  Let \(\idxsh : \univ\) and \(\mode : \idxsh \to \ilLexAcc\). A
  \emph{canonical meet} \(\bigland_{\idx : \idxsh} \mode(\idx)\) is a
  {\acrLAM} that is the meet of \(\mode(\idx)\)'s in \(\ilSU\). A
  \emph{canonical join} \(\biglor_{\idx : \idxsh} \mode(\idx)\) is a
  {\acrLAM} satisfying that a type \(\ty : \univ\) is
  \((\biglor_{\idx : \idxsh} \mode(\idx))\)-connected if and only if
  it is \(\mode(\idx)\)-connected for all \(\idx : \idxsh\). Note that
  a canonical join is the join in \(\ilModality\).
\end{defiC}

\begin{example}
  The \emph{top modality} \(\modeTop\), for which all the types are
  modal, is the canonical meet of the empty family. The \emph{bottom
    modality} \(\modeBottom\), for which only the contractible types
  are modal, is the canonical join of the empty family.
\end{example}

The canonical meet of an arbitrary family of {\acrLAMs} exists
\parencite[Theorem 3.29 and Remark
3.23]{rijke2020modalities}. Canonical joins are less understood than
canonical meets. One important case when canonical joins exist and can
be computed is the following.

\begin{definition}
  Let \(\mode\) and \(\modeI\) be {\acrLAMs}. \(\modeI\) is
  \emph{strongly disjoint from \(\mode\)} if any \(\mode\)-modal type
  is \(\modeI\)-connected or equivalently if
  \(\mode \le {}^{\orthMark}\modeI\) in \(\ilSU\).
\end{definition}

\begin{proposition}[{Fracture and gluing theorem}]
  \label{prop:join-strongly-disjoint}
  Let \(\mode\) and \(\modeI\) be {\acrLAMs} such that
  \(\mode \leq {}^{\orthMark}\modeI\).
  \begin{enumerate}
  \item The canonical join \(\mode \lor \modeI\) exists.
  \item A type \(\ty\) is \((\mode \lor \modeI)\)-modal if and only if
    the function
    \(\unitModality_{\modeI}(\ty) : \ty \to \opModality_{\modeI} \ty\)
    has \(\mode\)-modal fibers.
  \item
    \(\univ_{\mode \lor \modeI} \simeq \ilexists_{\ty : \univ_{\mode}}
    \ilexists_{\tyI : \univ_{\modeI}} \ty \to
    \opModality^{\modeI}_{\mode} \tyI\).
  \end{enumerate}
  In the special case when \(\mode = {}^{\orthMark} \modeI\), we have
  \(\mode \lor \modeI = \modeTop\).
\end{proposition}
\begin{proof}
  All but the accessibility of \(\mode \lor \modeI\) are proved by
  Rijke, Shulman, and Spitters \parencite[Theorem 3.50]{rijke2020modalities}. We will prove the
  accessibility of \(\mode \lor \modeI\) in
  \cref{prop:canonical-join-accessible} below using an open modality.
\end{proof}

The distributive law holds in some special cases.

\begin{fact}[{\cite[Theorem 3.30]{rijke2020modalities}}]
  \label{prop:meet-preserves-modal}
  Let \(\mode\) and \(\modeI\) be {\acrLAMs}. If
  \(\opModality_{\mode}\) preserves \(\modeI\)-modal types, then
  \begin{math}
    \opModality_{\mode \land \modeI} \ty =
    \opModality_{\mode} \opModality_{\modeI} \ty.
  \end{math}
\end{fact}

\begin{proposition}
  \label{prop:join-strongly-disjoint-distributive}
  Let \(\mode_{1}\), \(\mode_{2}\), and \(\mode_{3}\) be
  {\acrLAMs}. Suppose \(\mode_{1} \le {}^{\orthMark} \mode_{3}\) and
  \(\mode_{2} \le {}^{\orthMark} \mode_{3}\). Then
  \begin{math}
    (\mode_{1} \lor \mode_{3}) \land (\mode_{2} \lor \mode_{3})
    = (\mode_{1} \land \mode_{2}) \lor \mode_{3}.
  \end{math}
\end{proposition}
\begin{proof}
  Note that the right side exists as
  \(\mode_{1} \land \mode_{2} \le \mode_{1} \le {}^{\orthMark}
  \mode_{3}\). Let \(\ty : \univ\). By
  \cref{prop:join-strongly-disjoint}, \(\ty\) is
  \(((\mode_{1} \lor \mode_{3}) \land (\mode_{2} \lor
  \mode_{3}))\)-modal if and only if fibers of
  \(\unitModality_{\mode_{3}}(\ty)\) are both \(\mode_{1}\)-modal and
  \(\mode_{2}\)-modal, but this is equivalent to that \(\ty\) is
  \((\mode_{1} \land \mode_{2}) \lor \mode_{3}\)-modal again by
  \cref{prop:join-strongly-disjoint}.
\end{proof}

\begin{proposition}
  \label{prop:meet-preserves-modal-distribute-disjoint}
  Let \(\mode_{1}\), \(\mode_{2}\), and \(\mode_{3}\) be
  {\acrLAMs}. Suppose that \(\mode_{1} \le {}^{\orthMark}\mode_{2}\) and
  that \(\opModality_{\mode_{2}}\) preserves \(\mode_{3}\)-modal
  types. Then
  \((\mode_{1} \lor \mode_{2}) \land \mode_{3} = (\mode_{1} \land
  \mode_{3}) \lor (\mode_{2} \land \mode_{3})\).
\end{proposition}
\begin{proof}
  Note that the right side exists as
  \begin{math}
    \mode_{1} \land \mode_{3}
    \le \mode_{1}
    \le {}^{\orthMark} \mode_{2}
    \le {}^{\orthMark} (\mode_{2} \land \mode_{3}).
  \end{math}
  Let \(\ty\) be a
  \(((\mode_{1} \lor \mode_{2}) \land \mode_{3})\)-modal type. We show
  that \(\ty\) is
  \(((\mode_{1} \land \mode_{3}) \lor (\mode_{2} \land
  \mode_{3}))\)-modal, which is by \cref{prop:join-strongly-disjoint}
  equivalent to that
  \(\unitModality_{\mode_{2} \land \mode_{3}}(\ty) : \ty \to
  \opModality_{\mode_{2} \land \mode_{3}} \ty\) has
  \((\mode_{1} \land \mode_{3})\)-modal fibers. Since
  \(\opModality_{\mode_{2}}\) preserves \(\mode_{3}\)-modal types and
  since \(\ty\) is \(\mode_{3}\)-modal,
  \(\unitModality_{\mode_{2} \land \mode_{3}}(\ty)\) is equivalent to
  \(\unitModality_{\mode_{2}}(\ty) : \ty \to \opModality_{\mode_{2}}
  \ty\) by \cref{prop:meet-preserves-modal}. The fibers of
  \(\unitModality_{\mode_{2}}(\ty)\) are \(\mode_{1}\)-modal by
  \cref{prop:join-strongly-disjoint} and \(\mode_{3}\)-modal since
  both domain and codomain are \(\mode_{3}\)-modal.
\end{proof}

\subsection{Accessibility of the canonical join}
\label{sec:access-canon-join}

Let us fill the gap in the proof of
\cref{prop:join-strongly-disjoint}.

\begin{proposition}
  \label{prop:canonical-join-accessible}
  Let \(\mode\) and \(\modeI\) be {\acrLAMs} such that
  \(\mode \le {}^{\orthMark} \modeI\). Then the canonical join
  \(\mode \lor \modeI\) (in \(\ilLex\)) is accessible.
\end{proposition}

We have to find a null generator for \(\mode \lor \modeI\). A natural
guess is the following.

\begin{construction}
  Let \(\nullgen\) and \(\nullgenI\) be null generators. We define a
  null generator \(\nullgen \join \nullgenI\) by
  \begin{math}
    \idxNullGen_{\nullgen \join \nullgenI} \defeq
    \idxNullGen_{\nullgen} \times \idxNullGen_{\nullgenI}
  \end{math}
  and
  \begin{math}
    \tyNullGen_{\nullgen \join \nullgenI}(\idx, \idxI) \defeq
    \tyNullGen_{\nullgen}(\idx) \join \tyNullGen_{\nullgenI}(\idxI)
    \defeq \tyNullGen_{\nullgen}(\idx)
    +_{\tyNullGen_{\nullgen}(\idx) \times \tyNullGen_{\nullgenI}(\idxI)}
    \tyNullGen_{\nullgenI}(\idxI).
  \end{math}
\end{construction}

\begin{lemma}
  \label{lem:nullgen-join-connected}
  Let \(\mode\) and \(\modeI\) be {\acrLAMs}, and let \(\nullgen\) and
  \(\nullgenI\) be null generators for \(\mode\) and \(\modeI\),
  respectively. Then
  \(\tyNullGen_{\nullgen \join \nullgenI}(\idx, \idxI)\) is both
  \(\mode\)-connected and \(\modeI\)-connected for all
  \(\idx : \idxNullGen_{\nullgen}\) and
  \(\idxI : \idxNullGen_{\nullgenI}\).
\end{lemma}
\begin{proof}
  Recall that a function is \emph{\(\mode\)-connected} if its fibers
  are \(\mode\)-connected and that the class of \(\mode\)-connected
  functions is the left class of a (stable) orthogonal factorization
  system \parencite[Theorem 1.34]{rijke2020modalities}. Then the claim
  follows by the pushout stability and the right cancellability of
  \(\mode\)-connected and \(\modeI\)-connected functions.
\end{proof}

\Cref{lem:nullgen-join-connected} shows
\(\mode \lor \modeI \le \modeNull(\nullgen \join \nullgenI)\) for
arbitrary accessible modalities \(\mode\) and \(\modeI\) and for
arbitrary choices of \(\nullgen\) and \(\nullgenI\). We know neither
if the other direction holds in general for some choices of
\(\nullgen\) and \(\nullgenI\) nor if
\(\modeNull(\nullgen \join \nullgenI)\) is independent of \(\nullgen\)
and \(\nullgenI\). Note that Finster \parencite{finster2021left-note} observed
that \(\modeNull(\nullgen \join \nullgenI)\) is lex whenever
\(\modeNull(\nullgen)\) and \(\modeNull(\nullgenI)\) are lex. In the
special case when \(\mode \le {}^{\orthMark} \modeI\), the idea of the
proof of \(\mode \lor \modeI = \modeNull(\nullgen \join \nullgenI)\)
is to show that \(\modeI\) is an \emph{open modality} within the
subuniverse of \(\modeNull(\nullgen \join \nullgenI)\)-modal types.

\begin{lemma}
  \label{lem:disjoint-mode-open}
  Let \(\mode\) and \(\modeI\) be {\acrLAMs} such that
  \(\mode \le {}^{\orthMark} \modeI\). Then
  \(\modeI \le \modeOpen(\opModality_{\mode} \ilEmpty)\).
\end{lemma}
\begin{proof}
  This is because \(\opModality_{\mode} \ilEmpty\) is
  \(\modeI\)-connected by assumption.
\end{proof}

\begin{lemma}
  \label{lem:relative-partition-open-modality}
  Let \(\mode\) and \(\modeI\) be {\acrLAMs} such that
  \(\mode \le {}^{\orthMark} \modeI\). Suppose that \(\nullgen\) and
  \(\nullgenI\) are null generators for \(\mode\) and \(\modeI\),
  respectively, and that \(\nullgen\) admits a function
  \(\map : \opModality_{\mode} \ilEmpty \to \idxNullGen_{\nullgen}\)
  such that \(\ilEmpty \simeq \tyNullGen_{\nullgen}(\map(\idx))\) for all
  \(\idx : \opModality_{\mode} \ilEmpty\). Then
  \(\opModality_{\modeOpen(\opModality_{\mode} \ilEmpty)} \ty\) is
  \(\modeI\)-modal for any \(\modeNull(\nullgen \join \nullgenI)\)-modal
  type \(\ty\). Consequently, the canonical function
  \begin{math}
    \opModality_{\modeOpen(\opModality_{\mode} \ilEmpty)} \ty
    \to \opModality_{\modeI} \ty
  \end{math}
  induced by \cref{lem:disjoint-mode-open} is an equivalence for any
  \(\modeNull(\nullgen \join \nullgenI)\)-modal type \(\ty\).
\end{lemma}
\begin{proof}
  We show that
  \(\opModality_{\modeOpen(\opModality_{\mode} \ilEmpty)} \ty \defeq
  (\opModality_{\mode} \ilEmpty \to \ty)\) is \(\modeI\)-modal. Since
  \(\nullgenI\) is a null generator for \(\modeI\), it suffices to
  show that
  \(\tyNullGen_{\nullgenI}(\idxI) \relOrth (\opModality_{\mode}
  \ilEmpty \to \ty)\) for all \(\idxI : \idxNullGen_{\nullgenI}\). This
  is equivalent to that \(\tyNullGen_{\nullgenI}(\idxI) \relOrth \ty\)
  under an assumption \(\idx : \opModality_{\mode} \ilEmpty\). This holds since
  \begin{math}
    \tyNullGen_{\nullgenI}(\idxI) \simeq
    \ilEmpty \join \tyNullGen_{\nullgenI}(\idxI)
    \simeq \tyNullGen_{\nullgen}(\map(\idx)) \join \tyNullGen_{\nullgenI}(\idxI)
  \end{math}
  and since \(\ty\) is \(\modeNull(\nullgen \join \nullgenI)\)-modal.
\end{proof}

\begin{lemma}
  \label{lem:relative-partition-closed-modality}
  Let \(\mode\) and \(\modeI\) be {\acrLAMs} such that
  \(\mode \le {}^{\orthMark} \modeI\). Suppose that \(\nullgen\) and
  \(\nullgenI\) are null generators for \(\mode\) and \(\modeI\),
  respectively, and that \(\nullgenI\) has an element
  \(\idxI : \idxNullGen_{\nullgenI}\) such that
  \(\tyNullGen_{\nullgenI}(\idxI) \simeq \opModality_{\mode}
  \ilEmpty\). Then, if a type \(\ty\) is
  \(\modeNull(\nullgen \join \nullgenI)\)-modal and
  \(\modeOpen(\opModality_{\mode} \ilEmpty)\)-connected, then it is
  \(\mode\)-modal.
\end{lemma}
\begin{proof}
  We show that \(\tyNullGen_{\nullgen}(\idx) \relOrth \ty\) for all
  \(\idx : \idxNullGen_{\nullgen}\). By the definition of \({\join}\),
  we have the following pullback square.
  \begin{equation*}
    \begin{tikzcd}
      (\tyNullGen_{\nullgen}(\idx) \join \opModality_{\mode} \ilEmpty \to \ty)
      \arrow[r, "\simeq"]
      \arrow[d]
      \arrow[dr, pbMark] &
      (\tyNullGen_{\nullgen}(\idx) \to \ty)
      \arrow[d] \\
      (\opModality_{\mode} \ilEmpty \to \ty)
      \arrow[r, "\simeq"'] &
      (\tyNullGen_{\nullgen}(\idx) \to \opModality_{\mode} \ilEmpty \to \ty)
    \end{tikzcd}
  \end{equation*}
  Since \(\ty\) is
  \(\modeOpen(\opModality_{\mode} \ilEmpty)\)-connected, the domain
  and codomain of the bottom function are contractible, and thus the
  bottom function is an equivalence. It then follows that the top
  function is also an equivalence. Since \(\ty\) is
  \(\modeNull(\nullgen \join \nullgenI)\)-modal and since
  \(\tyNullGen_{\nullgenI}(\idxI) \simeq \opModality_{\mode} \ilEmpty\), we
  have
  \begin{math}
    \ty \simeq
    (\tyNullGen_{\nullgen}(\idx) \join \opModality_{\mode} \ilEmpty \to \ty)
    \simeq (\tyNullGen_{\nullgen}(\idx) \to \ty),
  \end{math}
  and thus \(\tyNullGen_{\nullgen}(\idx) \relOrth \ty\).
\end{proof}

\begin{proof}[Proof of \cref{prop:canonical-join-accessible}]
  Let \(\nullgen\) and \(\nullgenI\) be null generators for \(\mode\)
  and \(\modeI\), respectively. Note that the null generator obtained
  from a null generator for \(\mode\) by adjoining a family of
  \(\mode\)-connected types yields the same modality \(\mode\). Under
  an assumption \(\idx : \opModality_{\mode} \ilEmpty\), the empty
  type \(\ilEmpty\) becomes \(\mode\)-connected, and thus we may
  assume that \(\nullgen\) includes the type family
  \(\ilAbs (\blank : \opModality_{\mode} \ilEmpty). \ilEmpty\). Since
  \(\opModality_{\mode} \ilEmpty\) is \(\modeI\)-connected by
  assumption, we may assume that \(\nullgenI\) includes the type
  family \(\ilAbs (\blank : \ilUnit). \opModality_{\mode} \ilEmpty\).

  We show that
  \(\modeNull(\nullgen \join \nullgenI) = \mode \lor \modeI\). By
  \cref{lem:nullgen-join-connected},
  \(\mode \lor \modeI \le \modeNull(\nullgen \join \nullgenI)\). For the
  other direction, suppose that \(\ty\) is a
  \(\modeNull(\nullgen \join \nullgenI)\)-modal type. By
  \parencite[Theorem 3.50]{rijke2020modalities}, it suffices to show
  that
  \(\unitModality_{\mode}(\ty) : \ty \to \opModality_{\modeI} \ty\) has
  \(\mode\)-modal fibers. By
  \cref{lem:relative-partition-open-modality},
  \(\opModality_{\modeI} \ty \simeq
  \opModality_{\modeOpen(\opModality_{\mode} \ilEmpty)} \ty\). Then
  the fibers of \(\unitModality_{\mode}(\ty)\) are
  \(\modeOpen(\opModality_{\mode} \ilEmpty)\)-connected. Since both
  \(\ty\) and \(\opModality_{\modeI} \ty\) are
  \(\modeNull(\nullgen \join \nullgenI)\)-modal, the fibers of
  \(\unitModality_{\mode}(\ty)\) are also
  \(\modeNull(\nullgen \join \nullgenI)\)-modal. Thus, by
  \cref{lem:relative-partition-closed-modality},
  \(\unitModality_{\mode}(\ty)\) has \(\mode\)-modal fibers.
\end{proof}

As a by-product, we have the following.

\begin{corollary}
  \label{cor:partition-open-modality}
  Let \(\mode\) and \(\modeI\) be {\acrLAMs} such that
  \(\mode \le {}^{\orthMark} \modeI\). If
  \(\mode \lor \modeI = \modeTop\), then \(\mode\) and \(\modeI\) are
  the closed and open, respectively, modalities associated to the
  proposition \(\opModality_{\mode} \ilEmpty\). \qed
\end{corollary}

\subsection{Open and closed modalities}
\label{sec:open-clos-modal}

We collect some properties of \emph{open modalities} and \emph{closed
  modalities}.

\begin{proposition}
  \label{prop:closed-modality-embedding}
  The function \(\modeClosed : \ilProp \to \ilLexAcc\) is a
  contravariant full embedding of posets and takes arbitrary joins to
  canonical meets.
\end{proposition}
\begin{proof}
  \(\modeClosed\) takes joins to canonical meets by \parencite[Example
  3.27]{rijke2020modalities}. It then follows that \(\modeClosed\)
  contravariantly preserves ordering. To see that \(\modeClosed\)
  reflects ordering, let \(\propo\) and \(\propoI\) be propositions and
  suppose that \(\modeClosed(\propoI) \le \modeClosed(\propo)\). By
  definition, \(\propoI\) is \(\modeClosed(\propoI)\)-modal and thus
  \(\modeClosed(\propo)\)-modal, but this implies that
  \(\propo \to \propoI\).
\end{proof}

\begin{proposition}
  \label{prop:open-modality-embedding}
  The function \(\modeOpen : \ilProp \to \ilLexAcc\) is a covariant full
  embedding of posets and takes finite meets to canonical meets and
  arbitrary joins to canonical joins.
\end{proposition}
\begin{proof}
  Since \(\modeClosed(\propo) = {}^{\orthMark} \modeOpen(\propo)\), the
  function \(\modeOpen\) is a covariant full embedding of posets and
  takes arbitrary joins to canonical joins by
  \cref{prop:closed-modality-embedding}. \(\modeOpen\) takes finite
  meets to canonical meets by \parencite[Example
  3.26]{rijke2020modalities}.
\end{proof}

\begin{proposition}
  \label{prop:meet-open-closed}
  Let \(\propo\) be a proposition and \(\mode\) a
  {\acrLAM}. \(\opModality_{\modeOpen(\propo)}\) preserves
  \(\mode\)-modal types, and \(\opModality_{\mode}\) preserves
  \(\modeClosed(\propo)\)-types.
\end{proposition}
\begin{proof}
  Let \(\ty : \univ\). If \(\ty\) is \(\mode\)-modal, then
  \(\opModality_{\modeOpen(\propo)} \ty \defeq \propo \to \ty\) is
  \(\mode\)-modal by \parencite[Lemma 1.26]{rijke2020modalities}. If
  \(\ty\) is \(\modeClosed(\propo)\)-modal, then
  \(\opModality_{\mode} \ty\) is \(\modeClosed(\propo)\)-modal by
  \parencite[Lemma 1.27]{rijke2020modalities}.
\end{proof}

\begin{proposition}
  \label{prop:open-closed-decomposition}
  \begin{math}
    \mode =
    (\mode \land \modeClosed(\propo))
    \lor (\mode \land \modeOpen(\propo))
  \end{math}
  for any {\acrLAM} \(\mode\) and any proposition \(\propo\).
\end{proposition}

\begin{proof}
  \begin{math}
    \mode =
    \mode \land (\modeClosed(\propo) \lor \modeOpen(\propo)) =
    (\mode \land \modeClosed(\propo)) \lor (\mode \land \modeOpen(\propo)),
  \end{math}
  where the first identification is by
  \cref{prop:join-strongly-disjoint} and the second is by
  \cref{prop:meet-preserves-modal-distribute-disjoint,prop:meet-open-closed}.
\end{proof}

\section{Mode sketches}
\label{sec:mode-sketches}

We introduce \emph{mode sketches} as shapes of diagrams of
subuniverses definable internally to type theory. We work in homotopy
type theory through the section.

\subsection{Internal diagrams induced by modalities}
\label{sec:modes}

We consider postulating some {\acrLAMs} to encode some diagram of
subuniverses. The fundamental observation is that a pair of {\acrLAMs}
induces a canonical functor between them.

\begin{construction}
  Let \(\mode\) and \(\modeI\) be {\acrLAMs}. We define a function
  \(\opModality^{\modeI}_{\mode} : \univ_{\modeI} \to \univ_{\mode}\) to
  be the composite of the inclusion \(\univ_{\modeI} \subset \univ\) and
  \(\opModality_{\mode} : \univ \to \univ_{\mode}\).
\end{construction}

\begin{remark}
  We can say that \(\opModality^{\modeI}_{\mode}\) is a functor
  \emph{externally}: we can construct a function
  \begin{math}
    \ilforall_{\ty, \tyI : \univ_{\modeI}} (\ty \to \tyI) \to
    (\opModality^{\modeI}_{\mode} \ty \to \opModality^{\modeI}_{\mode}
    \tyI)
  \end{math}
  and every instance of the coherence laws. However, it is not known
  how to state that \(\opModality^{\modeI}_{\mode}\) is a functor
  internally to type theory, because defining the type of
  \((\infty, 1)\)-categories in plain homotopy type theory is still an open
  problem.
\end{remark}

We have two functors
\(\opModality^{\mode}_{\modeI} : \univ_{\mode} \to \univ_{\modeI}\) and
\(\opModality^{\modeI}_{\mode} : \univ_{\modeI} \to \univ_{\mode}\) for
every pair of {\acrLAMs} \(\mode\) and \(\modeI\), but we are often
interested in only one direction. It is thus useful to cut off one
direction by postulating that \(\mode \le {}^{\orthMark}\modeI\): by the
definition of \(\modeI\)-connectedness,
\(\opModality^{\mode}_{\modeI}\) becomes constant at the unit
type. The other direction
\(\opModality^{\modeI}_{\mode} : \univ_{\modeI} \to \univ_{\mode}\)
remains non-trivial. Therefore, a pair \((\mode, \modeI)\) of
{\acrLAMs} such that \(\mode \le {}^{\orthMark}\modeI\) encodes a
functor \(\univ_{\modeI} \to \univ_{\mode}\). When
\(\modeI \le {}^{\orthMark} \mode\) is also assumed, \(\univ_{\mode}\)
and \(\univ_{\modeI}\) are considered unrelated.

Given more than two {\acrLAMs}, we have canonical natural
transformations between the canonical functors.

\begin{construction}
  Let \(\mode_{0}, \mode_{1}, \mode_{2}\) be {\acrLAMs}. We define
  \begin{equation*}
    \unitModality_{\mode_{1}}^{\mode_{0}; \mode_{2}} :
    \ilforall_{\ty : \univ_{\mode_{2}}}
    \opModality^{\mode_{2}}_{\mode_{0}} \ty
    \to \opModality^{\mode_{1}}_{\mode_{0}}
    \opModality^{\mode_{2}}_{\mode_{1}} \ty
  \end{equation*}
  by
  \begin{math}
    \unitModality_{\mode_{1}}^{\mode_{0}; \mode_{2}}(\ty) \defeq
    \opModality_{\mode_{0}} \unitModality_{\mode_{1}}(\ty).
  \end{math}
  This is well-typed as follows. \(\ty : \univ_{\mode_{2}}\) and
  \(\opModality_{\mode_{1}} \ty : \univ_{\mode_{1}}\) are implicitly
  coerced along the inclusions \(\univ_{\mode_{2}} \subset \univ\) and
  \(\univ_{\mode_{1}} \subset \univ\) respectively, and then
  \(\opModality_{\mode_{0}} \unitModality_{\mode_{1}}(\ty)\) has type
  \(\opModality_{\mode_{0}} \ty \to \opModality_{\mode_{0}}
  \opModality_{\mode_{1}} \ty\). By the definition of
  \(\opModality^{\modeI}_{\mode}\), this type is definitionally equal
  to
  \(\opModality^{\mode_{2}}_{\mode_{0}} \ty \to
  \opModality^{\mode_{1}}_{\mode_{0}}
  \opModality^{\mode_{2}}_{\mode_{1}} \ty\). The family of functions
  \(\unitModality^{\mode_{0}; \mode_{2}}_{\mode_{1}}\) is natural in
  the sense that for any \(\ty, \tyI : \univ_{\mode_{2}}\) and
  \(\map : \ty \to \tyI\), we have a homotopy filling the following
  square.
  \begin{equation*}
    \begin{tikzcd}
      \opModality^{\mode_{2}}_{\mode_{0}} \ty
      \arrow[r, "\unitModality^{\mode_{0}; \mode_{2}}_{\mode_{1}}(\ty)"]
      \arrow[d, "\opModality^{\mode_{2}}_{\mode_{0}} \map"'] &
      [6ex]
      \opModality^{\mode_{1}}_{\mode_{0}} \opModality^{\mode_{2}}_{\mode_{1}} \ty
      \arrow[d, "\opModality^{\mode_{1}}_{\mode_{0}} \opModality^{\mode_{2}}_{\mode_{1}} \map"] \\
      \opModality^{\mode_{2}}_{\mode_{0}} \tyI
      \arrow[r, "\unitModality^{\mode_{0}; \mode_{2}}_{\mode_{1}}(\tyI)"'] &
      \opModality^{\mode_{1}}_{\mode_{0}} \opModality^{\mode_{2}}_{\mode_{1}} \tyI
    \end{tikzcd}
  \end{equation*}
\end{construction}

Let \(\mode_{0}, \mode_{1}, \mode_{2}, \mode_{3}\) be {\acrLAMs}. By
naturality, the following diagram commutes.
\begin{equation*}
  \begin{tikzcd}
    \opModality^{\mode_{3}}_{\mode_{0}}
    \arrow[r, "\unitModality^{\mode_{0}; \mode_{3}}_{\mode_{1}}"]
    \arrow[d, "\unitModality^{\mode_{0}; \mode_{3}}_{\mode_{2}}"'] &
    [8ex]
    \opModality^{\mode_{1}}_{\mode_{0}} \opModality^{\mode_{3}}_{\mode_{1}}
    \arrow[d, "\opModality^{\mode_{1}}_{\mode_{0}} \unitModality^{\mode_{1}; \mode_{3}}_{\mode_{2}}"] \\
    [2ex]
    \opModality^{\mode_{2}}_{\mode_{0}} \opModality^{\mode_{3}}_{\mode_{2}}
    \arrow[r, "\unitModality^{\mode_{0}; \mode_{2}}_{\mode_{1}} \opModality^{\mode_{3}}_{\mode_{2}}"'] &
    \opModality^{\mode_{1}}_{\mode_{0}} \opModality^{\mode_{2}}_{\mode_{1}} \opModality^{\mode_{3}}_{\mode_{2}}
  \end{tikzcd}
\end{equation*}
For more than four {\acrLAMs}, higher coherence laws are also
satisfied. Hence, a tuple
\begin{math}
  (\mode_{0}, \dots, \mode_{\nat})
\end{math}
of {\acrLAMs} such that
\(\mode_{\idx} \le {}^{\orthMark}\mode_{\idxI}\) for all
\(\idx < \idxI\) encodes an \(\nat\)-simplex with vertices
\(\univ_{\mode_{\idx}}\), edges
\(\opModality^{\mode_{\idxI}}_{\mode_{\idx}} : \univ_{\mode_{\idxI}}
\to \univ_{\mode_{\idx}}\) for \(\idx < \idxI\), triangles
\begin{equation*}
  \begin{tikzcd}
    \univ_{\mode_{\idx}} &
    [3ex] & [3ex]
    \univ_{\mode_{\idxII}}
    \arrow[dl, "\opModality^{\mode_{\idxII}}_{\mode_{\idxI}}"]
    \arrow[ll, "\opModality^{\mode_{\idxII}}_{\mode_{\idx}}"',
    "\phantom{a}"{name = a0}]
    \arrow[from = a0, dl, Rightarrow,
    "\unitModality^{\mode_{\idx};
      \mode_{\idxII}}_{\mode_{\idxI}}"{near start},
    end anchor = {[yshift = 1ex]}] \\
    & \univ_{\mode_{\idxI}}
    \arrow[ul, "\opModality^{\mode_{\idxI}}_{\mode_{\idx}}"]
  \end{tikzcd}
\end{equation*}
for \(\idx < \idxI < \idxII\), and higher homotopies.

Shapes other than simplices are expressed by postulating invertibility
of some of
\(\unitModality^{\mode_{\idx};
  \mode_{\idxII}}_{\mode_{\idxI}}\)'s. For example, let
\(\mode_{0}, \mode_{1}, \mode_{2}, \mode_{3}\) be {\acrLAMs} and
suppose that \(\mode_{\idx} \le {}^{\orthMark}\mode_{\idxI}\) for all
\(\idx < \idxI\), that \(\mode_{2} \le {}^{\orthMark}\mode_{1}\), and
that \(\unitModality^{\mode_{0}; \mode_{3}}_{\mode_{1}}\) is
invertible. We have a diagram
\begin{equation*}
  \begin{tikzcd}[column sep = 14ex, row sep = 8ex]
    & \univ_{\mode_{1}}
    \arrow[dl, "\opModality^{\mode_{1}}_{\mode_{0}}"'] \\
    \univ_{\mode_{0}} & &
    \univ_{\mode_{3}}
    \arrow[ul, "\opModality^{\mode_{3}}_{\mode_{1}}"']
    \arrow[dl, "\opModality^{\mode_{3}}_{\mode_{2}}"]
    \arrow[ll, "\opModality^{\mode_{3}}_{\mode_{0}}"{description},
    "\phantom{a}"'{name = a0}, "\phantom{a}"{name = a1}]
    \arrow[from = a0, to = ul, To,
    "\unitModality^{\mode_{0}; \mode_{3}}_{\mode_{1}}"',
    "\simeq",
    start anchor = {[yshift = 1ex]},
    end anchor = {[yshift = -1ex]}]
    \arrow[from = a1, to = dl, To,
    "\unitModality^{\mode_{0}; \mode_{3}}_{\mode_{2}}",
    start anchor = {[yshift = -1ex]},
    end anchor = {[yshift = 1ex]}] \\
    & \univ_{\mode_{2}}
    \arrow[ul, "\opModality^{\mode_{2}}_{\mode_{0}}"]
  \end{tikzcd}
\end{equation*}
which is equivalent to a diagram of the form
\begin{equation*}
  \begin{tikzcd}
    & \univ_{\mode_{1}}
    \arrow[dl, "\opModality^{\mode_{1}}_{\mode_{0}}"']
    \arrow[dd, To,
    start anchor = {[yshift = -2ex]},
    end anchor = {[yshift = 2ex]}] \\
    \univ_{\mode_{0}} & &
    \univ_{\mode_{3}}.
    \arrow[ul, "\opModality^{\mode_{3}}_{\mode_{1}}"']
    \arrow[dl, "\opModality^{\mode_{3}}_{\mode_{2}}"] \\
    & \univ_{\mode_{2}}
    \arrow[ul, "\opModality^{\mode_{2}}_{\mode_{0}}"]
  \end{tikzcd}
\end{equation*}
\noindent 
We cannot, however, naively postulate some properties of the functors
\(\opModality^{\modeI}_{\mode}\)'s such as conservativity, fullness,
faithfulness, adjointness, and invertibility. This is because the
internal statements of these conditions are too strong due to
stability under substitution, and indeed some ``no-go'' theorems on
internalizing properties of functors are known \parencite[Theorem
5.1]{licata2018internal}\parencite[Theorem
4.1]{shulman2018brouwer}.

\begin{remark}
  It is \emph{possible} to postulate arbitrary properties of
  \(\opModality^{\mode_{\idxI}}_{\mode_{\idx}}\)'s in the following
  way. We first postulate a ``base'' {\acrLAM} \(\modeBase\) and
  assume \(\modeBase \le {}^{\orthMark} \mode_{\idx}\) for all
  \(\idx\). The universe \(\univ_{\modeBase}\) is intended to be
  interpreted as the \((\infty, 1)\)-category of spaces, so statements in
  \(\univ_{\modeBase}\) will correspond to external statements. Since
  \(\opModality_{\modeBase} : \univ \to \univ_{\modeBase}\) preserves
  finite limits, it takes \((\infty, 1)\)-categories to
  \((\infty, 1)\)-categories and functors to functors. We can then
  postulate any property on the induced functor
  \(\opModality_{\modeBase} \univ_{\mode_{\idxI}} \to
  \opModality_{\modeBase} \univ_{\mode_{\idx}}\). In fact, cohesive
  homotopy type theory \parencite{schreiber2012quantum} was first
  formulated in a similar fashion where the \(\sharp\) modality plays the
  role of \(\modeBase\). However, since we only know that
  \(\opModality_{\modeBase} \univ_{\mode_{\idx}}\) is an
  \((\infty, 1)\)-category \emph{externally}, this approach is not so
  convenient to work with especially for formalization in proof
  assistants. For this and some other reasons, the newer version of
  cohesive homotopy type theory \parencite{shulman2018brouwer} is a
  proper extension of homotopy type theory. Nevertheless, this
  adding-base approach is attractive since it keeps type theory simple
  and works for any kind of diagram.
\end{remark}

\subsection{Mode sketches}
\label{sec:mode-sketches-1}

We introduce \emph{mode sketches} as shapes of diagrams definable by
the methodology explained in \cref{sec:modes}.

\begin{definition}
  \label{def:mode-sketch}
  A \emph{mode sketch} \(\modesketch\) consists of the following data:
  \begin{itemize}
  \item a decidable finite poset \(\idxModesketch_{\modesketch}\);
  \item a subset \(\triModesketch_{\modesketch}\) of triangles in
    \(\idxModesketch_{\modesketch}\).
  \end{itemize}
  Here, by a \emph{decidable} poset we mean a poset whose ordering
  relation \(\le\) is decidable. A type is \emph{finite} if it is merely
  equivalent to the coproduct of \(\nat\) copies of \(\ilUnit\) for
  some \(\nat : \ilNat\) \parencite[Definition
  16.3.1]{rijke2022introduction-arxiv}. The identity type on a finite
  type is decidable \parencite[Remark
  16.3.2]{rijke2022introduction-arxiv}. The strict ordering relation
  \(\idx < \idxI\) defined as
  \((\idx \le \idxI) \land (\idx \neq \idxI)\) is also decidable. By a
  \emph{triangle} in \(\idxModesketch_{\modesketch}\) we mean an
  ordered triple \((\idx_{0} < \idx_{1} < \idx_{2})\) of elements of
  \(\idxModesketch_{\modesketch}\). A triangle in
  \(\triModesketch_{\modesketch}\) is called \emph{thin}.
\end{definition}

\begin{remark}
  The definition of mode sketches also makes sense in the
  metatheory. Every mode sketch \(\modesketch\) in the metatheory can
  be encoded in type theory since it is finite.
\end{remark}

Let \(\modesketch\) be a mode sketch and
\(\mode : \modesketch \to \ilLexAcc\) a function. We consider the
following axioms.

\begin{axiom}
  \label{axm:mode-sketch-disjoint}
  \(\mode(\idx) \le {}^{\orthMark} \mode(\idxI)\) for any \(\idxI
  \not\le \idx\) in \(\modesketch\).
\end{axiom}

\begin{axiom}
  \label{axm:mode-sketch-invertible}
  For any thin triangle \((\idx_{0} < \idx_{1} < \idx_{2})\) in
  \(\modesketch\), the natural transformation
  \(\unitModality^{\mode(\idx_{0}); \mode(\idx_{2})}_{\mode(\idx_{1})}
  : \opModality^{\mode(\idx_{2})}_{\mode(\idx_{0})} \To
  \opModality^{\mode(\idx_{1})}_{\mode(\idx_{0})}
  \opModality^{\mode(\idx_{2})}_{\mode(\idx_{1})}\) is invertible.
\end{axiom}

\begin{remark}
  Assuming \cref{axm:mode-sketch-disjoint}, if \(\idx < \idxI\), then
  \(\mode(\idx) \le {}^{\orthMark} \mode(\idxI)\). The converse is not
  true: when neither \(\idx \le \idxI\) nor \(\idxI \le \idx\), we still
  get \(\mode(\idx) \le {}^{\orthMark} \mode(\idxI)\).
\end{remark}

\Cref{axm:mode-sketch-disjoint,axm:mode-sketch-invertible} are
motivated by the observation made in \cref{sec:modes}. That is, when
\(\idxI \not\leq \idx\), the functor in the direction
\(\univ_{\mode(\idx)} \to \univ_{\mode(\idxI)}\) is cut off.

\begin{remark}
  \label{rem:2-cat-mode-sketch}
  A mode sketch \(\modesketch\) is regarded as a presentation of an
  \((\infty, 2)\)-category \(\realize{\modesketch}\). The strict ordering
  relation generates \(1\)-cells \((\idx < \idxI) : \idx \to \idxI\),
  and the triangles \((\idx < \idxI < \idxII)\) generate \(2\)-cells
  in the direction
  
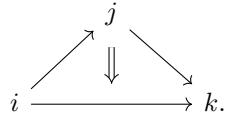
\begin{equation*}
    \begin{tikzcd}
      & \idxI
      \arrow[dr] \\
      \idx
      \arrow[ur]
      \arrow[rr, "\phantom{a}"{name = a0}]
      &
      \arrow[from = u, to = a0, To, start anchor = {[yshift = -1ex]}] &
      \idxII.
    \end{tikzcd}
  \end{equation*}
  When the triangle is thin, the corresponding \(2\)-cell is made
  invertible. Longer chains
  \((\idx_{0} < \idx_{1} < \dots < \idx_{\nat})\) present homotopies
  filling certain diagrams. A formal account is given in
  \cref{cst:oo-2-cat-from-mode-sketch}. A function
  \(\mode : \modesketch \to \ilLexAcc\) satisfying
  \cref{axm:mode-sketch-disjoint,axm:mode-sketch-invertible} is then
  considered as a diagram of subuniverses indexed over
  \(\realize{\modesketch}^{\opMark(1, 2)}\), the
  \((\infty, 2)\)-category obtained from \(\realize{\modesketch}\) by
  reversing the directions of \(1\)-cells and \(2\)-cells.
\end{remark}

\begin{example}
  Every decidable finite poset is a mode sketch where no triangle is
  thin. The \((\infty, 2)\)-category presented by it is obtained from the
  left adjoint of the Duskin nerve \parencite{duskin2001simplicial} by
  reversing \(2\)-cells.
\end{example}

\begin{example}
  \label{exm:mode-sketch-for-functors}
  The \emph{mode sketch for functors} is drawn as
  \begin{equation*}
    \begin{tikzcd}
      0 \arrow[r] & 1.
    \end{tikzcd}
  \end{equation*}
  \Cref{axm:mode-sketch-disjoint} asserts
  \begin{math}
    \mode(0) \le {}^{\orthMark} \mode(1).
  \end{math}
  \Cref{axm:mode-sketch-invertible} is empty since there is no
  triangle. Thus, we get the following diagram.
  \begin{equation*}
    \begin{tikzcd}
      \univ_{\mode(0)} &
      [2ex]
      \univ_{\mode(1)}
      \arrow[l, "\opModality^{\mode(1)}_{\mode(0)}"']
    \end{tikzcd}
  \end{equation*}
\end{example}

\begin{example}
  The \emph{mode sketch for triangles} is drawn as
  \begin{equation*}
    \begin{tikzcd}
      0
      \arrow[rr, ""'{name = a0}]
      \arrow[dr] &
      \arrow[from = a0, to = d, phantom, "\simeq"{description}] &
      2 \\
      & 1
      \arrow[ur]
    \end{tikzcd}
  \end{equation*}
  where ``\(\simeq\)'' indicates that the triangle is
  thin. \Cref{axm:mode-sketch-disjoint} asserts
  \(\mode(0) \le {}^{\orthMark} \mode(1)\),
  \(\mode(0) \le {}^{\orthMark} \mode(2)\), and
  \(\mode(1) \le {}^{\orthMark}
  \mode(2)\). \Cref{axm:mode-sketch-invertible} asserts that
  \(\unitModality^{\mode(0); \mode(2)}_{\mode(1)}\) is
  invertible. Thus, we have the following commutative triangle.
  \begin{equation*}
    \begin{tikzcd}
      \univ_{\mode(0)} & &
      \univ_{\mode(2)}
      \arrow[ll, "\opModality^{\mode(2)}_{\mode(0)}"']
      \arrow[ld, "\opModality^{\mode(2)}_{\mode(1)}"] \\
      & \univ_{\mode(1)}
      \arrow[ul, "\opModality^{\mode(1)}_{\mode(0)}"]
    \end{tikzcd}
  \end{equation*}
\end{example}

\subsection{Intended models, internally}
\label{sec:intend-models-intern}

Let \(\modesketch\) be a mode sketch. In this subsection, we see,
internally to type theory, what kind of an \(\infty\)-logos is a model of
\(\modesketch\). Here, by a model of \(\modesketch\) we mean an
\(\infty\)-logos that admits an interpretation of a postulated function
\(\mode : \modesketch \to \ilLexAcc\) satisfying
\cref{axm:mode-sketch-disjoint,axm:mode-sketch-invertible} and the
following additional axiom.

\begin{axiom}
  \label{axm:mode-sketch-top}
  The top modality is the canonical join
  \begin{math}
    \biglor_{\modesketch} \mode.
  \end{math}
\end{axiom}

\Cref{axm:mode-sketch-top} roughly asserts that the whole universe
\(\univ\) is reconstructed from the subuniverses
\(\univ_{\mode(\idx)}\)'s. This is meant to exclude models other than
intended models.

\begin{example}
  \label{exm:intended-model-mode-sketch-for-functors}
  Consider the case when \(\modesketch\) is the mode sketch for
  functors
  (\cref{exm:mode-sketch-for-functors}). \Cref{axm:mode-sketch-top}
  asserts \(\modeTop = \mode(0) \lor \mode(1)\). The equivalence
  \(\univ \simeq \ilexists_{\ty : \univ_{\mode(0)}} \ilexists_{\tyI :
    \univ_{\mode(1)}} \ty \to \opModality^{\mode(1)}_{\mode(0)} \tyI\)
  (\cref{prop:join-strongly-disjoint}) suggests that \(\univ\) is the
  so-called \emph{Artin gluing} for the functor
  \(\opModality^{\mode(1)}_{\mode(0)} : \univ_{\mode(1)} \to
  \univ_{\mode(0)}\). Therefore, our intended models of
  \(\modesketch\) are \(\infty\)-logoses obtained by the Artin gluing.
\end{example}

A generalization of the Artin gluing is \emph{oplax limits}. In the
setting of \cref{exm:intended-model-mode-sketch-for-functors},
\(\univ\) fits into the following \emph{universal oplax cone} over the
diagram
\(\univ_{\mode(0)} \xleftarrow{\opModality^{\mode(1)}_{\mode(0)}}
\univ_{\mode(1)}\).
\begin{equation}
  \labelX[diagram]{eq:gluing-as-oplax-limit}
  \begin{tikzcd}
    & \univ
    \arrow[dl, "\phantom{a}"{name = a0}]
    \arrow[dr, "\phantom{a}"'{name = a1}]
    \arrow[from = a0, to = a1, To] \\
    \univ_{\mode(0)} & &
    \univ_{\mode(1)}
    \arrow[ll, "\opModality^{\mode(1)}_{\mode(0)}"]
  \end{tikzcd}
\end{equation}
An oplax cone over a diagram is a kind of cone but every triangle
formed by two projections and a functor in the diagram is only filled
by a not necessarily invertible natural transformation in the
direction of \cref{eq:gluing-as-oplax-limit}. The universal oplax cone
or oplax limit is the terminal object in the \((\infty, 1)\)-category of
oplax cones.

\begin{example}
  \label{exm:intended-model-mode-sketch-for-trans}
  Consider the case when \(\modesketch\) is the mode sketch
  \(\{0 \to 1 \to 2\}\) with no thin triangle. \Cref{axm:mode-sketch-top}
  asserts \(\modeTop = \mode(0) \lor \mode(1) \lor
  \mode(2)\). Iterating \cref{prop:join-strongly-disjoint}, we see
  that every type \(\ty : \univ\) is fractured into
  \(\ty_{0} : \univ_{\mode(0)}\), \(\ty_{1} : \univ_{\mode(1)}\),
  \(\ty_{2} : \univ_{\mode(2)}\),
  \(\map_{01} : \ty_{0} \to \opModality^{\mode(1)}_{\mode(0)} \ty_{1}\),
  \(\map_{02} : \ty_{0} \to \opModality^{\mode(2)}_{\mode(0)} \ty_{2}\),
  \(\map_{12} : \ty_{1} \to \opModality^{\mode(2)}_{\mode(1)} \ty_{2}\),
  and
  \(\pth_{012} : \opModality^{\mode(1)}_{\mode(0)} \map_{12} \comp
  \map_{01} = \unitModality^{\mode(0); \mode(2)}_{\mode(1)}(\ty_{2})
  \comp \map_{02}\). Indeed, we have
  \begin{EqReasoning}
    \begin{align*}
      & \term{\univ} \\
      \simeq & \by{\cref{prop:join-strongly-disjoint} for \(\mode(0)\) and \(\mode(1) \lor \mode(2)\)} \\
      & \term{\ilexists_{\ty_{0} : \univ_{\mode(0)}}
        \ilexists_{\ty_{12} : \univ_{\mode(1) \lor \mode(2)}}
        \ty_{0} \to \opModality_{\mode_{0}} \ty_{12}} \\
      \simeq & \by{\cref{prop:join-strongly-disjoint} for \(\mode(1)\) and \(\mode(2)\)} \\
      & \term{\ilexists_{\ty_{0} : \univ_{\mode(0)}}
        \ilexists_{\ty_{1} : \univ_{\mode(1)}}
        \ilexists_{\ty_{2} : \univ_{\mode(2)}}
        \ilexists_{\map_{12} : \ty_{1} \to \opModality_{\mode(1)} \ty_{2}}
        \ty_{0} \to \opModality_{\mode(0)}(\ty_{1} \times_{\opModality_{\mode(1)} \ty_{2}} \ty_{2})}
    \end{align*}
  \end{EqReasoning}
  where the pullback is taken for
  \(\map_{12} : \ty_{1} \to \opModality_{\mode(1)} \ty_{2}\) and
  \(\unitModality_{\mode(1)}(\ty_{2}) : \ty_{2} \to
  \opModality_{\mode(1)} \ty_{2}\). Since \(\opModality_{\mode(0)}\)
  preserves pullbacks, the component
  \(\ty_{0} \to \opModality_{\mode(0)}(\ty_{1} \times_{\opModality_{\mode(1)}
    \ty_{2}} \ty_{2})\) corresponds to the components \(\map_{01}\),
  \(\map_{02}\), and \(\pth_{012}\). Then \(\univ\) is the oplax limit
  of the diagram
  \begin{equation}
    \labelX[diagram]{eq:canonical-triangle}
    \begin{tikzcd}[column sep = 12ex]
      \univ_{\mode(0)} & &
      \univ_{\mode(2)}.
      \arrow[ll, "\opModality^{\mode(2)}_{\mode(0)}"', "\phantom{a}"{name = a0}]
      \arrow[dl, "\opModality^{\mode(2)}_{\mode(1)}"]\\
      & \univ_{\mode(1)}
      \arrow[ul, "\opModality^{\mode(1)}_{\mode(0)}"]
      \arrow[from = a0, To, "\unitModality^{\mode(0); \mode(2)}_{\mode(1)}",
      end anchor = {[yshift = 1ex]}]
    \end{tikzcd}
  \end{equation}
  That is, we have projections \(\univ \to \univ_{\mode(\idx)}\) for
  all \(\idx\), natural transformations
  \begin{equation*}
    \begin{tikzcd}
      & \univ
      \arrow[dl, "\phantom{a}"{name = a0}]
      \arrow[dr, "\phantom{a}"'{name = a1}]
      \arrow[from = a0, to = a1, To] \\
      \univ_{\mode(\idx)} & &
      \univ_{\mode(\idxI)}
      \arrow[ll, "\opModality^{\mode(\idxI)}_{\mode(\idx)}"]
    \end{tikzcd}
  \end{equation*}
  for all \(\idx < \idxI\), and a homotopy
  \begin{equation*}
    \begin{tikzcd}
      & \univ
      \arrow[dl, "\phantom{a}"{name = a0, near end}]
      \arrow[dr, "\phantom{a}"'{name = a1, near end}]
      \arrow[dd, "\phantom{a}"'{name = a2}, "\phantom{a}"{name = a3}]
      \arrow[from = a0, to = a2, To]
      \arrow[from = a3, to = a1, To] \\
      \univ_{\mode(0)} & &
      \univ_{\mode(2)}
      \arrow[dl, "\opModality^{\mode(2)}_{\mode(1)}"] \\
      & \univ_{\mode(1)}
      \arrow[ul, "\opModality^{\mode(1)}_{\mode(0)}"]
    \end{tikzcd}
    =
    \begin{tikzcd}[column sep = 10ex]
      & \univ
      \arrow[dl, "\phantom{a}"{name = a0}]
      \arrow[dr, "\phantom{a}"'{name = a1}]
      \arrow[from = a0, to = a1, To] \\
      \univ_{\mode(0)} & &
      \univ_{\mode(2)},
      \arrow[ll, "\opModality^{\mode(2)}_{\mode(0)}"{description}, "\phantom{a}"{name = a2}]
      \arrow[dl, "\opModality^{\mode(2)}_{\mode(1)}"] \\
      & \univ_{\mode(1)}
      \arrow[ul, "\opModality^{\mode(1)}_{\mode(0)}"]
      \arrow[from = a2, To, "\unitModality^{\mode(0); \mode(2)}_{\mode(1)}",
      end anchor = {[yshift = 1ex]}]
    \end{tikzcd}
  \end{equation*}
  and these data form a universal oplax cone over
  \cref{eq:canonical-triangle}.

  Let us make the triangle \((0 < 1 < 2)\) thin so that the natural
  transformation \labelcref{eq:canonical-triangle} becomes
  invertible. In this setting, \(\univ\) is still the oplax limit of
  \cref{eq:canonical-triangle}, but the presentation can be simplified
  since the type of data \((\map_{02}, \pth_{012})\) is contractible.
\end{example}

For a general mode sketch \(\modesketch\), we apply
\cref{prop:join-strongly-disjoint} for a minimal element
\(\mode(\idx_{0})\) and the rest
\(\biglor_{\idx : \modesketch \setminus \idx_{0}} \mode(\idx)\) and repeat
this for \(\modesketch \setminus \idx_{0}\) to fracture types into modal
types. \Cref{exm:intended-model-mode-sketch-for-functors,exm:intended-model-mode-sketch-for-trans}
suggest that \(\univ\) is the oplax limit of the diagram formed by
\(\univ_{\mode(\idx)}\)'s explained in
\cref{rem:2-cat-mode-sketch}. Thus, our intended models of
\(\modesketch\) are oplax limits of \(\infty\)-logoses indexed over the
\((\infty, 2)\)-category presented by \(\modesketch\). The formal account
of this is described in \cref{sec:semant-mode-sketch}.

\section{Mode sketches and synthetic Tait computability}
\label{sec:mode-sketch-synth}

We give an alternative set of axioms for mode sketches and exhibit a
connection between mode sketches and \emph{synthetic Tait
  computability} of Sterling \parencite{sterling2021thesis}. The core axiom of
synthetic Tait computability is to postulate a proposition. The
proposition induces the open and closed modalities, and then every
type is fractured into an open type equipped with a closed type family
and behaves like a \emph{logical relation}. In this story, the open
and closed modalities seem more essential than the postulated
proposition, so we aim to formulate synthetic Tait computability
purely in terms of modalities. We work in homotopy type theory.

\subsection{Alternative mode sketch axioms}
\label{sec:altern-mode-sketch}

The \(\infty\)-logoses obtained by the Artin gluing can be characterized as
\(\infty\)-logoses equipped with a subterminal object
(cf.\@ \parencite[A4.5.6]{johnstone2002sketches}). We generalize this
from the Artin gluing to oplax limits indexed by mode sketches,
internally to type theory: the type of functions
\(\modesketch \to \ilLexAcc\) satisfying
\cref{axm:mode-sketch-disjoint,axm:mode-sketch-top} is equivalent to
the type of lattice morphisms \(\ilCosieve(\modesketch) \to \ilProp\)
(\cref{thm:equiv-two-mode-sketch-axioms}).

\begin{definition}
  A \emph{cosieve} on a decidable poset \(\idxsh\) is an upward-closed
  decidable subset of it. Let \(\ilCosieve(\idxsh)\) denote the poset
  of cosieves on \(\idxsh\) ordered by inclusion. Note that cosieves
  are closed under finite meets and joins, so \(\ilCosieve(\idxsh)\)
  is a lattice.
\end{definition}

\begin{notation}
  For \(\idx : \modesketch\), let \((\idx \downarrow \modesketch)\) denote the cosieve
  \(\{\idxI : \modesketch \mid \idx \le \idxI\}\) and
  \(\boundary (\idx \downarrow \modesketch)\) the cosieve \((\idx \downarrow \modesketch) \setminus \{\idx\}\).
\end{notation}

\begin{construction}
  \label{cst:mode-sketch-prop-to-mode}
  Let \(\propo : \ilCosieve(\modesketch) \to \ilProp\) be a function. We
  define a function
  \begin{math}
    \modeFromProp_{\propo} : \modesketch \to \ilLexAcc
  \end{math}
  by
  \begin{math}
    \modeFromProp_{\propo}(\idx) \defeq
    \modeOpen(\propo(\idx \downarrow \modesketch)) \land
    \modeClosed(\propo(\boundary (\idx \downarrow \modesketch))).
  \end{math}
\end{construction}

\begin{theorem}
  \label{thm:equiv-two-mode-sketch-axioms}
  \Cref{cst:mode-sketch-prop-to-mode} is restricted to an equivalence
  between the following types:
  \begin{enumerate}
  \item \label{item:mode-sketch-axiom-1} the type of lattice morphisms
    \(\propo : \ilCosieve(\modesketch) \to \ilProp\);
  \item \label{item:mode-sketch-axiom-2} the type of functions
    \(\mode : \modesketch \to \ilLexAcc\) satisfying
    \cref{axm:mode-sketch-disjoint,axm:mode-sketch-top}.
  \end{enumerate}
  We write
  \(\propCanonical_{\mode} : \ilCosieve(\modesketch) \to \ilProp\) for
  the lattice morphism corresponding to \(\mode\).
\end{theorem}
\noindent 
Before giving a proof of \cref{thm:equiv-two-mode-sketch-axioms}, let
us relate \cref{thm:equiv-two-mode-sketch-axioms} to \emph{synthetic
  Tait computability}
\parencite{sterling2021logical,sterling2021thesis,sterling2022sheaf}. The
core axiom of synthetic Tait computability is to postulate some
propositions. One can work with those propositions directly but also
with the induced open and closed
modalities. \Cref{thm:equiv-two-mode-sketch-axioms} says that
synthetic Tait computability can, in fact, be formulated completely in
terms of modalities. The simplest version of synthetic Tait
computability postulates a single proposition. The corresponding mode
sketch is \(\{0 \to 1\}\) as follows.

\begin{example}
  \label{exm:mode-sketch-and-stc}
  Let \(\modesketch\) be the mode sketch for functors
  (\cref{exm:mode-sketch-for-functors}). Then
  \begin{math}
    \ilCosieve(\modesketch) =
    \{\{\}, \{1\}, \{0, 1\}\}
  \end{math}
  is the free lattice generated by the single element
  \(\{1\}\). We thus have
  \begin{math}
    \{\text{lattice morphisms \(\ilCosieve(\modesketch) \to \ilProp\)}\}
    \simeq \ilProp.
  \end{math}
\end{example}

The rest of this subsection is devoted to the proof of
\cref{thm:equiv-two-mode-sketch-axioms}. We first show that
\cref{cst:mode-sketch-prop-to-mode} is restricted to a function
\begin{math}
  \labelcref{item:mode-sketch-axiom-1}
  \to \labelcref{item:mode-sketch-axiom-2}.
\end{math}

\begin{proposition}
  \label{prop:mode-sketch-axioms-comparison}
  Let \(\propo : \ilCosieve(\modesketch) \to \ilProp\) be a function.
  \begin{enumerate}
  \item If \(\propo\) preserves binary meets, then
    \(\modeFromProp_{\propo}\) satisfies
    \cref{axm:mode-sketch-disjoint}.
  \item If \(\propo\) preserves top elements and finite joins, then
    \(\modeFromProp_{\propo}\) satisfies \cref{axm:mode-sketch-top}.
  \end{enumerate}
\end{proposition}
\noindent 
We prepare a couple of lemmas.

\begin{lemma}
  \label{lem:join-mode-from-prop}
  If a function \(\propo : \ilCosieve(\modesketch) \to \ilProp\)
  preserves finite joins, then
  \(\modeOpen(\propo(\sieve)) = \biglor_{\sieve}
  \modeFromProp_{\propo}\) for any cosieve \(\sieve \subset \modesketch\).
\end{lemma}
\begin{proof}
  By induction on the size of \(\sieve\). By assumption and by
  \cref{prop:open-modality-embedding}, we have
  \(\modeOpen(\propo(\sieve)) = \biglor_{\idx : \sieve}
  \modeOpen(\propo(\idx \downarrow \modesketch))\). Thus, it is enough to show the case
  when \(\sieve\) is of the form \((\idx \downarrow \modesketch)\). By
  \cref{prop:open-closed-decomposition}, we have
  \begin{math}
    \modeOpen(\propo(\idx \downarrow \modesketch))
    = (\modeOpen(\propo(\idx \downarrow \modesketch)) \land \modeClosed(\propo(\boundary (\idx \downarrow \modesketch)))) \lor (\modeOpen(\propo(\idx \downarrow \modesketch)) \land \modeOpen(\propo(\boundary (\idx \downarrow \modesketch))))
    = \modeFromProp_{\propo}(\idx) \lor \modeOpen(\propo(\boundary (\idx \downarrow \modesketch))).
  \end{math}
  Then apply the induction hypothesis for \(\boundary (\idx \downarrow \modesketch)\).
\end{proof}

\begin{lemma}
  \label{prop:meet-connected-cancel-1}
  Let \(\mode\) and \(\modeI\) be {\acrLAMs}. If
  \(\opModality_{\modeI}\) preserves \(\mode\)-modal types, then
  \(\mode \land {}^{\orthMark}(\modeI \land \mode) = \mode \land
  {}^{\orthMark} \modeI\).
\end{lemma}
\begin{proof}
  By \cref{prop:meet-preserves-modal}.
\end{proof}

\begin{proof}[Proof of \cref{prop:mode-sketch-axioms-comparison}]
  If \(\propo\) preserves binary meets, then for any
  \(\idxI \not\le \idx\) in \(\modesketch\),
  \begin{EqReasoning}
    \begin{align*}
      & \term{\modeFromProp_{\propo}(\idx)} \\
      = & \by{definition} \\
      & \term{\modeOpen(\propo(\idx \downarrow \modesketch)) \land {}^{\orthMark} \modeOpen(\propo(\boundary (\idx \downarrow \modesketch)))} \\
      \le & \by{\((\idxI \downarrow \modesketch) \cap (\idx \downarrow \modesketch) \subset \boundary (\idx \downarrow \modesketch)\) as \(\idxI \not\le \idx\)} \\
      & \term{\modeOpen(\propo(\idx \downarrow \modesketch)) \land {}^{\orthMark} \modeOpen(\propo((\idxI \downarrow \modesketch) \cap (\idx \downarrow \modesketch)))} \\
      = & \by{\cref{prop:open-modality-embedding}, and \(\propo\) preserves binary meets} \\
      & \term{\modeOpen(\propo(\idx \downarrow \modesketch)) \land {}^{\orthMark} (\modeOpen(\propo(\idxI \downarrow \modesketch)) \land \modeOpen(\propo(\idx \downarrow \modesketch)))} \\
      = & \by{\cref{prop:meet-connected-cancel-1,prop:meet-open-closed}} \\
      & \term{\modeOpen(\propo(\idx \downarrow \modesketch)) \land {}^{\orthMark} \modeOpen(\propo(\idxI \downarrow \modesketch))} \\
      \le & \by{definition} \\
      & \term{{}^{\orthMark} \modeFromProp_{\propo}(\idxI)},
    \end{align*}
  \end{EqReasoning}
  and thus \cref{axm:mode-sketch-disjoint} is satisfied. If \(\propo\)
  preserves top elements and finite joins, then
  \(\modeFromProp_{\propo}\) satisfies \cref{axm:mode-sketch-top} by
  \cref{lem:join-mode-from-prop}.
\end{proof}
\noindent 
We then construct the inverse function
\begin{math}
  \labelcref{item:mode-sketch-axiom-2}
  \to \labelcref{item:mode-sketch-axiom-1}.
\end{math}
The key observation is that canonical joins
of \(\mode(\idx)\)'s exist and are well-behaved under
\cref{axm:mode-sketch-disjoint}.

\begin{proposition}
  \label{prop:mode-sketch-canonical-join}
  If a function \(\mode : \modesketch \to \ilLexAcc\) satisfies
  \cref{axm:mode-sketch-disjoint}, then the canonical join
  \begin{math}
    \biglor_{\sieve} \mode
  \end{math}
  exists for any decidable subset \(\sieve \subset \modesketch\).
\end{proposition}
\begin{proof}
  By induction on the size of \(\sieve\). If \(\sieve\) is empty, then
  \(\biglor_{\emptyset} \mode\) is the bottom modality. Suppose that
  \(\sieve\) is non-empty. Since \(\modesketch\) is finite, there is
  an element \(\idx_{0}\) minimal in \(\sieve\). Then
  \(\sieve \setminus \{\idx_{0}\}\) admits a canonical join by the induction
  hypothesis. Since \(\idx_{0}\) is minimal,
  \(\mode(\idx_{0}) \le {}^{\orthMark}\mode(\idx)\) for any
  \(\idx : \sieve \setminus \{\idx_{0}\}\) by \cref{axm:mode-sketch-disjoint},
  and thus
  \(\mode(\idx_{0}) \le {}^{\orthMark} (\biglor_{\sieve \setminus \{\idx_{0}\}}
  \mode)\). Then we have the canonical join
  \begin{math}
    \biglor_{\sieve} \mode \defeq
    \mode(\idx_{0}) \lor (\biglor_{(\sieve \setminus \{\idx_{0}\})} \mode)
  \end{math}
  by \cref{prop:join-strongly-disjoint}.
\end{proof}

\begin{lemma}
  \label{cor:join-strongly-disjoint-in-susigma}
  Let \(\mode_{0}\), \(\mode_{1}\), and \(\mode_{2}\) be {\acrLAMs}
  such that \(\mode_{\idx} \le {}^{\orthMark} \mode_{\idxI}\) for any
  \(\idx < \idxI\). Then
  \(\mode_{0} \lor \mode_{1} \le {}^{\orthMark} \mode_{2}\).
\end{lemma}
\begin{proof}
  Let \(\ty\) be a \((\mode_{0} \lor \mode_{1})\)-modal type. By
  \cref{prop:join-strongly-disjoint},
  \(\unitModality_{\mode_{1}}(\ty) : \ty \to \opModality_{\mode_{1}}
  \ty\) has \(\mode_{0}\)-modal fibers. Then, by assumption,
  \(\opModality_{\mode_{1}} \ty\) and the fibers of
  \(\unitModality_{\mode_{1}}(\ty)\) are made contractible by
  \(\opModality_{\mode_{2}}\). Thus, \(\opModality_{\mode_{2}} \ty\)
  is contractible.
\end{proof}

\begin{proposition}
  \label{prop:mode-sketch-partition-by-sieve}
  If a function \(\mode : \modesketch \to \ilLexAcc\) satisfies
  \cref{axm:mode-sketch-disjoint}, then
  \begin{math}
    \biglor_{\modesketch \setminus \sieve} \mode
    \le {}^{\orthMark} (\biglor_{\sieve} \mode)
  \end{math}
  for any cosieve \(\sieve \subset \modesketch\).
\end{proposition}
\begin{proof}
  Since \(\sieve\) is upward-closed, \(\idxI \not\le \idx\) for any
  \(\idx : \modesketch \setminus \sieve\) and \(\idxI : \sieve\). Thus, by
  \cref{axm:mode-sketch-disjoint},
  \(\mode(\idx) \le {}^{\orthMark} \mode(\idxI)\) for any
  \(\idx : \modesketch \setminus \sieve\) and \(\idxI : \sieve\). The claim
  follows from \cref{cor:join-strongly-disjoint-in-susigma} and the
  construction of the canonical join in
  \cref{prop:mode-sketch-canonical-join}.
\end{proof}

\begin{fact}[{\cite[Example 3.32]{rijke2020modalities}}]
  \label{prop:meet-strongly-disjoint}
  Let \(\mode\) and \(\modeI\) be {\acrLAMs}. If
  \(\mode \le {}^{\orthMark} \modeI\), then
  \(\mode \land \modeI = \modeBottom\).
\end{fact}

\begin{construction}
  Let \(\mode : \modesketch \to \ilLexAcc\) be a function satisfying
  \cref{axm:mode-sketch-disjoint,axm:mode-sketch-top}. We define a
  lattice morphism
  \begin{math}
    \propCanonical_{\mode} :
    \ilCosieve(\modesketch) \to \ilProp
  \end{math}
  by
  \begin{math}
    \propCanonical_{\mode}(\sieve) \defeq
    \opModality_{\biglor_{\modesketch \setminus \sieve} \mode} \ilEmpty.
  \end{math}
  By \cref{cor:partition-open-modality} and by
  \cref{prop:mode-sketch-partition-by-sieve},
  \(\propCanonical_{\mode}(\sieve)\) is the unique proposition such
  that
  \(\modeOpen(\propCanonical_{\mode}(\sieve)) = \biglor_{\sieve}
  \mode\). Because \(\sieve \mapsto \biglor_{\sieve} \mode\) preserves
  finite joins by definition, \(\propCanonical_{\mode}\) preserves
  finite joins. By \cref{axm:mode-sketch-top},
  \(\propCanonical_{\mode}\) preserves top elements. For preservation
  of binary meets, let \(\sieve_{1}\) and \(\sieve_{2}\) be cosieves
  on \(\modesketch\). We have to show that
  \begin{math}
    \biglor_{\sieve_{1} \cap \sieve_{2}} \mode
    = (\biglor_{\sieve_{1}} \mode) \land (\biglor_{\sieve_{2}} \mode).
  \end{math}
  Let \(\sieve_{3} \defeq \sieve_{1} \cap \sieve_{2}\), \(\sieve_{1}'
  \defeq \sieve_{1} \setminus \sieve_{3}\), and \(\sieve_{2}' \defeq
  \sieve_{2} \setminus \sieve_{3}\). By
  \cref{prop:mode-sketch-partition-by-sieve}, \(\biglor_{\sieve_{1}'}
  \mode \le
  {}^{\orthMark} (\biglor_{\sieve_{2}} \mode)\) and
  \(\biglor_{\sieve_{2}'} \mode \le {}^{\orthMark} (\biglor_{\sieve_{1}}
  \mode)\). Then,
  \begin{EqReasoning}
    \let\oldTerm\term
    \renewcommand{\term}[1]{\oldTerm{\textstyle#1}}
    \begin{align*}
      & \term{(\biglor_{\sieve_{1}} \mode) \land (\biglor_{\sieve_{2}} \mode)} \\
      = & \by{definition} \\
      & \term{((\biglor_{\sieve_{1}'} \mode) \lor (\biglor_{\sieve_{3}} \mode)) \land ((\biglor_{\sieve_{2}'} \mode) \lor (\biglor_{\sieve_{3}} \mode))} \\
      = & \by{\cref{prop:join-strongly-disjoint-distributive}} \\
      & \term{((\biglor_{\sieve_{1}'} \mode) \land (\biglor_{\sieve_{2}'} \mode)) \lor (\biglor_{\sieve_{3}} \mode)} \\
      = & \by{\cref{prop:meet-strongly-disjoint}} \\
      & \term{\biglor_{\sieve_{3}} \mode}.
    \end{align*}
  \end{EqReasoning}
\end{construction}

\begin{proof}[Proof of \cref{thm:equiv-two-mode-sketch-axioms}]
  It remains to show that the constructions \(\propo \mapsto
  \modeFromProp_{\propo}\) and \(\mode \mapsto \propCanonical_{\mode}\) are
  mutually inverses. \Cref{lem:join-mode-from-prop} implies that
  \(\propCanonical_{\modeFromProp_{\propo}} = \propo\). For the other
  identification, we have
  \belowdisplayskip=0pt
  \begin{EqReasoning}
    \let\oldTerm\term
    \renewcommand{\term}[1]{\oldTerm{\textstyle#1}}
    \begin{align*}
      & \term{\modeFromProp_{\propCanonical_{\mode}}(\idx)} \\
      = & \by{definition} \\
      & \term{\modeOpen(\propCanonical_{\mode}(\idx \downarrow \modesketch)) \land \modeClosed(\propCanonical_{\mode}(\boundary (\idx \downarrow \modesketch)))} \\
      = & \by{\cref{lem:join-mode-from-prop}} \\
      & \term{(\biglor_{(\idx \downarrow \modesketch)} \mode) \land \modeClosed(\propCanonical_{\mode}(\boundary (\idx \downarrow \modesketch)))} \\
      = & \by{\cref{prop:meet-preserves-modal-distribute-disjoint,prop:meet-open-closed}} \\
      & \term{\biglor_{\idxI : (\idx \downarrow \modesketch)} \mode(\idxI) \land \modeClosed(\propCanonical_{\mode}(\boundary (\idx \downarrow \modesketch)))} \\
      = & \by{\cref{lem:join-mode-from-prop}} \\
      & \term{\biglor_{\idxI : (\idx \downarrow \modesketch)} \mode(\idxI) \land (\bigland_{\idxII : \boundary(\idx \downarrow \modesketch)} {}^{\orthMark} \mode(\idxII))} \\
      = & \by{\(\mode(\idxII) \land {}^{\orthMark} \mode(\idxII) = \modeBottom\) for \(\idxII : \boundary(\idx \downarrow \modesketch)\) by \cref{prop:meet-strongly-disjoint}} \\
      & \term{\mode(\idx) \land (\bigland_{\idxII : \boundary(\idx \downarrow \modesketch)} {}^{\orthMark} \mode(\idxII))} \\
      = & \by{\cref{axm:mode-sketch-disjoint}} \\
      & \term{\mode(\idx)}. \qedhere 
    \end{align*}
  \end{EqReasoning}
\end{proof}

\subsection{Logical relations as types}
\label{sec:logical-relations-as}

We have seen in \cref{sec:altern-mode-sketch} that synthetic Tait
computability is reformulated in terms of {\acrLAMs}. The slogan of
synthetic Tait computability is ``logical relations as types''
\parencite{sterling2021logical}. This is also formulated purely in
terms of {\acrLAMs}.

\begin{fact}[{\cite[Theorem 3.11]{rijke2020modalities}}]
  \label{fact:universe-of-modal-types}
  For any {\acrLAM} \(\mode\), the universe of \(\mode\)-modal types
  \begin{math}
    \univ_{\mode} \defeq \{\ty : \univ \mid \ilIn_{\mode}(\ty)\}
  \end{math}
  is \(\mode\)-modal.
\end{fact}

\begin{proposition}[Fracture and gluing]
  \label{prop:fracture-and-gluing-alt}
  Let \(\mode\) and \(\modeI\) be {\acrLAMs} such that
  \(\mode \le {}^{\orthMark} \modeI\). Then we have an equivalence
  \begin{equation*}
    \univ_{\mode \lor \modeI}
    \simeq \ilexists_{\tyI : \univ_{\modeI}}\tyI \to \univ_{\mode}
  \end{equation*}
  whose right-to-left function sends a \((\tyI, \ty)\) to
  \(\ilexists_{\var : \tyI}\ty(\var)\).
\end{proposition}
\begin{proof}
  For any \(\tyI : \univ_{\modeI}\), we have
  \begin{EqReasoning}
    \begin{align*}
      & \term{\ilexists_{\ty : \univ_{\mode}}\ty \to \opModality_{\mode} \tyI} \\
      \simeq& \by{equivalence between fibrations and type families} \\
      & \term{\opModality_{\mode} \tyI \to \univ_{\mode}} \\
      \simeq& \by{\cref{fact:universe-of-modal-types}} \\
      & \term{\tyI \to \univ_{\mode}}.
    \end{align*}
  \end{EqReasoning}
  Then apply \cref{prop:join-strongly-disjoint}.
\end{proof}
\noindent 
\Cref{prop:fracture-and-gluing-alt} asserts that a type in
\(\univ_{\mode \lor \modeI}\) is a \(\modeI\)-modal type equipped with
a \(\mode\)-modal unary (proof-relevant) relation on it, so
\emph{types (in \(\univ_{\mode \lor \modeI}\)) are relations}. More
generally, for a mode sketch \(\modesketch\) and a function
\(\mode : \modesketch \to \ilLexAcc\) satisfying
\cref{axm:mode-sketch-disjoint}, types in
\(\univ_{\biglor_{\modesketch}\mode}\) are fractured into a sort of
generalized relations by iterated applications of
\cref{prop:fracture-and-gluing-alt}. Intuitively, the ordering on
\(\modesketch\) is understood as ``dependency'': every type
\(\ty : \univ_{\biglor_{\modesketch} \mode}\) is fractured into a
family of type families \(\{\ty_{\mode(\idx)}\}_{\idx : \modesketch}\)
such that \(\ty_{\mode(\idx)}\) depends on \(\ty_{\mode(\idxI)}\) for
all \(\idxI > \idx\). One may also regard the underlying finite poset
of \(\modesketch\) as a FOLDS signature \parencite{makkai1995folds}.

\begin{example}
  \label{exm:fracture-span}
  When \(\modesketch\) is the mode sketch
  \(\{0 \leftarrow 01 \rightarrow 1\}\), we have an equivalence
  \begin{align*}
    \begin{autobreak}
      \MoveEqLeft
      \univ_{\mode(01) \lor \mode(1) \lor \mode(0)} \simeq
      \ilexists_{\ty_{0} : \univ_{\mode(0)}}
      \ilexists_{\ty_{1} : \univ_{\mode(1)}}
      \ty_{0}
      \to \ty_{1}
      \to \univ_{\mode(01)}.
    \end{autobreak}
  \end{align*}
\end{example}

\begin{example}
  \label{exm:fracture-0-1-2}
  When \(\modesketch\) is the mode sketch \(\{0 \to 1 \to 2\}\) (with no
  thin triangle), we have an equivalence
  \begin{align*}
    \begin{autobreak}
      \MoveEqLeft
      \univ_{\mode(0) \lor \mode(1) \lor \mode(2)} \simeq
      \ilexists_{\ty_{2} : \univ_{\mode(2)}}
      \ilexists_{\ty_{1} : \ty_{2} \to \univ_{\mode(1)}}
      \ilforall_{\var_{2}}
      \ty_{1}(\var_{2})
      \to \univ_{\mode(0)}.
    \end{autobreak}
  \end{align*}
\end{example}
\noindent 
The equivalence in \cref{prop:fracture-and-gluing-alt} nicely
interacts with type constructors, and we derive the \emph{logical
  relation translation} (also called the parametricity translation) of
dependent type theory
\parencite{bernardy2012proofs,shulman2015inverse,uemura2017fibred,lasson2014canonicity}
as a \emph{theorem} in type theory.  Let \(\mode\) and \(\modeI\) be
{\acrLAMs} such that \(\mode \le {}^{\orthMark}\modeI\). Type
constructors in \(\univ_{\mode \lor \modeI}\) behave in the same way
as the definition of the logical relation translation of type
constructors \parencite[Section 3]{lasson2014canonicity} as follows.

\begin{itemize}
\item \(\univ_{\mode \lor \modeI} : \enlarge \univ_{\mode \lor
    \modeI}\) corresponds to the pair
  \begin{math}
    (\univ_{\modeI}, \ilAbs \tyI. \tyI \to \univ_{\mode});
  \end{math}
\item \(\ilUnit : \univ_{\mode \lor \modeI}\) corresponds to the pair
  \begin{math}
    (\ilUnit, \ilAbs \blank. \ilUnit);
  \end{math}
\item Suppose that \(\ty : \univ_{\mode \lor \modeI}\) corresponds
  to a pair \((\ty_{\modeI}, \ty_{\mode})\). Then \((\ty \to
  \univ_{\mode \lor \modeI}) : \enlarge \univ_{\mode \lor \modeI}\)
  corresponds to the pair
  \begin{equation*}
    (\ty_{\modeI} \to \univ_{\modeI},
    \ilAbs \tyI. \ilforall_{\var : \ty_{\modeI}}\ty_{\mode}(\var) \to
    \tyI(\var) \to \univ_{\mode}).
  \end{equation*}
  Indeed,
  \begin{EqReasoning}
    \begin{align*}
      & \term{\ty \to \univ_{\mode \lor \modeI}} \\
      \simeq& \by{fracture and gluing} \\
      & \term{(\ilexists_{\var : \ty_{\modeI}}\ty_{\mode}(\var)) \to (\ilexists_{\tyI : \univ_{\modeI}}\tyI \to \univ_{\mode})} \\
      \simeq& \by{\(\ilforall\) distributes over \(\ilexists\)} \\
      & \term{\ilexists_{\tyI : \ilforall_{\var : \ty_{\modeI}}\ty_{\mode}(\var) \to \univ_{\modeI}}\ilforall_{\var}\ilforall_{\varI}\tyI(\var, \varI) \to \univ_{\mode}} \\
      \simeq& \by{\(\univ_{\modeI} \simeq (\ty_{\mode}(\var) \to \univ_{\modeI})\) since \(\mode \le {}^{\orthMark}\modeI\)} \\
      & \term{\ilexists_{\tyI : \ty_{\modeI} \to \univ_{\modeI}}\ilforall_{\var}\ty_{\mode}(\var) \to \tyI(\var) \to \univ_{\mode}};
    \end{align*}
  \end{EqReasoning}
\item Suppose that \(\ty : \univ_{\mode \lor \modeI}\) corresponds
  to a pair \((\ty_{\modeI}, \ty_{\mode})\) and that
  \(\tyI : \ty \to \univ_{\mode \lor \modeI}\) corresponds to a pair
  \((\tyI_{\modeI}, \tyI_{\mode})\). Then
  \(\ilforall_{\var : \ty}\tyI(\var) : \univ_{\mode \lor \modeI}\)
  corresponds to the pair
  \begin{equation*}
    (\ilforall_{\var_{\modeI} : \ty_{\modeI}}\tyI_{\modeI}(\var_{\modeI}),
    \ilAbs \map. \ilforall_{\var_{\modeI}}\ilforall_{\var_{\mode} : \ty_{\mode}(\var_{\modeI})}\tyI_{\mode}(\var_{\modeI}, \var_{\mode}, \map(\var_{\modeI})))
  \end{equation*}
  by a similar calculation to the previous clause.
  \(\ilexists_{\var : \ty}\tyI(\var) : \univ_{\mode \lor \modeI}\)
  corresponds to the pair
  \begin{equation*}
    (\ilexists_{\var_{\modeI} : \ty_{\modeI}}\tyI_{\modeI}(\var_{\modeI}),
    \ilAbs (\el_{\modeI}, \elI_{\modeI}). \ilexists_{\var_{\mode} : \ty_{\mode}(\el_{\modeI})}\tyI_{\mode}(\el_{\modeI}, \var_{\mode}, \elI_{\modeI}));
  \end{equation*}
\item Suppose that \(\ty : \univ_{\mode \lor \modeI}\) corresponds
  to a pair \((\ty_{\modeI}, \ty_{\mode})\), that \(\el : \ty\)
  corresponds to a pair \((\el_{\modeI}, \el_{\mode})\), and that
  \(\el' : \ty\) corresponds to a pair
  \((\el'_{\modeI}, \el'_{\mode})\). Then
  \(\el = \el' : \univ_{\mode \lor \modeI}\) corresponds to the pair
  \begin{equation*}
    (\el_{\modeI} = \el'_{\modeI},
    \ilAbs \pth. \el_{\mode} =^{\ty_{\mode}}_{\pth} \el'_{\mode}).
  \end{equation*}
\end{itemize}
\noindent 
Thus, any type \(\ty : \univ_{\mode \lor \modeI}\) constructed using
these type constructors is fractured into a type
\(\ty_{\modeI} : \univ_{\modeI}\) and a type family
\(\ty_{\mode} : \ty_{\modeI} \to \univ_{\mode}\), and \(\ty_{\mode}\) is
equivalent to the logical relation translation of \(\ty_{\modeI}\). In
this sense, \emph{types in \(\univ_{\mode \lor \modeI}\) are logical
  relations}. The interaction of the equivalences in
\cref{exm:fracture-span,exm:fracture-0-1-2} and type constructors is
similarly calculated. We thus conclude that types in
\(\univ_{\biglor_{\modesketch} \mode}\) are generalized logical
relations.

\section{Higher category theory}
\label{sec:high-categ-theory}

We collect facts about higher categories needed to develop semantics
of mode sketches in \(\infty\)-logoses.

We work with a model-independent language of \((\infty, 1)\)-category
theory rather than choosing a specific model of
\((\infty, 1)\)-categories. An \emph{\((\infty, 1)\)-category}
\(\cat\) consists of a space \(\Obj(\cat)\) of \emph{objects} of
\(\cat\), a space \(\Map_{\cat}(\obj, \objI)\) of \emph{morphisms from
  \(\obj\) to \(\objI\)} for any objects \(\obj\) and \(\objI\), and
composition operators unital and associative up to coherent
homotopy. Concepts in category theory such as functors, natural
transformations, equivalences, adjoints, (co)limits, and Kan
extensions have \((\infty, 1)\)-categorical analogues. We refer the reader
to \parencite{lurie2009higher,cisinski2019higher,riehl2022elements}
for general \((\infty, 1)\)-category theory.

At least two Grothendieck universes \(\setuniv \in \enlarge \setuniv\)
are assumed to exist. By \emph{small} we mean \(\setuniv\)-small and
by \emph{large} we mean \(\enlarge \setuniv\)-small. For an
\((\infty,1)\)-category \(\cat\) of small objects of some kind, we write
\(\enlarge \cat\) for the \((\infty,1)\)-category of large objects of the
same kind. Let \(\Space\) denote the \((\infty, 1)\)-category of small
spaces. Let \(\Cat\) denote the \((\infty, 1)\)-category of small
\((\infty, 1)\)-categories. For \((\infty, 1)\)-categories \(\cat\) and
\(\catI\), let \(\Fun(\cat, \catI)\) denote the
\((\infty, 1)\)-category of functors from \(\cat\) to \(\catI\) and natural
transformations between them. For a functor
\(\obj : \idxsh \to \cat\), we write
\(\{\proj_{\idx} : \lim_{\idxI \in \idxsh} \obj_{\idxI} \to
\obj_{\idx}\}_{\idx \in \idxsh}\) for the limit cone if it exists and
\(\{\inc_{\idx} : \obj_{\idx} \to \colim_{\idxI \in \idxsh}
\obj_{\idxI}\}_{\idx \in \idxsh}\) for the colimit cocone if it exists.

\emph{Accessible} and \emph{presentable} \((\infty,1)\)-categories
\parencite[Chapter 5]{lurie2009higher} are important classes of
\((\infty,1)\)-categories. We do not need precise definitions of
them. Every presentable \((\infty,1)\)-category is accessible and has small
colimits and limits.

\subsection{$\infty$-logoses}
\label{sec:infty-topoi-infty}

We review the theory of \emph{\(\infty\)-logoses}, also known as
\emph{\(\infty\)-toposes}. The standard reference is \parencite[Chapter
6]{lurie2009higher}.

\begin{definition}
  An \emph{\(\infty\)-logos} is a presentable \((\infty, 1)\)-category
  \(\logos\) such that, for any small \((\infty, 1)\)-category
  \(\idxsh\), any natural transformation
  \(\map : \shI \To \sh : \idxsh \to \logos\) where all the naturality
  squares are pullbacks, and for any cocone over \(\map\) of the form
  \begin{equation}
    \labelX[diagram]{eq:descent-cocone}
    \begin{tikzcd}
      \shI_{\idx}
      \arrow[r, "\mapI_{\idx}"]
      \arrow[d, "\map_{\idx}"'] &
      \shI'
      \arrow[d, "\map'"] \\
      \sh_{\idx}
      \arrow[r, "\inc_{\idx}"'] &
      \colim_{\idx \in \idxsh} \sh_{\idx},
    \end{tikzcd}
  \end{equation}
  \(\shI'\) is the colimit of \(\shI_{\idx}\)'s if and only if
  \cref{eq:descent-cocone} is a pullback for every
  \(\idx \in \idxsh\).  A \emph{morphism of \(\infty\)-logoses} is a functor
  preserving small colimits and finite limits. Note that any morphism
  of \(\infty\)-logoses has a right adjoint by the adjoint functor theorem
  \parencite[Corollary 5.5.2.9]{lurie2009higher}. We write
  \(\Logos \subset \enlarge \Cat\) for the subcategory spanned by the
  \(\infty\)-logoses and the morphisms of \(\infty\)-logoses.
\end{definition}

The following are immediate from the definition.

\begin{proposition}
  \label{prop:logos-descent-along-conservative}
  Let
  \(\{\fun_{\idx} : \logos \to \logosI_{\idx}\}_{\idx \in \idxsh}\) be a
  family of functors between presentable \((\infty, 1)\)-categories
  preserving small colimits and finite limits. If all the
  \(\logosI_{\idx}\)'s are \(\infty\)-logoses and if
  \(\{\fun_{\idx}\}_{\idx \in \idxsh}\) is jointly conservative, then
  \(\logos\) is an \(\infty\)-logos. \qed
\end{proposition}

\begin{proposition}
  \label{prop:logos-init-strict}
  Let \(\logos\) be an \(\infty\)-logos and let \(\sh \in \logos\) be an
  object. If there is a map \(\sh \to \objInitial\), then
  \(\sh \simeq \objInitial\). \qed
\end{proposition}
\noindent 
{\acrLAMs} in homotopy type theory are expected to correspond to (lex,
accessible) localizations of \(\infty\)-logoses.

\begin{definition}
  A morphism of \(\infty\)-logoses is a \emph{localization} if its right
  adjoint is fully faithful. For an \(\infty\)-logos \(\logos\), let
  \(\LexAcc(\logos)\) denote the full subcategory of
  \((\Logos_{\logos /})^{\opMark}\) spanned by the localization
  morphisms \(\logos \to \logosI\). Note that \(\LexAcc(\logos)\) is a
  poset. We call an object in \(\LexAcc(\logos)\) a \emph{lex,
    accessible modality ({\acrLAM}) in \(\logos\)}.
\end{definition}

\begin{notation}
  Let \(\logos\) be an \(\infty\)-logos. For a {\acrLAM} \(\mode\) in
  \(\logos\), we write
  \(\opModality_{\mode} : \logos \to \logos_{\mode}\) for the
  localization corresponding to \(\mode\) and
  \(\unitModality_{\mode}\) for its unit. For two {\acrLAMs}
  \(\mode_{0}\) and \(\mode_{1}\) in \(\logos\), let
  \(\opModality^{\mode_{1}}_{\mode_{0}} : \logos_{\mode_{1}} \to
  \logos_{\mode_{0}}\) denote the restriction of
  \(\opModality_{\mode_{0}}\) along the inclusion
  \(\logos_{\mode_{1}} \subset \logos\). For three {\acrLAMs}
  \(\mode_{0}\), \(\mode_{1}\), and \(\mode_{2}\) in \(\logos\), let
  \(\unitModality^{\mode_{0}; \mode_{2}}_{\mode_{1}} :
  \opModality^{\mode_{2}}_{\mode_{0}} \To
  \opModality^{\mode_{1}}_{\mode_{0}}
  \opModality^{\mode_{2}}_{\mode_{1}} : \logos_{\mode_{2}} \to
  \logos_{\mode_{0}}\) denote the natural transformation defined by
  \((\unitModality^{\mode_{0}; \mode_{2}}_{\mode_{1}})_{A} =
  \opModality_{\mode_{0}} (\unitModality_{\mode_{1}})_{A}\) for
  \(A \in \logos_{\mode_{2}}\).
\end{notation}

\((\infty, 1)\)-categorical counterparts of open and closed modalities are
open and closed, respectively, localizations
(\cref{cst:open-localization}).

\begin{fact}[{\cite[Proposition 6.3.5.1]{lurie2009higher}}]
  \label{fact:etale-morphism-of-logoses}
  Let \(\logos\) be an \(\infty\)-logos. Then, for every object
  \(\sh \in \logos\), the slice \(\logos_{/ \sh}\) is an
  \(\infty\)-logos. Moreover, the pullback functor
  \(\sh^{\pbMark} : \logos \to \logos_{/ \sh}\) is a morphism of
  \(\infty\)-logoses.
\end{fact}

\begin{fact}[{\cite[Propositions 4.3.6 and 4.3.7]{anel2022left-exact}}]
  \label{fact:topological-localization}
  Let \(\logos\) be an \(\infty\)-logos and let \(\cls\) be a class of
  monomorphisms in \(\logos\). Let \(\logosI \subset \logos\) be the full
  subcategory spanned by those objects \(\sh\) such that
  \(\map^{\pbMark} : \Map(\shXI, \sh) \to \Map(\shX, \sh)\) is an
  equivalence for every pullback \(\map\) of a morphism in
  \(\cls\). Then \(\logosI\) is an \(\infty\)-logos, and the inclusion
  \(\logosI \to \logos\) has a left adjoint which is a localization of
  \(\infty\)-logoses.
\end{fact}

\begin{construction}
  \label{cst:open-localization}
  Let \(\logos\) be an \(\infty\)-logos and let \(\propo \in \logos\) be a
  \((-1)\)-truncated object. The \emph{open localization} associated
  to \(\propo\) is \(\propo^{\pbMark} : \logos \to \logos_{/
    \propo}\). This is indeed a morphism of \(\infty\)-logoses by
  \cref{fact:etale-morphism-of-logoses}, and its right adjoint is
  fully faithful since \(\propo\) is \((-1)\)-truncated. The
  \emph{closed localization} associated to \(\propo\) is the
  localization obtained by \cref{fact:topological-localization} for
  the singleton class of monomorphisms \(\{\objInitial \to
  \propo\}\). Let \(\logos \to \logosI\) be the closed localization
  associated to \(\propo\). Then an object \(\sh \in \logos\) belongs to
  \(\logosI\) if and only if
  \(\Map(\shX, \sh) \to \Map(\map^{\pbMark} \objInitial, A)\) is an
  equivalence for every \(\map : \shX \to \propo\). Since
  \(\map^{\pbMark} \objInitial \simeq \objInitial\) by
  \cref{prop:logos-init-strict},
  \(\Map(\map^{\pbMark} \objInitial, \sh) \simeq \objFinal\). Therefore,
  \(\sh\) belongs to \(\logosI\) if and only if
  \(\propo^{\pbMark} \sh \in \logos_{/ \propo}\) is the final object,
  which is equivalent to that the projection
  \(\propo \times \sh \to \propo\) is an equivalence.
\end{construction}

\subsection{The language of $(\infty, 2)$-category theory}
\label{sec:language-infty-2}

We formulate concepts in \((\infty, 2)\)-category theory using the
\emph{\((\infty, 1)\)-category \(\nPrefix{2}\Cat\) of
  \((\infty, 2)\)-categories} axiomatized and proved to be equivalent to
various models of \((\infty, 2)\)-categories by
Barwick and Schommer-Pries \parencite{barwick2021unicity}.

\begin{definition}
  A (strict) \(2\)-category \(\cat\) is said to be \emph{gaunt} if
  only invertible \(1\)-cells and \(2\)-cells are the identities.
\end{definition}

\begin{example}
  The \emph{walking \(\natI\)-cell} \(\walkingCell_{\natI}\) for
  \(0 \le \natI \le 2\) is the \(2\)-category freely generated by a single
  \(\natI\)-cell and is gaunt.
\end{example}

Among Barwick and Schommer-Pries's axioms, the following are
important to us.

\begin{fact}[{\cite[Basic Data]{barwick2021unicity}}]
  \(\nPrefix{2}\Cat\) is presentable and contains finitely presentable
  gaunt \(2\)-categories as a full subcategory.
\end{fact}

\begin{fact}[{\cite[Axiom C.2]{barwick2021unicity}}]
  \label{fact:cells-generator}
  \begin{math}
    \{\walkingCell_{0}, \walkingCell_{1}, \walkingCell_{2}\}
  \end{math}
  is a set of generators for \(\nPrefix{2}\Cat\).
\end{fact}

\begin{fact}[{\cite[Axiom C.3]{barwick2021unicity}}]
  The \((\infty, 1)\)-category \(\nPrefix{2}\Cat\) (more generally
  \(\nPrefix{2}\Cat_{/ \walkingCell_{\natI}}\) for
  \(0 \le \natI \le 2\)) is cartesian closed. For
  \(\cat, \catI \in \nPrefix{2}\Cat\), let \(\Fun(\cat, \catI)\) denote
  the exponential in \(\nPrefix{2}\Cat\).
\end{fact}

\begin{fact}[{\cite[Axiom C.4]{barwick2021unicity}}]
  \(\walkingCell_{0}\), \(\walkingCell_{1}\), and \(\walkingCell_{2}\)
  satisfy certain pushout formulas. We will recall them when needed.
\end{fact}

Let us fix terminology and notation.

\begin{definition}
  Objects in \(\nPrefix{2}\Cat\) are called
  \emph{\((\infty, 2)\)-categories}. For an \((\infty, 2)\)-category
  \(\cat\), morphisms in \(\nPrefix{2}\Cat\) to \(\cat\) from
  \(\walkingCell_{0}\), \(\walkingCell_{1}\), and \(\walkingCell_{2}\)
  are called \emph{objects} or \emph{\(0\)-cells}, \emph{morphisms} or
  \emph{\(1\)-cells}, and \emph{\(2\)-morphisms} or
  \emph{\(2\)-cells}, respectively, in \(\cat\). Morphisms in
  \(\nPrefix{2}\Cat\) are called \emph{functors}. For
  \((\infty, 2)\)-categories \(\cat\) and \(\catI\), \(1\)-cells in the
  \((\infty, 2)\)-category \(\Fun(\cat, \catI)\) are called \emph{natural
    transformations}.
\end{definition}

\Cref{fact:cells-generator} is particularly useful and mostly used in
the following form.

\begin{proposition}
  Let \(\cls\) be a class of small \((\infty, 2)\)-categories. If
  \(\cls\) is closed under small colimits (that is,
  \(\colim_{\idx \in \idxsh} \cat_{\idx} \in \cls\) whenever
  \(\cat_{\idx} \in \cls\) for every \(\idx \in \idxsh\), for any small
  diagram \(\cat : \idxsh \to \nPrefix{2}\Cat\)) and if \(\cls\)
  contains \(\walkingCell_{0}\), \(\walkingCell_{1}\), and
  \(\walkingCell_{2}\), then \(\cls\) contains all small
  \((\infty, 2)\)-categories. \qed
\end{proposition}
\noindent 
\Cref{cst:2-category-obj,cst:2-category-core,cst:2-category-dual}
below are easily justified by using the 2-fold complete Segal spaces
model \parencite{barwick2005thesis}.

\begin{construction}
  \label{cst:2-category-obj}
  Let \(\cat\) be an \((\infty, 2)\)-category. We define the space
  \(\Obj(\cat)\) of objects in \(\cat\) to be
  \(\Map_{\nPrefix{2}\Cat}(\walkingCell_{0}, \cat)\). For a pair of
  objects \((\obj, \objI)\), an \((\infty, 1)\)-category
  \(\Map_{\cat}(\obj, \objI)\) called the \emph{mapping
    \((\infty, 1)\)-category} is constructed. It is defined by the
  pullbacks
  \begin{equation*}
    \begin{tikzcd}
      \Obj(\Map_{\cat}(\obj, \objI))
      \arrow[r]
      \arrow[d]
      \arrow[dr, pbMark] &
      \Map_{\nPrefix{2}\Cat}(\walkingCell_{1}, \cat)
      \arrow[d, "{(\dom, \cod)}"] \\
      \objFinal
      \arrow[r, "{(\obj, \objI)}"'] &
      \Obj(\cat) \times \Obj(\cat)
    \end{tikzcd}
  \end{equation*}
  and
  \begin{equation*}
    \begin{tikzcd}
      \Map_{\Map_{\cat}(\obj, \objI)}(\mor, \morI)
      \arrow[r]
      \arrow[d]
      \arrow[dr, pbMark] &
      \Map_{\nPrefix{2}\Cat}(\walkingCell_{2}, \cat)
      \arrow[d, "{(\dom, \cod)}"] \\
      \objFinal
      \arrow[r, "{(\mor, \morI)}"'] &
      \Obj(\Map_{\cat}(\obj, \objI)) \times \Obj(\Map_{\cat}(\obj, \objI)).
    \end{tikzcd}
  \end{equation*}
  \(\cat\) also has a functorial composition operator between its
  mapping \((\infty, 1)\)-categories.
\end{construction}

\begin{construction}
  \label{cst:2-category-core}
  For an \((\infty, 2)\)-category \(\cat\), we define an \((\infty, 1)\)-category \(\Core_{\nMark{1}}(\cat)\) called the \emph{\((\infty, 1)\)-core of \(C\)} by
  \(\Obj(\Core_{\nMark{1}}(\cat)) = \Obj(\cat)\) and
  \(\Map_{\Core_{\nMark{1}}(\cat)}(\obj, \objI) =
  \Obj(\Map_{\cat}(\obj, \objI))\). This defines a functor
  \(\Core_{\nMark{1}} : \nPrefix{2}\Cat \to \Cat\). It is shown that
  \(\Core_{\nMark{1}}\) has a fully faithful left adjoint, and thus we
  regard \(\Cat\) as a coreflective full subcategory of
  \(\nPrefix{2}\Cat\). An \((\infty, 2)\)-category \(\cat\) is an
  \((\infty, 1)\)-category if and only if it is \emph{locally discrete} in
  the sense that \(\Map_{\cat}(\obj, \objI)\) is an
  \(\infty\)-groupoid for any \(\obj, \objI \in \cat\).
\end{construction}

\begin{definition}
  Let \(\cat\) be an \((\infty, 2)\)-category. A \emph{locally full
    subcategory} of \(\cat\) is an \((\infty, 2)\)-category \(\cat'\)
  equipped with a functor \(\fun : \cat' \to \cat\) such that
  \begin{math}
    \Obj(\cat') \to \Obj(\cat)
  \end{math}
  is mono and
  \begin{math}
    \Map_{\cat'}(\obj, \objI) \to \Map_{\cat}(\fun(\obj), \fun(\objI))
  \end{math}
  is fully faithful for any \(\obj, \objI \in \cat'\). A locally full
  subcategory of \(\cat\) is usually specified by a class
  \(\cat'_{0}\) of \(0\)-cells in \(\cat\) and a class \(\cat'_{1}\)
  of \(1\)-cells in \(\cat\) between objects in \(\cat'_{0}\) such
  that \(\cat'_{1}\) contains all the equivalences between objects in
  \(\cat'_{0}\) and is closed under composition.
\end{definition}

\noindent 
The \emph{unicity theorem} \parencite[Theorem 7.3]{barwick2021unicity}
asserts that there are exactly \((\Integer/2 \Integer)^{2}\)
automorphisms on \(\nPrefix{2}\Cat\). Those automorphisms are
\emph{opposite} constructions.

\begin{construction}
  \label{cst:2-category-dual}
  Let \(\cat\) be an \((\infty, 2)\)-category.
  \((\infty, 2)\)-categories \(\cat^{\opMark(1)}\) and
  \(\cat^{\opMark(2)}\) are defined by
  \begin{math}
    \Obj(\cat^{\opMark(1)}) = \Obj(\cat^{\opMark(2)}) = \Obj(\cat)
  \end{math}
  and
  \begin{align*}
    \Map_{\cat^{\opMark(1)}}(\obj, \objI)
    &= \Map_{\cat}(\objI, \obj) \\
    \Map_{\cat^{\opMark(2)}}(\obj, \objI)
    &= \Map_{\cat}(\obj, \objI)^{\opMark}.
  \end{align*}
  We abbreviate
  \((\cat^{\opMark(1)})^{\opMark(2)} \simeq
  (\cat^{\opMark(2)})^{\opMark(1)}\) as \(\cat^{\opMark(1, 2)}\).
\end{construction}

The following is an \((\infty, 2)\)-categorical version of the fact that a
natural transformation is invertible if it is point-wise
invertible. It is proved without using a specific model.

\begin{proposition}
  \label{prop:invertible-trans-pointwise}
  For any \((\infty, 2)\)-categories \(\cat\) and \(\catI\), the
  restriction functor
  \begin{equation*}
    \Core_{\nMark{1}}(\Fun(\cat, \catI)) \to
    \Core_{\nMark{1}}(\Fun(\Obj(\cat), \catI))
  \end{equation*}
  is conservative.
\end{proposition}
\begin{proof}
  Let \(\cls\) be the class of \((\infty, 2)\)-categories \(\cat\)
  such that the functor
  \begin{math}
    \Core_{\nMark{1}}(\Fun(\cat, \catI)) \to
    \Core_{\nMark{1}}(\Fun(\Obj(\cat), \catI))
  \end{math}
  is conservative for any \((\infty, 2)\)-category \(\catI\). We first
  show that \(\cls\) is closed under colimits. For a functor
  \(\cat : \idxsh \to \nPrefix{2}\Cat\), we have the following
  commutative diagram.
  \begin{equation*}
    \begin{tikzcd}
      \Core_{\nMark{1}}(\Fun(\colim_{\idx \in \idxsh} \cat_{\idx}, \catI))
      \arrow[r]
      \arrow[d, "\simeq"'] &
      \Core_{\nMark{1}}(\Fun(\Obj(\colim_{\idx \in \idxsh} \cat_{\idx}), \catI))
      \arrow[d] \\
      \lim_{\idx \in \idxsh} \Core_{\nMark{1}}(\Fun(\cat_{\idx}, \catI))
      \arrow[r] &
      \lim_{\idx \in \idxsh} \Core_{\nMark{1}}(\Fun(\Obj(\cat_{\idx}), \catI))
    \end{tikzcd}
  \end{equation*}
  The left functor is an equivalence as \(\Core_{\nMark{1}}\)
  preserves limits. If every \(\cat_{\idx}\) belongs to \(\cls\), then
  the bottom functor is conservative, and so is the top.

  It remains to show that \(\cls\) contains \(\walkingCell_{0}\),
  \(\walkingCell_{1}\), and \(\walkingCell_{2}\). Note that when
  \(\cat\) is locally discrete, we have
  \begin{math}
    \Core_{\nMark{1}}(\Fun(\cat, \catI))
    \simeq \Fun(\cat, \Core_{\nMark{1}}(\catI)),
  \end{math}
  and \((\infty, 1)\)-category theory applies. The only non-trivial case is
  when \(\cat = \walkingCell_{2}\). We recall the following pushout in
  \(\nPrefix{2}\Cat\) \parencite[Axiom C.4]{barwick2021unicity}.
  \begin{equation*}
    \begin{tikzcd}
      \walkingCell_{2}
      \arrow[r]
      \arrow[d]
      \arrow[dr, poMark] &
      \walkingCell_{1} +_{\walkingCell_{0}} \walkingCell_{2}
      \arrow[d] \\
      \walkingCell_{2} +_{\walkingCell_{0}} \walkingCell_{1}
      \arrow[r] &
      \walkingCell_{2} \times \walkingCell_{1}
    \end{tikzcd}
  \end{equation*}
  \(\walkingCell_{2} \times \walkingCell_{1}\) should be a cylinder
  \begin{equation*}
    \begin{tikzcd}
      \bullet
      \arrow[r]
      \arrow[d, bend right = 9ex, ""{name = a0}]
      \arrow[d, bend left = 9ex, ""'{name = a1}]
      \arrow[from = a0, to = a1, To] &
      \bullet
      \arrow[d, bend right = 9ex, ""{name = a2}]
      \arrow[d, bend left = 9ex, ""'{name = a3}]
      \arrow[from = a2, to = a3, To] \\
      \bullet
      \arrow[r] &
      \bullet,
    \end{tikzcd}
  \end{equation*}
  and the above pushout formula asserts that this is the case because
  the cylinder is obtained by gluing the lower-left part
  \begin{math}
    \walkingCell_{2} +_{\walkingCell_{0}} \walkingCell_{1} =
    \left(
      \begin{tikzcd}
        \bullet
        \arrow[d, bend right = 9ex, ""{name = a0}]
        \arrow[d, bend left = 9ex, ""'{name = a1}]
        \arrow[from = a0, to = a1, To] \\
        \bullet
        \arrow[r] &
        \bullet
      \end{tikzcd}
    \right)
  \end{math}
  and the upper-right part
  \begin{math}
    \walkingCell_{1} +_{\walkingCell_{0}} \walkingCell_{2} =
    \left(
      \begin{tikzcd}
        \bullet
        \arrow[r] &
        \bullet
        \arrow[d, bend right = 9ex, ""{name = a2}]
        \arrow[d, bend left = 9ex, ""'{name = a3}]
        \arrow[from = a2, to = a3, To] \\
        & \bullet
      \end{tikzcd}
    \right).
  \end{math}
  A functor
  \begin{math}
    \fun \in \Map_{\nPrefix{2}\Cat}(\walkingCell_{1}, \Fun(\walkingCell_{2}, \catI))
    \simeq \Map_{\nPrefix{2}\Cat}(\walkingCell_{2} \times \walkingCell_{1}, \catI)
  \end{math}
  is thus a diagram in \(\catI\) of the form
  \begin{equation*}
    \begin{tikzcd}
      \fun(0, 0)
      \arrow[r, "{\fun(0, 01)}"]
      \arrow[d, bend right = 9ex, ""{name = a0}]
      \arrow[d, bend left = 9ex, ""'{name = a1}]
      \arrow[from = a0, to = a1, To] &
      [3ex]
      \fun(0, 1)
      \arrow[d, bend right = 9ex, ""{name = a2}]
      \arrow[d, bend left = 9ex, ""'{name = a3}]
      \arrow[from = a2, to = a3, To] \\
      \fun(1, 0)
      \arrow[r, "{\fun(1, 01)}"'] &
      \fun(1, 1).
    \end{tikzcd}
  \end{equation*}
  When the morphisms \(\fun(0, 01)\) and \(\fun(1, 01)\) are
  invertible, then one can construct an inverse of \(\fun\) in
  \(\Fun(\walkingCell_{2}, \catI)\).
\end{proof}

\subsection{Scaled simplicial sets}
\label{sec:scal-simpl-sets}

We review one of models for \((\infty, 2)\)-categories, \emph{scaled
  simplicial sets} \parencite{lurie2009goodwillie}. This model is
convenient for presenting \((\infty, 2)\)-categories by combinatorial data,
provides computation of colimits of \((\infty, 2)\)-categories, and gives
the \((\infty, 2)\)-category of \((\infty, 1)\)-categories a useful universal
property.

\begin{definition}
  A \emph{scaled simplicial set} is a simplicial set \(\sh\) equipped
  with a class of \(2\)-simplices called \emph{thin} \(2\)-simplices
  such that all the degenerate \(2\)-simplices are thin.
\end{definition}

A scaled simplicial set \(\sh\) is thought of as a presentation of an
\((\infty, 2)\)-category: the \(0\)-cells are the \(0\)-simplices of
\(\sh\); the \(1\)-cells are generated by \(1\)-simplices of \(\sh\);
the \(2\)-cells are generated by \(2\)-simplices of \(\sh\) in the
following direction;
\begin{equation*}
  \begin{tikzcd}
    \el_{0}
    \arrow[rr, "\phantom{a}"'{name = a0}]
    \arrow[dr] &
    \arrow[from = a0, to = d, To, end anchor = {[yshift = 1ex]}] &
    \el_{2} \\
    & \el_{1}
    \arrow[ur]
  \end{tikzcd}
\end{equation*}
thin \(2\)-simplices are made invertible \(2\)-cells; higher simplices
presents homotopies filling certain diagrams.

\begin{proposition}
  There is a model structure on the category of scaled simplicial sets
  that presents \(\nPrefix{2}\Cat\).
\end{proposition}
\begin{proof}
  See \parencite{lurie2009goodwillie} for the model structure and
  comparison with other models for \((\infty,
  2)\)-categories. Barwick and Schommer-Pries \parencite{barwick2021unicity} show that the
  \((\infty, 1)\)-category presented by this model structure indeed
  satisfies their axioms.
\end{proof}

\begin{exaC}[{\cite[Example 1.1.5.9]{lurie2009higher}}]
  \label{exm:infinity-2-category-from-poset}
  Let \(\idxsh\) be a (decidable) poset and we regard it as a scaled
  simplicial set with no non-degenerate thin \(2\)-simplex. It
  presents a gaunt \(2\)-category \(\freeBracket{\idxsh}\) defined as
  follows. The objects of \(\freeBracket{\idxsh}\) are the elements of
  \(\idxsh\). The mapping category
  \(\Map_{\freeBracket{\idxsh}}(\idx, \idxI)\) for
  \(\idx, \idxI \in \idxsh\) is the poset of totally ordered subsets
  \(\idxshI \subset \idxsh\) with least element \(\idx\) and largest element
  \(\idxI\). A morphism from \(\idx\) to \(\idxI\) is thus a chain
  \begin{math}
    \idx = \idxII_{0} < \idxII_{1} < \dots < \idxII_{\nat} = \idxI.
  \end{math}
  It then follows that \(\Core_{\nMark{1}}(\freeBracket{\idxsh})\) is
  the free category over the strict ordering relation on
  \(\idxsh\). For any \(\idx < \idxI\) in \(\idxsh\), the chain
  \((\idx < \idxI)\) is the initial object in
  \(\Map_{\freeBracket{\idxsh}}(\idx, \idxI)\).
\end{exaC}

We can calculate colimits in \(\nPrefix{2}\Cat\) via homotopy colimits
of scaled simplicial sets. We do not need much details, but a useful
consequence is the following.

\begin{corollary}
  \label{lem:2-cat-colimit-surjective-on-objects}
  For any functor \(\cat : \idxsh \to \nPrefix{2}\Cat\), the map
  \begin{equation*}
    \colim_{\idx \in \idxsh} \Obj(\cat_{\idx})
    \to \Obj(\colim_{\idx \in \idxsh} \cat_{\idx})
  \end{equation*}
  is surjective. \qed
\end{corollary}
\noindent 
Lurie \parencite[Section 4.5]{lurie2009goodwillie} shows that the
\((\infty, 2)\)-category of \((\infty, 1)\)-categories is presented by the
scaled simplicial set classifying \emph{locally cocartesian
  fibrations}.

\begin{definition}
  Let \(\sh\) be a scaled simplicial set and \(\shI\) a simplicial
  set.  A map \(\mapProj : \shI \to \sh\) of simplicial sets is a
  \emph{locally cocartesian fibration} if the following are satisfied:
  \begin{enumerate}
  \item \(\mapProj\) is an inner fibration of simplicial sets;
  \item for every \(1\)-simplex \(\el : \stdsimp^{1} \to \sh\), the base
    change
    \(\el^{\pbMark} \mapProj : \el^{\pbMark} \shI \to \stdsimp^{1}\) is
    a cocartesian fibration of simplicial sets;
  \item for every thin \(2\)-simplex \(\el : \stdsimp^{2} \to \sh\), the
    base change \(\el^{\pbMark} \mapProj\) is a cocartesian fibration
    of simplicial sets.
  \end{enumerate}
\end{definition}

\begin{fact}[{\cite[Corollary 4.5.7]{lurie2009goodwillie}}]
  \label{thm:universal-locally-cocartesian-fibration}
  There exists a \emph{universal locally cocartesian fibration with
    small fibers}
  \begin{math}
    \proj : \ptCat^{\scaledMark} \to \Cat^{\scaledMark}
  \end{math}
  in the following sense:
  \begin{enumerate}
  \item \(\Cat^{\scaledMark}\) is a large fibrant scaled simplicial
    set, \(\ptCat^{\scaledMark}\) is a large simplicial set, and
    \(\proj\) is a locally cocartesian fibration whose fibers are
    small;
  \item for any locally cocartesian fibration with small fibers
    \(\mapProj\), there exists a unique, up to homotopy, homotopy
    pullback from \(\mapProj\) to \(\proj\).
  \end{enumerate}
\end{fact}

\begin{construction}
  Let \(\Cat^{\nMark{2}} \in \enlarge \nPrefix{2}\Cat\) be the
  \((\infty, 2)\)-category presented by \(\Cat^{\scaledMark}\). By the
  construction of \(\Cat^{\scaledMark}\) \parencite[Definition
  4.5.1]{lurie2009goodwillie}, we have
  \begin{math}
    \Obj(\Cat^{\nMark{2}}) \simeq \Obj(\Cat)
  \end{math}
  and
  \begin{math}
    \Map_{\Cat^{\nMark{2}}}(\cat, \catI) \simeq
    \Fun(\cat, \catI).
  \end{math}
\end{construction}

\begin{notation}
  The constructions of \(\Cat^{\scaledMark}\) and \(\Cat^{\nMark{2}}\)
  are parameterized by a universe. We use the notation
  \(\enlarge \Cat^{\scaledMark}\) and \(\enlarge \Cat^{\nMark{2}}\)
  for those constructions with respect to the larger universe
  \(\enlarge \setuniv\).
\end{notation}

The notion of locally cocartesian fibration is, however, specific to
the scaled simplicial sets model. A more model-independent notion is
as follows.

\begin{definition}
  \label{def:1-cocartesian-2-right-fibration}
  A functor \(\funProj : \catI \to \cat\) between
  \((\infty, 2)\)-categories is a \emph{\(1\)-cocartesian \(2\)-right
    fibration} if the following conditions are satisfied:
  \begin{enumerate}
  \item \label{item:2-right-fibration} \(\funProj\) is \emph{locally a
      right fibration} in the sense that for any objects
    \(\obj, \objI \in \catI\), the functor
    \begin{math}
      \funProj : \Map_{\catI}(\obj, \objI)
      \to \Map_{\cat}(\funProj(\obj), \funProj(\objI))
    \end{math}
    is a right fibration;
  \item \label{item:1-cartesian-fibration}
    \(\Core_{\nMark{1}}(\funProj) : \Core_{\nMark{1}}(\catI) \to
    \Core_{\nMark{1}}(\cat)\) is a cocartesian fibration.
  \end{enumerate}
  By duality, \emph{\(1\)-(cocartesian/cartesian) \(2\)-(right/left)
    fibrations} are defined. (These are called \emph{inner/outer
    cocartesian/cartesian fibrations} for the scaled simplicial set
  model \parencite{gagna2020fibrations-arxiv}.)
\end{definition}

\begin{remark}
  Assuming \cref{item:2-right-fibration} in
  \cref{def:1-cocartesian-2-right-fibration}, any cocartesian morphism
  \(\mor : \obj \to \objI\) in \(\Core_{\nMark{1}}(\catI)\) is also
  cocartesian in the \(\Cat\)-enriched sense: for any object
  \(\objII \in \catI\), the square
  \begin{equation*}
    \begin{tikzcd}
      \Map_{\catI}(\objI, \objII)
      \arrow[r, "\blank \comp \mor"]
      \arrow[d, "\funProj"'] &
      [4ex]
      \Map_{\catI}(\obj, \objII)
      \arrow[d, "\funProj"] \\
      \Map_{\cat}(\funProj(\objI), \funProj(\objII))
      \arrow[r, "\blank \comp \funProj(\mor)"'] &
      \Map_{\cat}(\funProj(\obj), \funProj(\objII))
    \end{tikzcd}
  \end{equation*}
  is a pullback in \(\Cat\), because \(\Obj : \Cat \to \Space\)
  reflects pullbacks of right fibrations.
\end{remark}

\begin{proposition}
  \label{prop:universal-1-cocart-2-right-fib}
  \(\Cat^{\nMark{2}}\) is part of a universal \(1\)-cocartesian
  \(2\)-right fibration with small fibers
  \(\ptCat^{\nMark{2}, \rightMark} \to \Cat^{\nMark{2}}\). For a functor
  \(\cat : \idxsh \to \Cat^{\nMark{2}}\), let
  \(\El_{\idxsh}(\cat) \to \idxsh\) denote the corresponding
  \(1\)-cocartesian \(2\)-right fibration.
\end{proposition}
\begin{proof}
  This follows from the equivalence between locally cocartesian
  fibrations and inner cartesian fibrations of categories enriched
  over marked simplicial sets given by Gagna, Harpaz, and Lanari \parencite[Propositions 2.4.1
  and 3.1.3]{gagna2020fibrations-arxiv}.
\end{proof}
\noindent 
An object in \(\ptCat^{\nMark{2}, \rightMark}\) is a pair
\((\cat, \obj)\) consisting of an \((\infty, 1)\)-category \(\cat\) and an
object \(\obj \in \cat\). A morphism
\((\cat, \obj) \to (\catI, \objI)\) in
\(\ptCat^{\nMark{2}, \rightMark}\) is a pair \((\fun, \mor)\)
consisting of a functor \(\fun : \cat \to \catI\) and a morphism
\(\mor : \fun(\obj) \to \objI\). A \(2\)-morphism \((\fun, \mor) \To
(\funI, \morI) : (\cat, \obj) \to (\catI, \objI)\) in
\(\ptCat^{\nMark{2}, \rightMark}\) is a pair \((\trans, \pth)\)
consisting of a natural transformation \(\trans : \fun \To \funI\) and
a path \(\pth : \mor \sim \morI \comp \trans_{\obj}\).

For the purpose of \cref{sec:oplax-limits}, we introduce a variant of
\(\ptCat^{\nMark{2}, \rightMark}\).

\begin{construction}
  The equivalence of \((\infty, 1)\)-categories
  \(\Cat \ni \cat \mapsto \cat^{\opMark} \in \Cat\) extends to an equivalence of
  \((\infty, 2)\)-categories
  \begin{equation*}
    (\blank)^{\opMark} : \Cat^{\nMark{2}}
    \simeq (\Cat^{\nMark{2}})^{\opMark(2)}.
  \end{equation*}
  Let
  \(\ptCat^{\nMark{2}, \leftMark} = ((\blank)^{\opMark})^{\pbMark}
  (\ptCat^{\nMark{2}, \rightMark})^{\opMark(2)}\). By construction,
  the functor
  \((\ptCat^{\nMark{2}, \leftMark})^{\opMark(1, 2)} \to
  (\Cat^{\nMark{2}})^{\opMark(1, 2)}\) classifies \(1\)-cartesian
  \(2\)-right fibrations with small fibers. For a functor
  \(\cat : \idxsh^{\opMark(1, 2)} \to \Cat^{\nMark{2}}\), let
  \(\El_{\idxsh}(\cat) \to \idxsh\) denote the corresponding
  \(1\)-cartesian \(2\)-right fibration.
\end{construction}

In other words, \(\ptCat^{\nMark{2}, \leftMark}\) is the fiberwise
opposite of \(\ptCat^{\nMark{2}, \rightMark}\). Thus, the objects of
\(\ptCat^{\nMark{2}, \leftMark}\) are the same as
\(\ptCat^{\nMark{2}, \rightMark}\), but a morphism
\((\cat, \obj) \to (\catI, \objI)\) in
\(\ptCat^{\nMark{2}, \leftMark}\) is a pair \((\fun, \mor)\)
consisting of a functor \(\fun : \cat \to \catI\) and a morphism
\(\mor : \objI \to \fun(\obj)\).

\subsection{Oplax limits}
\label{sec:oplax-limits}

\emph{Oplax limits} in general \((\infty, 2)\)-categories are defined by
Gagna, Harpaz, and Lanari \parencite{gagna2020fibrations-arxiv}. In this paper, we only need
oplax limits of \((\infty, 1)\)-categories which have the following simple
construction \parencite[Example 5.3.12]{gagna2020fibrations-arxiv}.

\begin{construction}
  Let \(\idxsh\) be a small \((\infty, 2)\)-category and
  \(\cat : \idxsh^{\opMark(1, 2)} \to \Cat^{\nMark{2}}\) a
  functor. The \emph{oplax limit of \(\cat\)} is defined to be the
  pullback
  \begin{equation*}
    \begin{tikzcd}
      \opLaxLim_{\idx \in \idxsh} \cat_{\idx}
      \arrow[r, dotted]
      \arrow[d, dotted]
      \arrow[dr, pbMark] &
      \Fun(\idxsh, (\ptCat^{\nMark{2}, \leftMark})^{\opMark(1, 2)})
      \arrow[d] \\
      \walkingCell_{0}
      \arrow[r, "\cat"'] &
      \Fun(\idxsh, (\Cat^{\nMark{2}})^{\opMark(1, 2)}).
    \end{tikzcd}
  \end{equation*}
  Equivalently, it is the pullback
  \begin{equation*}
    \begin{tikzcd}
      \opLaxLim_{\idx \in \idxsh} \cat_{\idx}
      \arrow[r, dotted]
      \arrow[d, dotted]
      \arrow[dr, pbMark] &
      \Fun(\idxsh, \El_{\idxsh}(\cat))
      \arrow[d] \\
      \walkingCell_{0}
      \arrow[r, "\id_{\idxsh}"'] &
      \Fun(\idxsh, \idxsh).
    \end{tikzcd}
  \end{equation*}
  In other words, \(\opLaxLim_{\idx \in \idxsh} \cat_{\idx}\) is the
  \((\infty, 2)\)-category of sections of
  \(\El_{\idxsh}(\cat) \to \idxsh\). Any functor
  \(\fun : \idxshI \to \idxsh\) induces a functor
  \begin{math}
    \fun^{\pbMark} :
    \opLaxLim_{\idx \in \idxsh} \cat_{\idx}
    \to \opLaxLim_{\idxI \in \idxshI} \cat_{\fun(\idxI)}.
  \end{math}
\end{construction}

\begin{remark}
  \(\opLaxLim_{\idx \in \idxsh} \cat_{\idx}\) is small and locally
  discrete. Indeed, the class of small \((\infty, 2)\)-categories
  \(\idxsh\) such that \(\opLaxLim_{\idx \in \idxsh} \cat_{\idx}\) is
  small and locally discrete for any
  \(\cat : \idxsh^{\opMark(1, 2)} \to \Cat^{\nMark{2}}\) is closed under
  small colimits, and the cases when
  \(\idxsh = \walkingCell_{0}, \walkingCell_{1}, \walkingCell_{2}\)
  are directly calculated in
  \cref{exm:oplax-limit-cell-0,exm:oplax-limit-cell-1,exm:oplax-limit-cell-2}
  below. Hence, we regard \(\opLaxLim_{\idx \in \idxsh} \cat_{\idx}\) as
  a small \((\infty, 1)\)-category.
\end{remark}

Concretely, an object \(\obj\) in
\(\opLaxLim_{\idx \in \idxsh} \cat_{\idx}\) consists of: an object
\(\obj_{\idx} \in \cat_{\idx}\) for any object \(\idx \in \idxsh\); a
morphism
\(\obj_{\moridx} : \obj_{\idx} \to \cat_{\moridx}(\obj_{\idxI})\) for
any morphism \(\moridx : \idx \to \idxI\) in \(\idxsh\); some coherence
data. A morphism \(\mor : \obj \to \objI\) in
\(\opLaxLim_{\idx \in \idxsh} \cat_{\idx}\) consists of: a morphism
\(\mor_{\idx} : \obj_{\idx} \to \objI_{\idx}\) for any object
\(\idx \in \idxsh\); a homotopy \(\mor_{\moridx}\) filling the square
\begin{equation*}
  \begin{tikzcd}
    \obj_{\idx}
    \arrow[r, "\mor_{\idx}"]
    \arrow[d, "\obj_{\moridx}"'] &
    [4ex]
    \objI_{\idx}
    \arrow[d, "\objI_{\moridx}"] \\
    \cat_{\moridx}(\obj_{\idxI})
    \arrow[r, "\cat_{\moridx}(\mor_{\idxI})"'] &
    \cat_{\moridx}(\objI_{\idxI})
  \end{tikzcd}
\end{equation*}
for any morphism \(\moridx : \idx \to \idxI\) in \(\idxsh\); some
coherence data.

\begin{example}
  \label{exm:oplax-limit-discrete}
  When \(\idxsh\) is an \(\infty\)-groupoid, we have
  \begin{math}
    \opLaxLim_{\idx \in \idxsh} \cat_{\idx}
    \simeq \prod_{\idx \in \idxsh} \cat_{\idx}.
  \end{math}
  For a general \(\idxsh\), we have the forgetful functor
  \begin{math}
    \opLaxLim_{\idx \in \idxsh} \cat_{\idx}
    \to \opLaxLim_{\idx \in \Obj(\idxsh)} \cat_{\idx}
    \simeq \prod_{\idx \in \idxsh} \cat_{\idx}
  \end{math}
  which is conservative by \cref{prop:invertible-trans-pointwise}.
\end{example}

\begin{example}
  \label{exm:oplax-limit-cell-0}
  When \(\idxsh = \walkingCell_{0}\), a functor
  \(\walkingCell_{0}^{\opMark(1, 2)} \to \Cat^{\nMark{2}}\)
  corresponds to an \((\infty, 1)\)-category \(\cat\), and its oplax
  limit is \(\cat\) itself.
\end{example}

\begin{example}
  \label{exm:oplax-limit-cell-1}
  When \(\idxsh = \walkingCell_{1}\), a functor
  \(\cat : \walkingCell_{1}^{\opMark(1, 2)} \to \Cat^{\nMark{2}}\)
  corresponds to a functor \(\fun : \cat_{1} \to \cat_{0}\). Its oplax
  limit is the \((\infty, 1)\)-category of triples \((\obj_{0}, \obj_{1},
  \mor)\) consisting of objects \(\obj_{0} \in \cat_{0}\) and \(\obj_{1}
  \in \cat_{1}\) and a morphism \(\mor : \obj_{0} \to \fun(\obj_{1})\). In
  other words, we have the following pullback.
  \begin{equation*}
    \begin{tikzcd}
      \opLaxLim_{\idx \in \walkingCell_{1}} \cat_{\idx}
      \arrow[r]
      \arrow[d]
      \arrow[dr, pbMark] &
      \cat_{0}^{\to}
      \arrow[d, "\cod"] \\
      \cat_{1}
      \arrow[r, "\fun"'] &
      \cat_{0}
    \end{tikzcd}
  \end{equation*}
  This oplax limit is called the \emph{Artin gluing} for \(\fun\) and
  denoted by \(\Glue(\fun)\).
\end{example}

\begin{example}
  \label{exm:oplax-limit-cell-2}
  When \(\idxsh = \walkingCell_{2}\), a functor
  \(\cat : \walkingCell_{2}^{\opMark(1, 2)} \to \Cat^{\nMark{2}}\)
  corresponds to a natural transformation
  \(\trans : \fun_{1} \To \fun_{0} : \cat_{1} \to \cat_{0}\). Since
  \((\ptCat^{\nMark{2}, \leftMark})^{\opMark(1, 2)} \to
  (\Cat^{\nMark{2}})^{\opMark(1, 2)}\) is locally a right fibration, a
  section of it over \(\cat\) is completely determined by the
  restriction along the codomain inclusion
  \(\walkingCell_{1} \to \walkingCell_{2}\). Therefore,
  \begin{math}
    \opLaxLim_{\idx \in \walkingCell_{2}} \cat_{\idx}
    \simeq \Glue(\fun_{1}).
  \end{math}
\end{example}

\begin{example}
  \label{exm:oplax-limit-colimit}
  Let \(\idxshI\) be a small \((\infty, 1)\)-category and
  \(\idxsh : \idxshI \to \nPrefix{2}\Cat\) a functor. For any functor
  \(\cat : (\colim_{\idxI \in \idxshI} \idxsh_{\idxI})^{\opMark(1, 2)}
  \to \Cat^{\nMark{2}}\), we have a canonical equivalence
  \begin{equation*}
    \opLaxLim_{\idx \in \colim_{\idxI \in \idxshI} \idxsh_{\idxI}} \cat_{\idx}
    \simeq \lim_{\idxI \in \idxshI} \opLaxLim_{\idx \in \idxsh_{\idxI}} \cat_{\inc_{\idxI}(\idx)}.
  \end{equation*}
  Hence, arbitrary oplax limits are constructed from
  \cref{exm:oplax-limit-cell-0,exm:oplax-limit-cell-1,exm:oplax-limit-cell-2}
  using small limits.
\end{example}

We show that \(\infty\)-logoses are closed under oplax limits. For this, we
consider a weaker notion of morphism of \(\infty\)-logoses.

\begin{definition}
  A functor between accessible categories is \emph{accessible} if it
  preserves small \(\card\)-filtered colimits for some regular
  cardinal \(\card\).
\end{definition}

\begin{notation}
  We define
  \(\Logos_{\LexAccMark}^{\nMark{2}} \subset
  \enlarge \Cat^{\nMark{2}}\) to be the locally full subcategory
  spanned by the \(\infty\)-logoses and the lex, accessible functors
  between \(\infty\)-logoses.
\end{notation}

\begin{theorem}
  \label{prop:logos-accessible-oplax-limit}
  Let \(\idxsh\) be a small \((\infty, 2)\)-category and
  \(\logos : \idxsh^{\opMark(1, 2)} \to
  \Logos_{\LexAccMark}^{\nMark{2}}\) a functor. Then
  \(\opLaxLim_{\idx \in \idxsh} \logos_{\idx}\) is an
  \(\infty\)-logos. Moreover, the forgetful functor
  \begin{math}
    \opLaxLim_{\idx \in \idxsh} \logos_{\idx}
    \to \prod_{\idx \in \idxsh} \logos_{\idx}
  \end{math}
  preserves small colimits and finite limits.
\end{theorem}

The rest of this subsection is devoted to the proof of
\cref{prop:logos-accessible-oplax-limit}.

\begin{fact}[{\cite[Proposition 6.3.2.3]{lurie2009higher}}]
  \label{lem:logos-limit}
  \(\Logos \subset \enlarge \Cat\) is closed under small limits.
\end{fact}

\begin{lemma}
  \label{lem:gluing-logos}
  Let \(\fun : \logos_{1} \to \logos_{0}\) be a lex, accessible functor
  between \(\infty\)-logoses. Then \(\Glue(\fun)\) is an
  \(\infty\)-logos, and the projections \(\Glue(\fun) \to \logos_{0}\) and
  \(\Glue(\fun) \to \logos_{1}\) preserve small colimits and finite
  limits.
\end{lemma}
\begin{proof}
  We use a characterization of presentability: an
  \((\infty, 1)\)-category is presentable if and only if it is accessible
  and has small colimits \parencite[Definition
  5.5.0.1]{lurie2009higher}. It follows from \parencite[Propositions
  5.4.4.3 and 5.4.6.6]{lurie2009higher} that \(\Glue(\fun)\) is
  accessible. By construction, \(\Glue(\fun)\) is the
  \((\infty, 1)\)-category of triples \((\sh_{0}, \sh_{1}, \map)\)
  consisting of objects \(\sh_{0} \in \logos_{0}\) and
  \(\sh_{1} \in \logos_{1}\) and a map
  \(\map : \sh_{0} \to \fun(\sh_{1})\). It follows from this description
  that the projection \(\Glue(\fun) \to \logos_{0} \times \logos_{1}\)
  creates small colimits and finite limits. In particular,
  \(\Glue(\fun)\) admits small colimits and thus is presentable. Since
  the projection \(\Glue(\fun) \to \logos_{0} \times \logos_{1}\) is
  conservative by \cref{exm:oplax-limit-discrete}, \(\Glue(\fun)\) is
  an \(\infty\)-logos by \cref{prop:logos-descent-along-conservative}.
\end{proof}

\begin{proof}[Proof of \cref{prop:logos-accessible-oplax-limit}]
  Let \(\cls\) be the class of small \((\infty, 2)\)-categories
  \(\idxsh\) such that for any functor
  \(\logos : \idxsh^{\opMark(1, 2)} \to
  \Logos_{\LexAccMark}^{\nMark{2}}\), the oplax limit
  \(\opLaxLim_{\idx \in \idxsh} \logos_{\idx}\) is an \(\infty\)-logos
  and the forgetful functor
  \(\opLaxLim_{\idx \in \idxsh} \logos_{\idx} \to \prod_{\idx \in
    \idxsh} \logos_{\idx}\) preserves small colimits and finite
  limits. It is enough to show that \(\cls\) is closed under small
  colimits and contains \(\walkingCell_{0}\), \(\walkingCell_{1}\),
  and \(\walkingCell_{2}\).

  Let \(\idxsh : \idxshI \to \nPrefix{2}\Cat\) be a functor from a
  small \((\infty, 1)\)-category \(\idxshI\) and suppose that every
  \(\idxsh_{\idxI}\) belongs to \(\cls\). Let
  \(\logos : (\colim_{\idxI \in \idxshI}\idxsh_{\idxI})^{\opMark(1, 2)} \to
  \Logos_{\LexAccMark}^{\nMark{2}}\) be a functor. As in
  \cref{exm:oplax-limit-colimit}, we have
  \begin{equation*}
    \opLaxLim_{\idx \in \colim_{\idxI \in \idxshI} \idxsh_{\idxI}} \logos_{\idx}
    \simeq \lim_{\idxI \in \idxshI} \opLaxLim_{\idx \in \idxsh_{\idxI}} \logos_{\inc_{\idxI}(\idx)}.
  \end{equation*}
  Since \(\idxsh_{\idxI} \in \cls\), small colimits and finite limits in
  \(\opLaxLim_{\idx \in \idxsh_{\idxI}} \logos_{\inc_{\idxI}(\idx)}\)
  are computed in
  \(\prod_{\idx \in \idxsh_{\idxI}} \logos_{\inc_{\idxI}(\idx)}\). Thus, for
  any morphism \(\idxI_{1} \to \idxI_{2}\) in \(\idxshI\), the functor
  \begin{math}
    \opLaxLim_{\idx \in \idxsh_{\idxI_{2}}} \logos_{\inc_{\idxI_{1}}(\idx)}
    \to \opLaxLim_{\idx \in \idxsh_{\idxI_{1}}} \logos_{\inc_{\idxI_{2}}(\idx)}
  \end{math}
  preserves small colimits and finite limits.
  It then follows from \cref{lem:logos-limit} that
  \begin{math}
    \lim_{\idxI \in \idxshI} \opLaxLim_{\idx \in \idxsh_{\idxI}} \logos_{\inc_{\idxI}(\idx)}
  \end{math}
  is an \(\infty\)-logos. Consider the following commutative square.
  \begin{equation*}
    \begin{tikzcd}
      \opLaxLim_{\idx \in \colim_{\idxI \in \idxshI} \idxsh_{\idxI}} \logos_{\idx}
      \arrow[r, "\simeq"]
      \arrow[d] &
      \lim_{\idxI \in \idxshI} \opLaxLim_{\idx \in \idxsh_{\idxI}} \logos_{\inc_{\idxI}(\idx)}
      \arrow[d] \\
      \prod_{\idx \in \Obj(\colim_{\idxI \in \idxshI} \idxsh_{\idxI})} \logos_{\idx}
      \arrow[r] &
      \lim_{\idxI \in \idxshI} \prod_{\idx \in \Obj(\idxsh_{\idxI})} \logos_{\inc_{\idxI}(\idx)}
    \end{tikzcd}
  \end{equation*}
  We have seen that the top functor is an equivalence. The right
  functor preserves small colimits and finite limits as every
  \(\idxsh_{\idxI}\) belongs to \(\cls\). The bottom functor is
  equivalent to the restriction along
  \begin{math}
    \colim_{\idxI \in \idxshI} \Obj(\idxsh_{\idxI}) \to
    \Obj(\colim_{\idxI \in \idxshI} \idxsh_{\idxI})
  \end{math}
  and thus preserves small colimits and finite limits. By
  \cref{lem:2-cat-colimit-surjective-on-objects}, the bottom functor
  is conservative. We thus conclude that the left functor preserves
  small colimits and finite limits. Hence,
  \(\colim_{\idxI \in \idxshI} \idxsh_{\idxI}\) belongs to \(\cls\).

  \noindent 
  \(\walkingCell_{0}\) belongs to \(\cls\) by
  \cref{exm:oplax-limit-cell-0}. \(\walkingCell_{1}\) belongs to
  \(\cls\) by
  \cref{exm:oplax-limit-cell-1,lem:gluing-logos}. \(\walkingCell_{2}\)
  belongs to \(\cls\) by \cref{exm:oplax-limit-cell-2} and by the case
  of \(\walkingCell_{1}\).
\end{proof}

\subsection{Oplax natural transformations}
\label{sec:oplax-transf}

An alternative description of oplax limits is that they are
\((\infty, 1)\)-categories of \emph{oplax natural transformations}
(\cref{prop:oplax-limit-universal-property}).

\begin{consC}[{\cite[Definition 2.1]{gagna2020gray}}]
  Let \(\sh\) and \(\shI\) be scaled simplicial sets. The \emph{Gray
    product} \(\sh \tensorGray \shI\) is the scaled simplicial set
  whose underlying simplicial set is the cartesian product of \(\sh\)
  and \(\shI\) and whose \(2\)-simplex
  \((\el, \elI) : \stdsimp^{2} \to \sh \times \shI\) is thin if both
  \(\el\) and \(\elI\) are thin and either \(\el\) degenerates along
  \(\stdsimp^{\{1, 2\}}\) or \(\elI\) degenerates along
  \(\stdsimp^{\{0, 1\}}\).
\end{consC}

\begin{fact}[{\cite[Theorem 2.17]{gagna2020gray}}]
  \label{prop:gray-tensor-quillen-bifunctor}
  The Gray product is part of a left Quillen bifunctor on
  scaled simplicial sets. Consequently, it induces a functor
  \begin{equation*}
    \blank \tensorGray \blank :
    \nPrefix{2}\Cat \times \nPrefix{2}\Cat
    \to \nPrefix{2}\Cat
  \end{equation*}
  preserving small colimits on each variable.
\end{fact}

\begin{construction}
  By \cref{prop:gray-tensor-quillen-bifunctor} and by the adjoint
  functor theorem, for any \((\infty, 2)\)-category \(\cat\), the functors
  \((\blank \tensorGray \cat)\) and \((\cat \tensorGray \blank)\) have
  right adjoints \(\Fun(\cat, \blank)_{\LaxMark}\) and
  \(\Fun(\cat, \blank)_{\opLaxMark}\), respectively.
\end{construction}

\begin{example}
  \begin{math}
    \cat \tensorGray \walkingCell_{0} \simeq \cat,
  \end{math}
  and thus
  \begin{math}
    \Obj(\Fun(\cat, \catI)_{\opLaxMark}) \simeq
    \Map_{\nPrefix{2}\Cat}(\cat, \catI).
  \end{math}
  Dually,
  \(\Obj(\Fun(\cat, \catI)_{\LaxMark}) \simeq \Map_{\nPrefix{2}\Cat}(\cat,
  \catI)\).
\end{example}

\begin{example}
  Let \(\cat\) be a scaled simplicial set and \(\mor : \obj \to \objI\)
  a \(1\)-simplex in \(\cat\). Consider the following \(2\)-simplices
  in \(\cat \tensorGray \stdsimp^{1}\).
  \begin{equation}
    \label{eq:gray-triangles-1}
    \begin{tikzcd}
      (\obj, 0)
      \arrow[r, "{(\obj, 0 \le 1)}"]
      \arrow[d, "{(\mor, 0)}"']
      \arrow[dr, "{(\mor, 0 \le 1)}"{description}] &
      [4ex]
      (\obj, 1)
      \arrow[d, "{(\mor, 1)}"] \\
      (\objI, 0)
      \arrow[r, "{(\objI, 0 \le 1)}"'] &
      (\objI, 1)
    \end{tikzcd}
  \end{equation}
  By definition, the lower \(2\)-simplex is thin, but the upper one is
  not (unless \(\mor\) is degenerate). Hence, these \(2\)-simplices
  compose and yields a \(2\)-cell
  \begin{equation}
    \label{eq:gray-square-1}
    \begin{tikzcd}
      (\obj, 0)
      \arrow[r, "{(\obj, 0 \le 1)}"]
      \arrow[d, "{(\mor, 0)}"'] &
      [4ex]
      (\obj, 1)
      \arrow[d, "{(\mor, 1)}"]\\
      (\objI, 0)
      \arrow[r, "{(\objI, 0 \le 1)}"']
      \arrow[ur, To, end anchor = {[xshift = -1ex, yshift = -1ex]},
      start anchor = {[xshift = 1ex, yshift = 1ex]}] &
      (\objI, 1)
    \end{tikzcd}
  \end{equation}
  in the \((\infty, 2)\)-category presented by
  \(\cat \tensorGray \stdsimp^{1}\). Then, a \(1\)-cell \(\trans :
  \fun \to \funI\) in \(\Fun(\cat, \catI)_{\opLaxMark}\) assigns: a
  \(1\)-cell \(\trans_{\obj} : \fun(\obj) \to \funI(\obj)\) to every
  \(0\)-cell \(\obj \in \cat\); a \(2\)-cell
  \begin{equation*}
    \begin{tikzcd}
      \fun(\obj)
      \arrow[r, "\trans_{\obj}"]
      \arrow[d, "\fun(\mor)"'] &
      [2ex]
      \funI(\obj)
      \arrow[d, "\funI(\mor)"] \\
      [2ex]
      \fun(\objI)
      \arrow[r, "\trans_{\objI}"']
      \arrow[ur, To, dotted, "\trans_{\mor}",
      start anchor = {[xshift = 1ex, yshift = 1ex]},
      end anchor = {[xshift = -1ex, yshift = -1ex]}] &
      \funI(\objI)
    \end{tikzcd}
  \end{equation*}
  to every \(1\)-cell \(\mor : \obj \to \objI\); and coherence data to
  higher cells. Such a structure is called an \emph{oplax natural
    transformation from \(\fun\) to \(\funI\)}. Dually, \(1\)-cells in
  \(\Fun(\cat, \catI)_{\LaxMark}\) are called \emph{lax natural
    transformations}. For a lax natural transformation
  \(\trans : \fun \to \funI\), the \(2\)-cell \(\trans_{\mor}\) is in
  the opposite direction
  \(\funI(\mor) \comp \trans_{\obj} \To \trans_{\objI} \comp
  \fun(\mor)\).
\end{example}

\begin{remark}
  By definition, we have a map
  \(\sh \tensorGray \shI \to \sh \times \shI\) of scaled simplicial sets
  which exhibits \(\sh \times \shI\) as the one obtained from
  \(\sh \tensorGray \shI\) by making the upper \(2\)-simplex in the
  diagram of the form \labelcref{eq:gray-triangles-1} thin. Thus, for
  \((\infty, 2)\)-categories \(\cat\) and \(\catI\), the cartesian product
  \(\cat \times \catI\) is obtained from the Gray product
  \(\cat \tensorGray \catI\) by making the \(2\)-cell in the diagram
  of the form \labelcref{eq:gray-square-1} invertible. By an adjoint
  argument, \(\Fun(\cat, \catI)\) is regarded as the locally full
  subcategory of \(\Fun(\cat, \catI)_{\LaxMark}\) whose \(1\)-cells
  are the lax natural transformations \(\trans\) such that the
  \(2\)-cell \(\trans_{\mor}\) is invertible for any \(1\)-cell
  \(\mor\) in \(\cat\). We may also regard \(\Fun(\cat, \catI)\) as a
  locally full subcategory of \(\Fun(\cat, \catI)_{\opLaxMark}\) in
  the same way.
\end{remark}

Lax natural transformations correspond to functors between
\(1\)-cocartesian \(2\)-right fibrations.

\begin{proposition}
  \label{prop:lax-univalence}
  Let \(\idxsh\) be an \((\infty, 2)\)-category and let
  \(\cat, \catI : \idxsh \to \Cat^{\nMark{2}}\) be functors. We have an
  equivalence
  \begin{equation*}
    \Map_{\Core_{\nMark{1}}(\Fun(\idxsh
      ,\Cat^{\nMark{2}})_{\LaxMark})}(\cat, \catI)
    \simeq \Map_{\nPrefix{2}\Cat_{/ \idxsh}}(\El_{\idxsh}(\cat), \El_{\idxsh}(\catI))
  \end{equation*}
  natural in \(\idxsh\).
\end{proposition}

\Cref{prop:lax-univalence} follows from the following special cases
which are already known.

\begin{fact}[{\cite[Corollary 4.4.3]{gagna2020fibrations-arxiv}}]
  \label{lem:lax-univalence-1}
  Let \(\idxsh\) be an \((\infty, 2)\)-category. For any functor
  \(\catI : \idxsh \to \Cat^{\nMark{2}}\), we have an equivalence
  \begin{equation*}
    \Map_{\Fun(\idxsh, \Cat^{\nMark{2}})_{\LaxMark}}(\lambda \blank. \objFinal, \catI)
    \simeq \Map_{\nPrefix{2}\Cat^{\nMark{2}}_{/ \idxsh}}(\idxsh, \El_{\idxsh}(\catI)).
  \end{equation*}
\end{fact}

\begin{fact}[{\cite[Theorem E]{haugseng2021mate-arxiv}}]
  \label{lem:lax-univalence-2}
  Let \(\idxsh\) be an \((\infty, 1)\)-category. The map
  \begin{math}
    (\cat : \idxsh \to \Cat^{\nMark{2}})
    \mapsto \El_{\idxsh}(\cat)
  \end{math}
  induces an equivalence between
  \(\Fun(\idxsh, \Cat^{\nMark{2}})_{\LaxMark}\) and the full
  subcategory of \(\Cat^{\nMark{2}}_{/ \idxsh}\) spanned by the
  cocartesian fibrations over \(\idxsh\). Dually,
  \(\Fun(\idxsh^{\opMark}, \Cat^{\nMark{2}})_{\opLaxMark}\) is
  equivalent to the full subcategory of
  \(\Cat^{\nMark{2}}_{/ \idxsh}\) spanned by the cartesian fibrations
  over \(\idxsh\).
\end{fact}

\begin{proof}[Proof of \cref{prop:lax-univalence}]
  We have a lax natural transformation
  \(\unit : (\lambda \blank. \objFinal) \to \cat
  \restrict_{\El_{\idxsh}(\cat)}\) corresponding to the diagonal
  functor
  \(\El_{\idxsh}(\cat) \to \El_{\idxsh}(\cat) \times_{\idxsh}
  \El_{\idxsh}(\cat)\) by \cref{lem:lax-univalence-1}. The
  precomposition with \(\unit\) induces a map
  \begin{align}
    \label{eq:lax-univalence-canonical-map}
    \begin{autobreak}
      \MoveEqLeft
      \Map_{\Core_{\nMark{1}}(\Fun(\idxsh, \Cat^{\nMark{2}})_{\LaxMark})}(\cat, \catI)
      \to \Map_{\Core_{\nMark{1}}(\Fun(\El_{\idxsh}(\cat), \Cat^{\nMark{2}})_{\LaxMark})}(\lambda \blank. \objFinal, \catI \restrict_{\El_{\idxsh}(\cat)}),
    \end{autobreak}
  \end{align}
  and the codomain is by \cref{lem:lax-univalence-1} equivalent to
  \begin{math}
    \Map_{\nPrefix{2}\Cat_{/ \El_{\idxsh}(\cat)}}(\El_{\idxsh}(\cat), \El_{\idxsh}(\cat) \times_{\idxsh} \El_{\idxsh}(\catI))
    \simeq \Map_{\nPrefix{2}\Cat_{/ \idxsh}}(\El_{\idxsh}(\cat), \El_{\idxsh}(\catI)).
  \end{math}
  One can verify that the map
  \labelcref{eq:lax-univalence-canonical-map} is an equivalence by
  reducing it to the cases when \(\cat\) is locally discrete
  (\cref{lem:lax-univalence-2}) and when \(\cat = \walkingCell_{2}\).
\end{proof}

A dual argument shows the following.

\begin{corollary}
  \label{prop:oplax-univalence}
  Let \(\idxsh\) be an \((\infty, 2)\)-category and let
  \(\cat, \catI : \idxsh^{\opMark(1, 2)} \to \Cat^{\nMark{2}}\) be
  functors. We have an equivalence
  \begin{equation*}
    \Map_{\Core_{\nMark{1}}(\Fun(\idxsh^{\opMark(1, 2)}
      ,\Cat^{\nMark{2}})_{\opLaxMark})}(\cat, \catI)
    \simeq \Map_{\nPrefix{2}\Cat_{/ \idxsh}}(\El_{\idxsh}(\cat), \El_{\idxsh}(\catI))
  \end{equation*}
  natural in \(\idxsh\). \qed
\end{corollary}

Natural transformations correspond to functors preserving cocartesian
morphisms.

\begin{proposition}
  \label{prop:strong-univalence}
  The equivalence in \cref{prop:lax-univalence} is restricted to an
  equivalence between
  \begin{math}
    \Map_{\Core_{\nMark{1}}(\Fun(\idxsh, \Cat^{\nMark{2}}))}(\cat, \catI)
  \end{math}
  and the space of functors
  \(\El_{\idxsh}(\cat) \to \El_{\idxsh}(\catI)\) over \(\idxsh\)
  preserving cocartesian \(1\)-cells.
\end{proposition}
\begin{proof}
  The equivalence between these mapping spaces is due to
  Lurie \parencite[Theorem 3.8.1]{lurie2009goodwillie}. One can see that it
  coincides with the equivalence in \cref{prop:lax-univalence}.
\end{proof}
\noindent 
Oplax limits are \((\infty, 1)\)-categories of oplax natural
transformations in the following sense.

\begin{proposition}
  \label{prop:oplax-limit-universal-property}
  Let \(\cat : \idxsh^{\opMark(1, 2)} \to \Cat^{\nMark{2}}\) be a
  functor. For any \((\infty, 1)\)-category \(\catI\), we have a natural
  equivalence
  \begin{equation*}
    \Map_{\Cat}(\catI, \opLaxLim_{\idx \in \idxsh} \cat_{\idx})
    \simeq \Map_{\Core_{\nMark{1}}(\Fun(\idxsh^{\opMark(1, 2)},
      \Cat^{\nMark{2}})_{\opLaxMark})}(\lambda \blank. \catI, \cat).
  \end{equation*}
\end{proposition}
\begin{proof}
  \begin{EqReasoning}
    \begin{align*}
      & \term{\Map_{\Cat}(\catI, \opLaxLim_{\idx \in \idxsh} \cat_{\idx})} \\
      \simeq & \by{definition} \\
      & \term{\Map_{\nPrefix{2}\Cat_{/ \idxsh}}(\idxsh \times \catI, \El_{\idxsh}(\cat))} \\
      \simeq & \by{\cref{prop:oplax-univalence}} \\
      & \term{\Map_{\Core_{\nMark{1}}(\Fun(\idxsh^{\opMark(1, 2)}, \Cat^{\nMark{2}})_{\opLaxMark})}(\lambda \blank. \catI, \cat).} \qedhere 
    \end{align*}
  \end{EqReasoning}
\end{proof}
\noindent 
A useful source of oplax natural transformations is the following
\emph{mate correspondence} whose special case when \(\cat\) is an
\((\infty, 1)\)-category is shown by Haugseng et al.\@ \parencite[Corollary
F]{haugseng2021mate-arxiv} in a stronger form of an equivalence
between \((\infty, 1)\)-categories of (op)lax natural transformations.

\begin{proposition}
  \label{prop:mate-correspondence}
  Let \(\cat\) be an \((\infty, 2)\)-category and
  \(\fun, \funI : \cat \to \Cat^{\nMark{2}}\) a functor. We have an
  equivalence natural in \(\cat\) between the following space:
  \begin{itemize}
  \item the space of oplax natural transformations \(\trans : \fun \to
    \funI\) such that \(\trans_{\obj}\) is a left adjoint for every
    \(0\)-cell \(\obj \in \cat\);
  \item the space of lax natural transformations \(\transI : \funI \to
    \fun\) such that \(\transI_{\obj}\) is a right adjoint for every
    \(0\)-cell \(\obj \in \cat\).
  \end{itemize}
  Moreover, when an oplax natural transformation \(\trans\)
  corresponds to a lax natural transformation \(\transI\) via this
  equivalence, \(\trans_{\obj} \adj \transI_{\obj}\) for any
  \(\obj \in \cat\).
\end{proposition}
\begin{proof}
  By
  \begin{EqReasoning}
    \begin{align*}
      & \term{\{\walkingCell_{1} \to \Fun(\cat, \Cat^{\nMark{2}})_{\opLaxMark} \mid \text{left adjoint at every \(\obj \in \cat\)}\}} \\
      \simeq & \by{transpose} \\
      & \term{\{\cat \to \Fun(\walkingCell_{1}, \Cat^{\nMark{2}})_{\LaxMark} \mid \text{valued in left adjoints}\}} \\
      \simeq & \by{\cref{lem:lax-univalence-2}} \\
      & \term{\{\cat \to \Cat^{\nMark{2}}_{/ \walkingCell_{1}} \mid \text{valued in \emph{bicartesian fibrations}}\}} \\
      \simeq & \by{\cref{lem:lax-univalence-2}} \\
      & \term{\{\cat \to \Fun(\walkingCell_{1}^{\opMark}, \Cat^{\nMark{2}})_{\opLaxMark} \mid \text{valued in right adjoints}\}} \\
      \simeq & \by{transpose} \\
      & \term{\{\walkingCell_{1}^{\opMark} \to \Fun(\cat, \Cat^{\nMark{2}})_{\LaxMark} \mid \text{right adjoint at every \(\obj \in \cat\)}\}.}
    \end{align*}
  \end{EqReasoning}
  Recall that a \emph{bicartesian fibration} is a functor between
  \((\infty, 1)\)-categories that is is both a cocartesian fibration and a
  cartesian fibration. Because adjunctions are bicartesian fibrations
  over \(\walkingCell_{1}\) \parencite[Definition
  5.2.2.1]{lurie2009higher}, the middle equivalences hold.
\end{proof}

\section{Semantics of mode sketches}
\label{sec:semant-mode-sketch}

We show that models of a mode sketch \(\modesketch\) are equivalent to
diagrams of \(\infty\)-logoses indexed over \(\modesketch\).

\begin{definition}
  Let \(\modesketch\) be a mode sketch. A \emph{model of
    \(\modesketch\)} is an \(\infty\)-logos \(\logos\) equipped with a
  function \(\mode : \modesketch \to \LexAcc(\logos)\) satisfying
  semantic counterparts of
  \cref{axm:mode-sketch-disjoint,axm:mode-sketch-invertible,axm:mode-sketch-top},
  that is:
  \begin{enumerate}[label=\Alph*'., ref=\Alph*']
  \item \label[axiom]{axm:model-disjoint}
    the functor \(\opModality^{\mode(\idx)}_{\mode(\idxI)} :
    \logos_{\mode(\idx)} \to \logos_{\mode(\idxI)}\) is constant at
    \(\objFinal\) for any \(\idxI \not\le \idx\) in \(\modesketch\);
  \item \label[axiom]{axm:model-invertible} for any thin triangle
    \((\idx_{0} < \idx_{1} < \idx_{2})\) in \(\modesketch\), the
    natural transformation
    \(\unitModality^{\mode(\idx_{0});
      \mode(\idx_{2})}_{\mode(\idx_{1})} :
    \opModality^{\mode(\idx_{2})}_{\mode(\idx_{0})} \To
    \opModality^{\mode(\idx_{1})}_{\mode(\idx_{0})}
    \opModality^{\mode(\idx_{2})}_{\mode(\idx_{1})}\) is invertible;
  \item \label[axiom]{axm:model-top} a map \(\map\) in \(\logos\) is
    an equivalence whenever \(\opModality_{\mode(\idx)} \map\) is for
    every \(\idx \in \modesketch\).
  \end{enumerate}
  A \emph{morphism \((\logos, \mode) \to (\logos', \mode')\) of models of \(\modesketch\)} is a morphism of \(\infty\)-logoses \(\fun : \logos \to \logos'\) such that, for every \(\idx \in \modesketch\), there exists a morphism of \(\infty\)-logoses \(\fun_{\idx} : \logos_{\mode(i)} \to \logos'_{\mode'(i)}\) making the following diagram commute.
  \begin{equation*}
    \begin{tikzcd}
      \logos
      \arrow[r, "\fun"]
      \arrow[d, "\opModality_{\mode(\idx)}"'] &
      \logos'
      \arrow[d, "\opModality_{\mode'(\idx)}"] \\
      \logos_{\mode(\idx)}
      \arrow[r, "\fun_{\idx}"', dotted] &
      \logos'_{\mode'(\idx)}
    \end{tikzcd}
  \end{equation*}
  Note that such a morphism \(F_{\idx}\) is unique since \(\opModality_{\mode(\idx)}\) is a localization.
  The models of \(\modesketch\) and their morphisms form an \((\infty, 1)\)-category \(\Model(\modesketch)\) whose cells of dimension \(\ge 2\) are inherited from \(\Logos\).
\end{definition}

\begin{remark}
  \Cref{axm:model-disjoint,axm:model-invertible} are straightforward
  interpretations of
  \cref{axm:mode-sketch-disjoint,axm:mode-sketch-invertible}, but
  \cref{axm:model-top} might look different from
  \cref{axm:mode-sketch-top}. This is because an interpretation of a
  type-theoretic axiom must has the stability under base change, which
  corresponds to the stability under substitution in type theory. A
  naive interpretation of \cref{axm:mode-sketch-top} would be that an
  object \(\sh \in \logos\) is contractible whenever
  \(\opModality_{\mode(\idx)} \sh\) is for every
  \(\idx \in \modesketch\), but this is not stable under base change in
  that it implies nothing about validity in slices
  \(\logos_{/ \shX}\). In contrast,
  \cref{axm:model-disjoint,axm:model-invertible,axm:model-top} are
  stable under base change in the following sense. Every {\acrLAM}
  \(\mode\) in \(\logos\) induces a {\acrLAM} \(\mode_{\shX}\) in each
  slice \(\logos_{/ \shX}\) determined by
  \(\opModality_{\mode_{\shX}} \sh \simeq
  (\unitModality_{\mode})_{\shX}^{\pbMark} \opModality_{\mode} \sh\)
  for \(\sh \in \logos_{/ \shX}\). A function
  \(\mode : \modesketch \to \LexAcc(\logos)\) then induces a function
  \(\mode_{\shX} : \modesketch \to \LexAcc(\logos_{/ \shX})\) for every
  \(\shX \in \logos\) by \(\mode_{\shX}(\idx) =
  \mode(\idx)_{\shX}\). One can verify that if \(\mode\) satisfies
  \cref{axm:model-disjoint,axm:model-invertible,axm:model-top}, then
  so does \(\mode_{\shX}\). In fact, \cref{axm:model-top} is
  equivalent to that the naive interpretation of
  \cref{axm:mode-sketch-top} holds in all the slices
  \(\logos_{/ \shX}\).
\end{remark}

\begin{construction}
  \label{cst:oo-2-cat-from-mode-sketch}
  Let \(\modesketch\) be a mode sketch. We construct an
  \((\infty, 2)\)-category \(\realize{\modesketch}\) as follows. We regard
  the underlying poset \(\idxModesketch_{\modesketch}\) as a
  simplicial set by taking its nerve. The set
  \(\triModesketch_{\modesketch}\) of thin triangles makes
  \(\idxModesketch_{\modesketch}\) a scaled simplicial set. Let
  \(\freeBracket{\idxModesketch_{\modesketch},
    \triModesketch_{\modesketch}}\) denote the
  \((\infty, 2)\)-category presented by it. We set
  \(\realize{\modesketch} = \freeBracket{\idxModesketch_{\modesketch},
    \triModesketch_{\modesketch}}^{\opMark(2)}\).
\end{construction}

\begin{theorem}
  \label{thm:main-theorem}
  For any mode sketch \(\modesketch\), we have an equivalence between
  the following \((\infty, 1)\)-categories:
  \begin{itemize}
  \item the \((\infty, 1)\)-category \(\Model(\modesketch)\) of models of \(\modesketch\);
  \item the \((\infty, 1)\)-category \(\myDiagram(\modesketch) \subset \Core_{\nMark{1}}(\Fun(\realize{\modesketch}^{\opMark(1,2)}, \enlarge \Cat^{\nMark{2}})_{\opLaxMark})\) whose objects are the functors \(\realize{\modesketch}^{\opMark(1, 2)}
    \to \Logos_{\LexAccMark}^{\nMark{2}}\) and morphisms \(\logosI \to \logosI'\) are the oplax natural transformations \(\trans : \logosI \to \logosI'\) whose components \(\trans_{\idx} : \logosI_{\idx} \to \logosI'_{\idx}\) are all morphisms of \(\infty\)-logoses.
  \end{itemize}
  Moreover, when a model \((\logos, \mode)\) of \(\modesketch\)
  corresponds to a functor
  \(\logosI : \realize{\modesketch}^{\opMark(1, 2)} \to
  \Logos_{\LexAccMark}^{\nMark{2}}\), the following hold.
  \begin{enumerate}
  \item
    \begin{math}
      \logos \simeq \opLaxLim_{\idx \in \realize{\modesketch}} \logosI_{\idx}
    \end{math}
  \item \(\logosI_{i} \simeq \logos_{\mode(i)}\) for every
    \(\idx \in \modesketch\).
  \end{enumerate}
\end{theorem}
\noindent 
The rest of this section is devoted to the proof of
\cref{thm:main-theorem}. In \cref{sec:models-mode-sketches}, we give a
construction of a model of \(\modesketch\) from a functor
\(\realize{\modesketch}^{\opMark(1, 2)} \to
\Logos^{\nMark{2}}_{\LexAccMark}\). In \cref{sec:fracture-gluing}, we
give an inverse construction.

\subsection{Models of mode sketches in oplax limits}
\label{sec:models-mode-sketches}

We first show that the oplax limit of a functor
\(\realize{\modesketch}^{\opMark(1, 2)} \to
\Logos_{\LexAccMark}^{\nMark{2}}\) is part of a model of
\(\modesketch\). We fix a functor
\(\logosI : \realize{\modesketch}^{\opMark(1, 2)} \to
\Logos^{\nMark{2}}_{\LexAccMark}\).

\begin{construction}
  \label{cst:prop-in-oplax-limit}
  For a cosieve \(\sieve\) on \(\modesketch\), we define an object
  \(\propCanonicalI(\sieve) \in \opLaxLim_{\idx \in \realize{\modesketch}}
  \logosI_{\idx}\) by
  \begin{equation*}
    \propCanonicalI(\sieve)_{\idx} =
    \left\{
      \begin{array}{ll}
        \objFinal & \text{if \(\idx \in \sieve\)} \\
        \objInitial & \text{otherwise}.
      \end{array}
    \right.
  \end{equation*}
  The other components are uniquely determined by the universal
  properties of initial and final objects. This determines a lattice
  morphism \(\propCanonicalI\) from cosieves on \(\modesketch\) to
  \((-1)\)-truncated objects in \(\opLaxLim_{\idx \in
    \realize{\modesketch}} \logosI_{\idx}\).
\end{construction}

\begin{notation}
  Let \(\sieve \subset \modesketch\) be a subset. We regard \(\sieve\) as a
  mode sketch with the structure inherited from \(\modesketch\). Let
  \begin{equation*}
    \proj_{\sieve} : \opLaxLim_{\idx \in \realize{\modesketch}} \logosI_{\idx}
    \to \opLaxLim_{\idx \in \realize{\sieve}} \logosI_{\idx}
  \end{equation*}
  denote the restriction functor.
\end{notation}

\begin{lemma}
  \label{lem:localization-projection-sieve}
  For any cosieve \(\sieve\) on \(\modesketch\), the restriction functor
  \begin{math}
    \proj_{\modesketch \setminus \sieve}
  \end{math}
  is the closed localization associated to \(\propCanonicalI(\sieve)\).
\end{lemma}
\begin{proof}
  Let \(\opLaxLim_{\idx \in \realize{\modesketch}} \logosI_{\idx} \to \logosII\)
  denote the closed localization associated to
  \(\propCanonicalI(\sieve)\). Recall that \(\logosII\) is the full
  subcategory of \(\opLaxLim_{\idx \in \realize{\modesketch}} \logosI_{\idx}\)
  spanned by those objects \(\sh\) such that
  \(\propCanonicalI(\sieve) \times \sh \simeq \propCanonicalI(\sieve)\). By the
  definition of \(\propCanonicalI(\sieve)\), this condition is
  equivalent to that \(\sh_{\idx} \simeq \objFinal\) for all
  \(\idx \in \sieve\). Then \(\proj_{\modesketch \setminus \sieve}\) induces an
  equivalence
  \(\logosII \simeq \opLaxLim_{\idx \in \realize{\modesketch \setminus \sieve}}
  \logosI_{\idx}\).
\end{proof}

\begin{lemma}
  \label{cor:localization-projection-cosieve}
  For any cosieve \(\sieve\) on \(\modesketch\), the restriction
  functor
  \begin{math}
    \proj_{\sieve}
  \end{math}
  is the open localization associated to \(\propCanonicalI(\sieve)\).
\end{lemma}
\begin{proof}
  An object
  \(\sh \in(\opLaxLim_{\idx \in \realize{\modesketch}} \logosI_{\idx})_{/
    \propCanonicalI(\sieve)}\) must satisfy that
  \(\sh_{\idx} \simeq \objInitial\) for all
  \(\idx \in \modesketch \setminus \sieve\) by the definition of
  \(\propCanonicalI(\sieve)\) and by \cref{prop:logos-init-strict}. Then
  \(\proj_{\sieve}\) induces an equivalence
  \begin{math}
    (\opLaxLim_{\idx \in \realize{\modesketch}} \logosI_{\idx})_{/ \propCanonicalI(\sieve)}
    \simeq (\opLaxLim_{\idx \in \realize{\sieve}} \logosI_{\idx})_{/ \proj_{\sieve}(\propCanonicalI(\sieve))}
    \simeq \opLaxLim_{\idx \in \realize{\sieve}} \logosI_{\idx}.
  \end{math}
\end{proof}

\begin{proposition}
  \label{prop:localization-projection}
  For any \(\idx \in \modesketch\), the projection
  \begin{math}
    \proj_{\idx} : \opLaxLim_{\idxI \in \realize{\modesketch}} \logosI_{\idxI}
    \to \logosI_{\idx}
  \end{math}
  is a localization.
\end{proposition}
\begin{proof}
  \(\proj_{\idx}\) factors as
  \begin{equation*}
    \opLaxLim_{\idxI \in \realize{\modesketch}} \logosI_{\idxI}
    \xrightarrow{\proj_{(\idx \downarrow \modesketch)}}
    \opLaxLim_{\idxI \in \realize{(\idx \downarrow \modesketch)}} \logosI_{\idxI}
    \xrightarrow{\proj_{(\idx \downarrow \modesketch) \setminus \boundary (\idx \downarrow \modesketch)}}
    \logosI_{\idx},
  \end{equation*}
  Thus, it is a composite of localizations by
  \cref{lem:localization-projection-sieve,cor:localization-projection-cosieve}.
\end{proof}

\begin{construction}
  \label{cst:mode-from-diagram}
  For \(\idx \in \modesketch\), we define
  \(\mymode_{\logosI}(\idx) \in \LexAcc(\opLaxLim_{\idx \in
    \realize{\modesketch}} \logosI_{\idx})\) to be the {\acrLAM}
  corresponding to the localization
  \(\proj_{\idx} : \opLaxLim_{\idx \in \realize{\modesketch}}
  \logosI_{\idx} \to \logosI_{\idx}\)
  (\cref{prop:localization-projection}).
\end{construction}

\begin{proposition}
  \label{thm:intended-model-mode-sketch}
  The pair
  \((\opLaxLim_{\idx \in \realize{\modesketch}} \logosI_{\idx},
  \mymode_{\logosI})\) is a model of \(\modesketch\) for any functor
  \(\logosI : \realize{\modesketch}^{\opMark(1, 2)} \to
  \Logos_{\LexAccMark}^{\nMark{2}}\).
\end{proposition}
\noindent 
\Cref{thm:intended-model-mode-sketch} breaks into three parts
(\cref{thm:intended-model-mode-sketch-c,thm:intended-model-mode-sketch-a,thm:intended-model-mode-sketch-b}). \Cref{axm:model-top}
is immediate from the construction.

\begin{proposition}
  \label{thm:intended-model-mode-sketch-c}
  \(\mymode_{\logosI}\) satisfies \cref{axm:model-top}. \qed
\end{proposition}

For \cref{axm:model-disjoint,axm:model-invertible}, we calculate
\(\opModality^{\mymode_{\logosI}(\idx)}_{\mymode_{\logosI}(\idxI)}\)
and
\(\unitModality^{\mymode_{\logosI}(\idxII);
  \mymode_{\logosI}(\idx)}_{\mymode_{\logosI}(\idxI)}\).

\begin{notation}
  For \(\idxI < \idx\) in \(\modesketch\), let \((\idxI < \idx)\)
  denote the associated generating \(1\)-cell in
  \(\realize{\modesketch}\). For \(\idxII < \idxI < \idx\) in
  \(\modesketch\), let \((\idxII < \idxI < \idx)\) denote the
  associated generating \(2\)-cell
  \((\idxI < \idx) \comp (\idxII < \idxI) \To (\idxII < \idx)\) in
  \(\realize{\modesketch}\).
\end{notation}

\begin{lemma}
  \label{lem:oplax-limit-ignore-thin}
  Let \(\modesketch'\) be the mode sketch with the same underlying
  poset as \(\modesketch\) but with no thin triangle. We have an
  equivalence
  \begin{math}
    \opLaxLim_{\idx \in \realize{\modesketch}} \logosI_{\idx}
    \simeq \opLaxLim_{\idx \in \realize{\modesketch'}} \logosI_{\idx}
  \end{math}
\end{lemma}
\begin{proof}
  This is because \(\El_{\realize{\modesketch}}(\logosI) \to
  \realize{\modesketch}\) is locally a right fibration and thus
  locally conservative.
\end{proof}

\begin{lemma}
  \label{lem:realization-mode-sketch-free}
  Suppose that \(\modesketch\) has no thin triangle. Then
  \(\Core_{\nMark{1}}(\realize{\modesketch})\) is freely generated by
  the strict ordering relation, and \((\idxI < \idx)\) is the final
  object in \(\Map_{\realize{\modesketch}}(\idxI, \idx)\) for any
  \(\idxI < \idx\).
\end{lemma}
\begin{proof}
  By \cref{exm:infinity-2-category-from-poset}.
\end{proof}

\begin{lemma}
  \label{lem:localization-projection-morphism-part-1}
  For any \(\idx \in \modesketch\), the right adjoint \(\idx_{\pbMark}\)
  is given by the following formula for \(\sh \in \logosI_{\idx}\) and
  \(\idxI \in \modesketch\).
  \begin{equation*}
    \idx_{\pbMark}(\sh)_{\idxI} \simeq
    \left\{
      \begin{array}{ll}
        \logosI_{(\idxI < \idx)}(\sh) & \text{if \(\idxI < \idx\)} \\
        \sh & \text{if \(\idxI = \idx\)} \\
        \objFinal & \text{otherwise}
      \end{array}
    \right.
  \end{equation*}
\end{lemma}
\begin{proof}
  We first see that we may assume without loss of generality that
  \(\idx\) is the largest element of \(\modesketch\). Otherwise,
  factor \(\proj_{\idx}\) as
  \begin{math}
    \opLaxLim_{\idxI \in \realize{\modesketch}} \logosI_{\idxI}
    \xrightarrow{\proj_{(\modesketch \downarrow \idx)}}
    \opLaxLim_{\idxI \in \realize{(\modesketch \downarrow \idx)}} \logosI_{\idxI}
    \xrightarrow{\proj_{(\idx \downarrow \modesketch)}}
    \logosI_{\idx},
  \end{math}
  where
  \((\modesketch \downarrow \idx) = \{\idxI \in \modesketch \mid \idxI \le \idx\}\). The first
  functor is a closed localization by
  \cref{lem:localization-projection-sieve} because
  \(\modesketch \setminus (\modesketch \downarrow \idx)\) is a cosieve, and the second functor is
  an open localization by
  \cref{cor:localization-projection-cosieve}. The right adjoint of
  \(\proj_{(\modesketch \downarrow \idx)}\) is then defined by extending
  \(\sh \in \opLaxLim_{\idxI \in \realize{(\modesketch \downarrow \idx)}} \logosI_{\idxI}\) by the final
  objects at all \(\idxI \in \modesketch \setminus (\modesketch \downarrow \idx)\). Therefore, the
  problem is reduced to the calculation of the right adjoint of
  \(\proj_{(\idx \downarrow \modesketch)}\), and in this case \(\idx\) is the largest
  element of \((\modesketch \downarrow \idx)\).

  Let \(\idx_{\pbMark}'(\sh)_{\idxI}\) be defined by the displayed
  formula. We turn \(\idx_{\pbMark}'(\sh)\) into an object of
  \(\opLaxLim_{\idxI \in \realize{\modesketch}} \logosI_{\idxI}\). By
  \cref{lem:oplax-limit-ignore-thin}, we assume that \(\modesketch\)
  has no thin triangle. Since
  \(\El_{\realize{\modesketch}}(\logosI) \to \realize{\modesketch}\) is
  locally a right fibration, it follows from
  \cref{lem:realization-mode-sketch-free} that an object
  \(\shI \in \opLaxLim_{\idxI \in \realize{\modesketch}} \logosI_{\idxI}\)
  is completely determined by \(\shI_{\idxI} \in \logosI_{\idxI}\) for
  all \(\idxI \in \modesketch\) and
  \(\shI_{(\idxII < \idxI)} : \shI_{\idxII} \to \logosI_{(\idxII <
    \idxI)}(\shI_{\idxI})\) for all \(\idxII < \idxI\) in
  \(\modesketch\). We can then extend \(\idx_{\pbMark}'(\sh)\) as
  follows.
  \begin{equation*}
    \idx_{\pbMark}'(\sh)_{(\idxII < \idxI)} =
    \left\{
      \begin{array}{ll}
        \logosI_{(\idxII < \idxI < \idx)}(\sh) : \logosI_{(\idxII < \idx)}(\sh) \to \logosI_{(\idxII < \idxI)}(\logosI_{(\idxI < \idx)}(\sh)) & \text{if \(\idxI < \idx\)} \\
        \id & \text{if \(\idxI = \idx\)}
      \end{array}
    \right.
  \end{equation*}

  Since \(\idx_{\pbMark}'(\sh)_{\idx} \simeq \sh\) by construction, we have
  a unique map \(\map : \idx_{\pbMark}'(\sh) \to \idx_{\pbMark}(\sh)\)
  whose \(\idx\)-th component is the identity on \(\sh\). To see that
  \(\map\) is invertible, it suffices to construct a retraction
  \(\mapI\) of \(\map\). Indeed, if \(\mapI \comp \map \simeq \id\), then
  the \(\idx\)-th component of \(\mapI\) must be the identity, and
  thus \(\map \comp \mapI \simeq \id\) follows by adjointness. Let \(\idxI
  \in \modesketch\). If \(\idxI = \idx\), then we must define
  \(\mapI_{\idx} = \id\). Suppose that \(\idxI < \idx\) and consider
  the following commutative diagram.
  \begin{equation*}
    \begin{tikzcd}
      \idx_{\pbMark}'(\sh)_{\idxI}
      \arrow[r, "\map_{\idxI}"]
      \arrow[d, "\idx_{\pbMark}'(\sh)_{(\idxI < \idx)}"', "\simeq"] &
      [6ex]
      \idx_{\pbMark}(\sh)_{\idxI}
      \arrow[d, "\idx_{\pbMark}(\sh)_{(\idxI < \idx)}"] \\
      \logosI_{(\idxI < \idx)}(\idx_{\pbMark}'(\sh)_{\idx})
      \arrow[r, "\logosI_{(\idxI < \idx)}(\map_{\idx})"', "\simeq"] &
      \logosI_{(\idxI < \idx)}(\idx_{\pbMark}(\sh)_{\idx})
    \end{tikzcd}
  \end{equation*}
  The left and bottom maps are invertible by definition. Hence, we
  have a unique retraction \(\mapI_{\idxI}\) of \(\map_{\idxI}\)
  commuting with \(\idx_{\pbMark}(\sh)_{\idxI < \idx}\). This defines
  a retraction of \(\map\).
\end{proof}

\begin{lemma}
  \label{prop:localization-projection-morphism-part}
  For any \(\idxI < \idx\) in \(\modesketch\), the functor
  \begin{math}
    \opModality^{\mymode_{\logosI}(\idx)}_{\mymode_{\logosI}(\idxI)} :
    \logosI_{\idx} \to \logosI_{\idxI}
  \end{math}
  is equivalent to \(\logosI_{(\idxI < \idx)}\).
\end{lemma}
\begin{proof}
  By \cref{lem:localization-projection-morphism-part-1}.
\end{proof}

\begin{proposition}
  \label{thm:intended-model-mode-sketch-a}
  \(\mymode_{\logosI}\) satisfies \cref{axm:model-disjoint}.
\end{proposition}
\begin{proof}
  By \cref{lem:localization-projection-morphism-part-1}.
\end{proof}

\begin{lemma}
  \label{prop:localization-projection-triangle-part}
  For any \(\idxII < \idxI < \idx\) in \(\modesketch\), the natural
  transformation
  \(\unitModality^{\mymode_{\logosI}(\idxII);
    \mymode_{\logosI}(\idx)}_{\mymode_{\logosI}(\idxI)}\) is
  equivalent to \(\logosI_{(\idxII < \idxI < \idx)}\).
\end{lemma}
\begin{proof}
  This is a consequence of
  \cref{lem:localization-projection-morphism-part-1}.
\end{proof}

\begin{proposition}
  \label{thm:intended-model-mode-sketch-b}
  \(\mymode_{\logosI}\) satisfies \cref{axm:model-invertible}.
\end{proposition}
\begin{proof}
  By \cref{prop:localization-projection-triangle-part}.
\end{proof}

\begin{proof}[Proof of \cref{thm:intended-model-mode-sketch}]
  By
  \cref{thm:intended-model-mode-sketch-c,thm:intended-model-mode-sketch-a,thm:intended-model-mode-sketch-b}.
\end{proof}

\begin{construction}
  We extend the construction \(\logosI \mapsto (\opLaxLim_{\idx \in \realize{\modesketch}} \logosI_{\idx}, \mymode_{\logosI})\) to a functor \(\mymode : \myDiagram(\modesketch) \to \Model(\modesketch)\) as follows.
  Characterized by the universal property (\cref{prop:oplax-limit-universal-property}), the oplax limit construction extends to a functor
  \begin{equation*}
    \Core_{\nMark{1}}(\Fun(I^{\opMark(1, 2)}, \Cat^{\nMark{2}})_{\opLaxMark}) \to \Cat^{\nMark{2}}.
  \end{equation*}
  It then restricts to a functor
  \begin{math}
    \myDiagram(\modesketch) \to \Logos
  \end{math}
  by \cref{prop:logos-accessible-oplax-limit}.
  It further lifts to a functor
  \begin{math}
    \myDiagram(\modesketch) \to \Model(\modesketch)
  \end{math}
  by the construction of \(\mymode_{\logosI}\) (\cref{cst:mode-from-diagram}).
\end{construction}

\subsection{Fracture and gluing}
\label{sec:fracture-gluing}

We show that any model of a mode sketch \(\modesketch\) induces a
functor
\(\realize{\modesketch}^{\opMark(1, 2)} \to
\Logos^{\nMark{2}}_{\LexAccMark}\) and that this gives an inverse of
the construction given in \cref{sec:models-mode-sketches}. This is an
externalization and generalization of the \emph{fracture and gluing
  theorem} (\cref{prop:join-strongly-disjoint}).

\begin{construction}
  Let \((\logos, \mode)\) be a model of \(\modesketch\). We define \(\mylogos_{\mode}\) to be
  the full subcategory of
  \(\realize{\modesketch}^{\opMark(1,2)} \times \logos\) spanned by those
  objects \((\idx, \sh)\) such that \(\sh\) belongs to
  \(\logos_{\mode(\idx)}\).
\end{construction}

\begin{proposition}
  \label{prop:mylogos-locally-cocartesian-fibration}
  For any model \((\logos, \mode)\) of \(\modesketch\), the projection
  \(\mylogos_{\mode} \to \realize{\modesketch}^{\opMark(1,2)}\) is a
  \(1\)-cocartesian \(2\)-right fibration.
\end{proposition}
\begin{proof}
  We work with the scaled simplicial sets model.
  It suffices to show that the pullback \(\mylogos'_{\mode}\) of \(\mylogos_{\mode}\) along the fibrant replacement \(\idxModesketch_{\modesketch}^{\opMark} \to \realize{\modesketch}^{\opMark(1, 2)}\) is a locally cocartesian fibration.

  Let
  \((\idx_{0} \le \idx_{1}) : \stdsimp^{1} \to
  \idxModesketch_{\modesketch}^{\opMark}\) be a map which corresponds
  to an ordered pair \((\idx_{0} \le \idx_{1})\) in
  \(\idxModesketch_{\modesketch}^{\opMark}\). A morphism
  \((\idx_{0}, \sh_{0}) \to (\idx_{1}, \sh_{1})\) in
  \(\mylogos'_{\mode}\) over \((\idx_{0} \le \idx_{1})\) is a map
  \(\sh_{0} \to \sh_{1}\) in \(\logos\), but it corresponds to a map
  \(\opModality_{\mode(\idx_{1})} \sh_{0} \to \sh_{1}\) in
  \(\logos_{\mode(\idx_{1})}\). Hence,
  \((\idx_{0} \le \idx_{1})^{\pbMark} \mylogos'_{\mode} \to \stdsimp^{1}\)
  is the Grothendieck construction for the diagram
  \begin{math}
    \logos_{\mode(\idx_{0})}
    \xrightarrow{\opModality^{\mode(\idx_{0})}_{\mode(\idx_{1})}}
    \logos_{\mode(\idx_{1})}
  \end{math}
  and thus a cocartesian fibration.

  Let
  \((\idx_{0} \le \idx_{1} \le \idx_{2}) : \stdsimp^{2} \to
  \idxModesketch_{\modesketch}^{\opMark}\) be a thin \(2\)-simplex. By
  \cref{axm:mode-sketch-invertible}, the canonical natural
  transformation
  \begin{equation}
    \label{eq:canonical-trans-1}
    \begin{tikzcd}
      \logos_{\mode(\idx_{0})}
      \arrow[rr, "\opModality^{\mode(\idx_{0})}_{\mode(\idx_{2})}",
      "\phantom{a}"'{name = a0}]
      \arrow[from = a0, dr, To, end anchor = {[yshift = 1ex]},
      "\simeq"]
      \arrow[dr, "\opModality^{\mode(\idx_{0})}_{\mode(\idx_{1})}"'] & &
      \logos_{\mode(\idx_{2})} \\
      & \logos_{\mode(\idx_{1})}
      \arrow[ur, "\opModality^{\mode(\idx_{1})}_{\mode(\idx_{2})}"']
    \end{tikzcd}
  \end{equation}
  is invertible, and \((\idx_{0} \le \idx_{1} \le \idx_{2})^{\pbMark}
  \mylogos'_{\mode} \to \stdsimp^{2}\) is the Grothendieck construction
  for the diagram \labelcref{eq:canonical-trans-1} and thus a
  cocartesian fibration.
\end{proof}

\begin{construction}
  \label{cst:diagram-from-model-of-modesketch}
  Let \((\logos, \mode)\) be a model of \(\modesketch\). By
  \cref{prop:mylogos-locally-cocartesian-fibration,prop:universal-1-cocart-2-right-fib},
  the projection
  \(\mylogos_{\mode} \to \realize{\modesketch}^{\opMark(1,2)}\) is
  classified by a functor
  \begin{equation*}
    \mylogosI_{\mode} : \realize{\modesketch}^{\opMark(1, 2)} \to \enlarge \Cat^{\nMark{2}}.
  \end{equation*}
  By construction,
  \(\mylogosI_{\mode}(\idx) \simeq \logos_{\mode(\idx)}\). As we have seen
  in the proof of \cref{prop:mylogos-locally-cocartesian-fibration},
  \(\mylogosI_{\mode}\) maps a \(1\)-cell \(\idx_{0} \le \idx_{1}\) in
  \(\realize{\modesketch}\) to
  \begin{math}
    \opModality^{\mode(\idx_{1})}_{\mode(\idx_{0})} :
    \logos_{\mode(\idx_{1})} \to \logos_{\mode(\idx_{0})}
  \end{math}
  which is lex and accessible (but need not preserve all
  colimits). Therefore, \(\mylogosI_{\mode}\) factors through
  \(\Logos_{\LexAccMark}^{\nMark{2}}\).
\end{construction}

We have constructed back and forth constructions between the \((\infty,1)\)-category of models of
\(\modesketch\) and the \((\infty,1)\)-category of functors
\(\realize{\modesketch}^{\opMark(1, 2)} \to
\Logos^{\nMark{2}}_{\LexAccMark}\).
We turn these constructions into an adjunction (\cref{cst:diagram-model-adj}) and then show that its unit and counit are invertible (\cref{thm:fracture-gluing,thm:gluing-fracture}).

\begin{construction}
  \label{cst:diagram-model-adj}
  Let \((\logos, \mode)\) be a model of \(\modesketch\) and \(\logosI : \realize{\modesketch}^{\opMark(1, 2)} \to \Logos^{\nMark{2}}_{\LexAccMark}\) a functor.
  We construct an equivalence natural in \(\logosI \in \myDiagram(\modesketch)\)
  \begin{equation}
    \label{eq:mylogos-universal-property}
    \Map_{\myDiagram(\modesketch)}(\mylogosI_{\mode}, \logosI) \simeq \Map_{\Model(\modesketch)}((\logos, \mode), (\opLaxLim_{\idx \in \realize{\modesketch}} \logosI_{\idx}, \mymode_{\logosI}))
  \end{equation}
  as follows.
  By \cref{prop:mate-correspondence}, a morphism \(\mylogosI_{\mode} \to \logosI\) corresponds to a lax natural transformation \(\logosI \to \mylogosI_{\mode}\) whose components are right adjoints of morphisms of \(\infty\)-logoses.
  It corresponds by \cref{prop:lax-univalence} to a map \(\El_{\realize{\modesketch}^{\opMark(1,2)}}(\logosI) \to \mylogos_{\mode}\) over \(\realize{\modesketch}^{\opMark(1,2)}\) whose fibers are right adjoints of morphisms of \(\infty\)-logoses.
  By the definition of \(\mylogos_{\mode}\), it corresponds to a map \(\El_{\realize{\modesketch}^{\opMark(1,2)}}(\logosI) \to \realize{\modesketch}^{\opMark(1,2)} \times \logos\) over \(\realize{\modesketch}^{\opMark(1,2)}\) whose fiber over \(\idx \in \realize{\modesketch}\) is a right adjoint of a morphism of \(\infty\)-logos that factors through \(\logos_{\mode(\idx)}\).
  Again by \cref{prop:lax-univalence,prop:mate-correspondence}, it corresponds to an oplax natural transformation \((\lambda \blank.\logos) \to \logosI\) whose component at \(\idx \in \realize{\modesketch}\) is a morphism of \(\infty\)-logoses that extends along \(\opModality_{\mode(i)} : \logos \to \logos_{\mode(i)}\).
  By \cref{prop:oplax-limit-universal-property}, it corresponds to a morphism \((\logos, \mode) \to (\opLaxLim_{\idx \in \realize{\modesketch}} \logosI_{\idx}, \mymode_{\logosI})\) in \(\Model(\modesketch)\).
  All of these correspondences are stated in the form of equivalence of spaces natural in \(\logosI \in \myDiagram(\modesketch)\), and thus we obtain \cref{eq:mylogos-universal-property}.

  By \cref{eq:mylogos-universal-property}, the functor \(\mymode : \myDiagram(\modesketch) \to \Model(\modesketch)\) has the left adjoint \((\logos, \mode) \mapsto \mylogosI_{\mode}\).
  Let
  \begin{math}
    \myfun_{\mode} : (\logos, \mode) \to (\opLaxLim_{\idx \in \realize{\modesketch}} \mylogosI_{\mode}(\idx), \mymode_{\mylogosI_{\mode}})
  \end{math}
  and
  \begin{math}
    \mytrans_{\logosI} : \mylogosI_{\mymode_{\logosI}} \to \logosI
  \end{math}
  be the unit and counit, respectively, of the adjunction.
\end{construction}

\begin{lemma}
  \label{lem:gluing-subtopos}
  Let \(\logos\) be an \(\infty\)-logos and let \(\mode\) and
  \(\modeI\) be {\acrLAMs} in \(\logos\). Suppose that
  \(\opModality^{\modeI}_{\mode}\) is constant at \(\objFinal\). Then
  the functor \(\logos \to \Glue(\opModality^{\mode}_{\modeI})\) that
  sends \(\sh \in \logos\) to
  \((\opModality_{\modeI} \sh, \opModality_{\mode} \sh,
  \opModality_{\modeI} (\unitModality_{\mode})_{\sh}) \in
  \Glue(\opModality^{\mode}_{\modeI})\) is a localization.
\end{lemma}
\begin{proof}
  Observe that the right adjoint of the functor
  \(\logos \to \Glue(\opModality^{\mode}_{\modeI})\) sends
  \((\shI', \shI, \mapI) \in \Glue(\opModality^{\mode}_{\modeI})\) to
  the pullback
  \begin{equation*}
    \begin{tikzcd}
      \unit_{\modeI}^{\pbMark} \shI'
      \arrow[r]
      \arrow[d]
      \arrow[dr, pbMark]
      & \shI'
      \arrow[d, "\mapI"] \\
      \shI
      \arrow[r, "\unit_{\modeI}"']
      & \opModality^{\mode}_{\modeI} \shI.
    \end{tikzcd}
  \end{equation*}
  \(\opModality_{\modeI}\) inverts \(\unitModality_{\modeI}\), and
  thus
  \(\opModality_{\modeI} \unitModality_{\modeI}^{\pbMark} \shI' \simeq
  \shI'\). Since \(\opModality^{\modeI}_{\mode}\) is constant at
  \(\objFinal\), it sends \(\mapI\) to the identity on \(\objFinal\),
  and thus
  \(\opModality_{\mode} \unitModality_{\modeI}^{\pbMark} \shI' \simeq
  \shI\). Therefore, the counit for this adjunction is invertible.
\end{proof}

\begin{proposition}
  \label{thm:fracture-gluing}
  The unit
  \begin{math}
    \myfun = \myfun_{\mode} : (\logos, \mode) \to (\opLaxLim_{\idx \in \realize{\modesketch}} \mylogosI_{\mode}(\idx), \mymode_{\mylogosI_{\mode}})
  \end{math}
  is an equivalence in \(\Model(\modesketch)\) for any model \((\logos, \mode)\) of \(\modesketch\).
\end{proposition}
\begin{proof}
  Since \(\logos_{\mode(\idx)} \simeq \mylogosI_{\mode}(\idx)\), it remains to show that the underlying functor of \(\myfun\) is an equivalence.
  We show that, for any cosieve \(\sieve\) on \(\modesketch\), the
  composite
  \begin{equation*}
    \myfun_{\sieve} : \logos
    \xrightarrow{\myfun}
    \opLaxLim_{\idx \in \realize{\modesketch}} \mylogosI_{\mode}(\idx)
    \xrightarrow{\proj_{\sieve}}
    \opLaxLim_{\idx \in \realize{\sieve}} \mylogosI_{\mode}(\idx)
  \end{equation*}
  is a localization. In particular, \(\myfun\) itself is a
  localization. Then, \(\myfun\) is an equivalence because it is
  conservative by \cref{axm:model-top}. We proceed by induction on the
  size of \(\sieve\). The case when \(\sieve\) is empty is
  trivial. Suppose that \(\sieve\) is inhabited. There is an element
  \(\idx_{0} \in \sieve\) minimal in \(\sieve\). By induction
  hypothesis, \(\myfun_{\sieve \setminus \{\idx_{0}\}}\) is a localization,
  and let \(\modeI\) be the corresponding {\acrLAM}. By
  \cref{axm:model-disjoint}, it follows that
  \(\opModality^{\modeI}_{\mode(\idx_{0})}\) is constant at
  \(\objFinal\). By \cref{lem:gluing-subtopos},
  \(\Glue(\opModality^{\modeI}_{\mode(\idx_{0})})\) is a localization
  of \(\logos\). Again by \cref{lem:gluing-subtopos}, we have the
  localization
  \(\opLaxLim_{\idx \in \realize{\sieve}} \mylogosI_{\mode}(\idx) \to
  \Glue(\opModality^{\modeI}_{\mode(\idx_{0})})\), but this is also
  conservative and thus an equivalence. Therefore, \(\myfun_{\sieve}\)
  is a localization.
\end{proof}

\begin{proposition}
  \label{thm:gluing-fracture}
  For any functor
  \(\logosI : \realize{\modesketch}^{\opMark(1, 2)} \to
  \Logos^{\nMark{2}}_{\LexAccMark}\), the unit \(\mytrans = \mytrans_{\logosI} : \mylogosI_{\mymode_{\logosI}} \to \logosI\) is an equivalence in \(\myDiagram(\modesketch)\).
\end{proposition}
\begin{proof}
  It suffices to show that \(\sigma\) is an invertible natural transformation.
  It suffices to show that the corresponding lax natural transformation \(\mytrans' : \logosI \to \mylogosI_{\mymode_{\logosI}}\) by \cref{prop:mate-correspondence} is an invertible natural transformation.
  To see that \(\mytrans'\) is a natural transformation, by
  \cref{prop:strong-univalence}, it suffices to check that the
  corresponding map
  \(\fun : \El_{\realize{\modesketch}^{\opMark(1,2)}}(\logosI) \to
  \mylogos_{\mymode_{\logosI}}\) over \(\realize{\modesketch}^{\opMark(1,2)}\) preserves cocartesian morphisms.
  Unfolding the definition (\cref{cst:diagram-model-adj}), \(\fun : \El_{\realize{\modesketch}^{\opMark(1,2)}}(\logosI) \to
  \mylogos_{\mymode_{\logosI}} \subset \realize{\modesketch}^{\opMark(1,2)} \times \opLaxLim_{\idxI \in \realize{\modesketch}} \logosI_{\idxI}\) sends an object \((\idx, \sh)\) to \((\idx, \sh) \in \realize{\modesketch}^{\opMark(1,2)} \times \logosI_{\idx} \subset \realize{\modesketch}^{\opMark(1,2)} \times \opLaxLim_{\idxI \in \realize{\modesketch}} \logosI_{\idxI}\) and a morphism \((\idx \ge \idx', \map) : (\idx, \sh) \to (\idx', \sh')\) in \(\El_{\realize{\modesketch}^{\opMark(1,2)}}(\logosI)\) to \((\idx \ge \idx', \map') : (\idx, \sh) \to (\idx', \sh')\) in \(\realize{\modesketch}^{\opMark(1,2)} \times \opLaxLim_{\idxI \in \realize{\modesketch}} \logosI_{\idxI}\), where \(\map'\) is the composite
  \begin{math}
    \sh \xrightarrow{\unitModality_{\mymode_{\logosI}(\idx')}} \opModality_{\mymode_{\logosI}(\idx')} \sh \simeq \logosI_{(\idx' \le \idx)}(\sh) \xrightarrow{\map} \sh'.
  \end{math}
  When \((\idx \ge \idx', \map)\) is cocartesian, \(f\) is invertible, and then \((\idx \ge \idx', \map')\) is a cocartesian morphism in \(\mylogos_{\mymode_{\logosI}}\).
  By construction,
  \(\mytrans'\) is point-wise invertible and thus invertible by
  \cref{prop:invertible-trans-pointwise}.
\end{proof}

\begin{proof}[Proof of \cref{thm:main-theorem}]
  The functor \(\mymode : \myDiagram(\modesketch) \to \Model(\modesketch)\) gives an equivalence by \cref{cst:diagram-model-adj,thm:fracture-gluing,thm:gluing-fracture}.
\end{proof}

\section*{Acknowledgements}
\label{sec:acknowledgements}

The author thanks Jonathan Sterling for useful conversations on the current work.
The author was supported by KAW Grant ``Type Theory for Mathematics and Computer Science'' investigated by Thierry Coquand and Peter LeFanu Lumsdaine.

\bibliographystyle{alphaurl}
\bibliography{my-references}

\end{document}